# A Taxonomy of C-systems °•

WALTER A. CARNIELLI   CLE and IFCH, Unicamp, Brazil
carniell@cle.unicamp.br

JOÃO MARCOS   RUG, Ghent, Belgium, and IFCH, Unicamp, Brazil
vegetal@cle.unicamp.br

**Abstract**
The logics of formal inconsistency (**LFI**s) are paraconsistent logics which permit us to internalize the concepts of consistency or inconsistency inside our object language, introducing new operators to talk about them, and allowing us, in principle, to logically separate the notions of contradictoriness and of inconsistency. We present the formal definitions of these logics in the context of General Abstract Logics, argue that they in fact represent the majority of all paraconsistent logics existing up to this point, if not the most exceptional ones, and we single out a subclass of them called **C**-systems, as the **LFI**s that are built over the positive basis of some given consistent logic. Given precise characterizations of some received logical principles, we point out that the gist of paraconsistent logic lies in the Principle of Explosion, rather than in the Principle of Non-Contradiction, and we also sharply distinguish these two from the Principle of Non-Triviality, considering the next various weaker formulations of explosion, and investigating their interrelations. Subsequently, we present the syntactical formulations of some of the main **C**-systems based on classical logic, showing how several well-known logics in the literature can be recast as such a kind of **C**-systems, and carefully study their properties and shortcomings, showing for instance how they can be used to faithfully reproduce all classical inferences, despite being themselves only fragments of classical logic, and venturing some comments on their algebraic counterparts. We also define a particular subclass of the **C**-systems, the **dC**-systems, as the ones in which the new operators of consistency and inconsistency can be dispensed. A survey of some general methods adequate to provide these logics with suitable interpretations, both in terms of valuation semantics and of possible-translations semantics, is to be found in a follow-up, the paper [42]. This study is intended both to fully present and characterize, from scratch, the field into which it inserts, hinting of course to the connections with other studies by several authors, as well as to set some open problems, and to point to a few directions of continuation, establishing on the way a unifying theoretical framework for further investigation for researchers involved with the foundations of paraconsistent logic.

---

°• Carnielli acknowledges financial support from CNPq / Brazil and from the A. von Humboldt Foundation, and thanks colleagues from the Advanced Reasoning Forum present at the Bucharest meeting in 2000 for the opportunity of discussing some aspects of this work. Those discussions gave rise to the pamphlet [40], an embryonic version of the present paper. Marcos acknowledges, first, the financial support received from CNPq / Brazil, and, later, from a Dehousse doctoral grant in Ghent, Belgium. Both authors acknowledge support also from a CAPES / DAAD grant for a ProBrAl project Campinas / Karlsruhe, and are indebted to all the colleagues present at the II World Congress on Paraconsistency, and especially to Newton da Costa, for his achievements, ideas, and his enthusiasm, both enduring and contagious. And, of course, to the excellent comments by Chris Mortensen, Dirk Batens and Jean-Yves Béziau on a beta version of this paper, as well as to all the patient people who have waited long enough for this study to be concluded. As we already clarified above, a complete semantical study of the systems here presented is soon to be found in [42] (as an outcome of [76]). Dividing it into two papers, we have tried to keep the length and termination of this study a bit more reasonable. All comments are welcome in the meanwhile.

# 1 THOU SHALT NOT TRIVIALIZE!

> On account of the classical principle of [non-]contradiction, a proposition and its negation cannot be both simultaneously true; thanks to this, it is not possible that a theory which is valid under the philosophical (or logical) point of view includes internal contradictions. To suppose the contrary would seemingly constitute a philosophical error.
> —Newton C. A. da Costa, [46], p. 6–7, 1958.

In the dawn of the XXI century, debates on the statute of contradiction in logic, philosophy and mathematics are still likely to raise the most diverse and animated sentiments. And this is an old story, whose first dramatic strokes can be traced back to authors as early as Aristotle (for the defense of non-contradiction), or Heraclitus (for the contrary position). Be that as it may, the fact is that in the beginning of the last century essentially the same dispute was still taking place, this time contraposing Russell to Meinong. And so it could still proceed, for centuries, if only the philosophical aspects of the dispute were touched. Even on more technical grounds, logicians of caliber, such as Alfred Tarski, would eventually speculate about that (cf. [106]):

> I do not think that our attitude towards an inconsistent theory would change even if we decided for some reason to weaken our system of logic so as to deprive ourselves of the possibility of deriving every sentence from any two contradictory sentences. It seems to me that the real reason of our attitude is a different one: We know (if only intuitively) that an inconsistent theory must contain false sentences; and we are not inclined to regard as acceptable any theory which has been shown to contain such sentences.

Against such suspicions, the philosopher Wittgenstein, who had devoted almost half of his late work to the philosophy of mathematics and used to refer to it as his 'main contribution' (cf. the entry *Mathematics*, in [63]), would have had something to say. Indeed, he often felt puzzled about 'the superstitious fear and awe of mathematicians in face of contradiction' (cf. [109], Ap.III–17), and asked himself: 'Contradiction. Why just this *one* spectre? This is surely much suspect.' (id., IV–56). His point was that 'it is one thing to use a mathematical technique consisting in the avoidance of contradiction, and another thing to philosophize against contradiction in mathematics' (id., IV–55), and that it was necessary to remove the 'metaphysical thorn' stuck here (id., VII–12). In this respect, the philosopher described his own objective as precisely that of altering the *attitude* of mathematicians concerning contradictions (id., III–82).

The above passage from da Costa's [46] could also be directed upon criticizing a position such as the above one of Tarski. The presupposition to be challenged here, of course, is that of an inconsistent theory obligatorily containing false sentences. Thus, if models may be described of structures in which some (but not all) contradictory sentences are simultaneously true, we will have a technical point against such suspicions of impossibility or implausibility of maintaining contradictory sentences inside of some theory and still being able to perform reasonable inferences from that, *instead* of being able to derive arbitrarily other sentences. This is sure to make a point in conferring to the task of studying the behavior of contradictory yet non-trivial theories —the task of *paraconsistency*— some respectability. And it *is* indeed possible to assign models for inconsistent non-trivial theories, even if these were to be regarded by some as epistemologically puzzling, or ontologically perplexing! Obtaining models and understanding their role is certainly an extraordinarily important mathematical enterprise: Enourmous efforts from the most brilliant minds and more than twenty centuries were required until mathematicians would allow themselves to consider models in which, given a straight line $S$ and a point $P$ outside of it, one could draw not just one line, but infinite, or no parallel lines to $S$ passing through $P$, as in the well-known case of non-Euclidean geometries. In the present case, then, the problem will not be that of *validating falsities*, but that of *extending our notion of truth* (an idea further explored, for instance, in [28]).

At that same decisive moment, in the first half of the last century, there were in fact these other people like Łukasiewicz or Vasiliev who were soon proposing relativizations of the idea of non-contradiction, offering formal interpretations to formal systems in which this idea did not hold, and in which contradictions could make sense. And in between the 40s and the 60s the world would finally be watching the birth of the first real operative systems of paraconsistent logic (cf. Jaśkowski's [67], Nelson's [86], and da Costa's [49]). But the paleontology of paraconsistent logic will not be our main subject here —for that we prefer to redirect the reader to some of the following articles [6], [59], [55], and those in section 1 of [95], plus the book [26].

**1.1 Contradictory theories do exist.** Be them a consequence of the only correct description of a contradictory world (as assumed in [90]), be them just a temporary state of our knowledge, or again the outcome of a particular language that we have chosen to describe the world, the result of conflicting observational criteria, superpositions of worldviews, or simply, in science, because they result from the best theories available at a given moment (cf. [14]), contradictions are presumably unavoidable in our theories. Even if contradictory theories were to appear only by mistake, or perhaps by some Janus-like crooked behavior of their proposers, it is hard to see, given for instance results such as Gödel's incompleteness theorems, how contradictions could be prevented from even being taken into consideration. So it should be clear that the point here is not about the *existence* of contradictory theories, but about *what we should do* with them! Should these theories be allowed to explode and derive anything else, as in classical logic, or rather should we try to substitute the underlying logic, in (potentially) critical situations, in order to still be able to draw (if only temporarily, if you want) reasonable conclusions from those theories?

At this point it is interesting to consider the following motto set down by Newton da Costa, one of the founders of modern paraconsistent logic (cf. [47]):

> From the syntactical-semantical standpoint, every mathematical theory is admissible, unless it is trivial.

Da Costa designated that motto 'Principle of Tolerance in Mathematics', in analogy to the 'syntactical' principle proposed before by Carnap (cf. [35], p.52). According to this, the dividing line in between systems worthy of investigation and those that do not 'make a difference' (cf. [55]), nor convey any information (cf. [14]), should be drawn around non-triviality, rather than in the vicinity of non-contradictoriness. This will give us the first key to paraconsistency: if there are no contradictions around, then everything is under control, once we are inside of a consistent environment; but if contradictions are allowed, non-triviality should be the aim —but then what we must control is the *explosive* character of our underlying logic. Indeed, inside of a consistent logic we know that contradictions are dangerous in a theory precisely because they will give sufficient reason for this theory to explode, deducing anything else!

So, given a logic whose language includes a negation symbol ¬, let's call *contra-dictory* a theory from which some formula $A$ and its negation ¬$A$ can be derived by way of the underlying logic. Let's also call a theory *trivial* if any formula $B$ can be derived from it by way of the underlying logic, and call a theory *explosive* if the addition to it of any contradiction $A$ and ¬$A$ is sufficient to make it trivial. The underlying logic, in its own right, will also be called *contradictory*, *trivial*, or *explosive* if, respectively, all of the theories about which it can talk are contradictory, trivial, or explosive. To be sure, any trivial theory / logic will also turn to be contradictory, whenever there is a negation available (anything is derived from it, in particular all pairs of formulas of the form $A$ and ¬$A$). Inside classical or intuitionistic logic, and, in a general way, inside any 'consistent' logic (this will be defined in what follows), the contradictory and the trivial theories simply coincide, by way of their explosive character. *Paraconsistent logics* were then proposed to be the logics to underlie those contradictory theories which were still to be kept non-trivial, and what those logics must of course effect to such an end is weakening or annulling the explosive character of these theories.[1] So, all at once, paraconsistency comes and provides a sharp distinction in between the logical notions of contradictoriness, explosiveness, and triviality.

Anyone working as a knowledge engineer, assembling and managing knowledge databases, will be perfectly aware that gathering inconsistent information is the rule rather than the exception. And again, either if you assume, by some sort of methodological requirement, inconsistent theories to be problematic (cf. [14]) or not (cf. [90]), this does not prevent you from also assuming them to be, in general, quite *informative*, and wanting to *reason* from them in a sensible way. Consider, for instance, this very simple situation (cf. [40]) in which you ask two people, in the due course of an investigation, a 'yes-no' question such as 'Does Dick live in Arizona?', so that what will result will be exacly one of the three following different possible scenarios: they might both say 'yes', they might both say 'no', or else one of them might say 'yes' while the other says 'no'. Now, it happens that in neither situation you may be sure about where Dick really lives (unless you trust some of the interviewees more than the other), but only in the last scenario, where an inconsistency appears, *are you sure* to have received wrong information from one of your sources!

Our next point is that also the logical notions of *inconsistency* and of *contradictoriness* can and should be *distinguished* in a purely abstract way. Distinctions have already been proposed, in the literature, among the notions of *paradoxical* and of *antinomical* theories (cf., for instance, Arruda's [6], p.3, or da Costa's [51], p.194), the paradoxical ones being identified with those theories in which inconsistencies could occur without necessarily leading to trivialization, and the antinomical ones identified with those in which any occurring contradiction turns out to be fatal, as in the case of Russell's antinomy in naive set theory. Let us, here, insist on this distinction for a moment and stretch it a bit further. One first difficulty to be confronted with is that of some English technical terms: It is such a pity that techniques and results such as Hilbert's witch-hunt programme in search of a *Widerspruchfreiheitbeweis* for Arith-

---

[1] Surprising as it may seem, this would also have been the advice given by Wittgenstein on how to proceed in the presence of contradictions: 'The contradiction does not even falsify anything. Let it lie. Do not go there.' (cf. [110], XIV, p.138) For the relations and non-relations between Wittgenstein and the paraconsistent enterprise the reader may consult, for instance, [75], [64] or [77].

metic were to be eventually translated into the search for a 'consistency proof', given that what it literally means is something much more precise, namely a 'proof of freedom from contradictions'! More often than not, German language indeed shows itself to be exceedingly precise, so that we should rather stick here to the literal meaning of *Widerspruchfreiheit* as non-contradictoriness, and associate inconsistency, if we may, with something like the term *Unbeständigkeit* (or any other synonym of *Inkonsistenz* together perhaps with some terms opposed to *Beschaffenheit* and to *Widerstandsfähigkeit*). Now, antinomies will be related to the presence of 'strong' contradictions —those with explosive behavior—, while paradoxes will be related to the presence of inconsistencies, which do not necessarily depend on negation, such as in the case of the well-known Curry's paradox (cf. [45]). Let us try to summarize this whole story in a picture (maybe you do not agree on our choice of names, but we beg you to stick to our terminology for the moment):

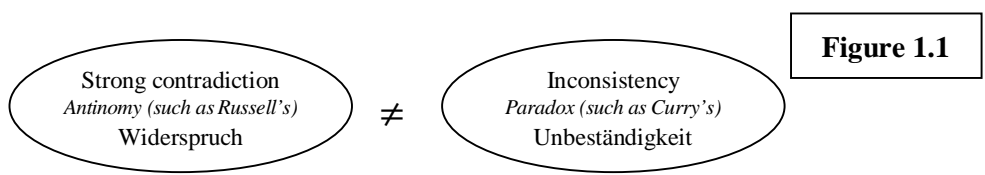

**Figure 1.1**

The above distinctions, of course, are more illustrative than formal (nothing prevents you, for instance, from thinking of Russell's antinomy as something not as destructive as it was, if you just change the underlying logic of its theory so as to make it only paradoxical, but the distinction between antinomies as involving the notion of strong negation, on the one hand, and inconsistencies, on the other hand, as something more general and in principle independent of negation, should be taken more seriously).

Now, whatever an inconsistency might *mean*, be it more general or not than a simple contradiction, we may certainly presuppose that a contradiction is at least an *example* of an inconsistency, be it the only possible one or not. Traditionally, as we have noted a few paragraphs above, the contradictoriness of a given theory / logic was to be identified with the fact that it derives at least some pairs of formulas of the form $A$ and $\neg A$, while inconsistency was usually talked about as a model-theoretic property to be guaranteed so that our theories can make sense and talk about 'real existing structures'. Of course, any trivial theory / logic, thus, given our assumption above that contradictions entail inconsistencies, will also be both contradictory and inconsistent. Now, if explosiveness does not hold, as we shall see, in the scope of paraconsistent logics, there is in principle no reason to suppose that the converse would also be the case, and that a contradiction would always lead to trivialization. How to reconcile these concepts then? Da Costa's idea, when proposing his first paraconsistent calculi (cf. [49]), was that the 'consistency' and the 'classic-like behavior' (he called that 'well-behavior') of a given formula, as a sufficient requisite to guarantee its explosive character, could be represented as simply another formula of its underlying logic (he chose, for his first calculus, $C_1$, to represent the consistency of a formula $A$ by the formula $\neg(A \wedge \neg A)$, and referred to this last formula —to be intuitively read as saying that 'it is not the case that both $A$ and $\neg A$ are true'—, as a realization of the 'Principle of Non-Contradiction', conventions that we will here, in general, *not* follow —neither will we follow, necessarily, the identification of con-sistency with 'classic-like behavior'). In fact, our proposal here, inspired by da Costa's idea, is exactly that of introducing consistency as a *primitive notion* of our logics: the paraconsistent logics which internalize the notion of consistency so as to introduce it already at the object level will be called *logics of formal inconsistency* (**LFI**s). And, given a consistent logic **L**, the **LFI**s which extend the positive basis of **L** will be said to constitute **C**-*systems based on* **L**. Our main aim in this paper, besides making all the above definitions and their multiple shades and interrelations entirely clear, will be that of studying a large class of **C**-systems based on classical logic (of which the calculi $C_n$ of da Costa will be but very particular examples).

On what concerns this story about regarding consistency as a primitive notion, the status of points, lines and planes in geometry may immediately be thought of, but the case of (imaginary) complex numbers seems to make an even better comparison: even if we do not know what they are, and may even suspect there is little sense in insisting on which way they can exist in the 'real' world, the most important aspect is that it is possible to calculate with them. Girolamo Cardano, who first had the idea of computing with such numbers, seems to have seen this point clearly —he failed, however, to acknowledge the importance of this; in 1545 he wrote in his *Ars Magna* (cf. [87]):

> Dismissing mental tortures, and multiplying $5+\sqrt{-15}$ by $5-\sqrt{-15}$, we obtain $25-(-15)$. Therefore the product is 40. …and thus far does arithmetical subtlety go, of which this, the extreme, is, as we have said, so subtle that it is useless.

His discovery, that one could operate with a mathematical concept independent of what our intuition would say and that usefulness (or something else) could be a guiding criterion for accepting or rejecting experimentation with mathematical objects, definitely contributed to the proof of the Fundamental Theorem of Algebra by C. F. Gauss in 1799, before which complex numbers were not fully accepted.

To make matters clear, the basic idea behind the internalization of consistency inside our logics will be, in general, accomplished by the addition of a unary connective expressing consistency (and usually also another connective to express inconsistency), plus the following important assumption, that *consistency* is exactly what a theory might be lacking in order to deliver triviality when exposed to a contradiction.[2] Recapitulating: as we said before, triviality entails contradictoriness (if a negation is present), and contradictoriness entails inconsistency (or, to be more precise, contradictoriness entails 'non-consistency', for it may happen, as we will see, that consistency and inconsistency are not exactly dual in some of our logics, if we take both notions as primitive); now we just add to this the assumption that contradictoriness *plus* consistency implies triviality! We are in fact introducing, in this way, a novel definition of consistency, more fine-grained than the usual model-theoretic one: for a large class of logics (see FACT 2.14(ii)) it will turn out that consistency may be identified with the presence of *both* non-contradictoriness and explosive features.

---

[2] It is interesting to notice, by the way, that this assumption is remarkably compatible with Jaśkowski's intuition on the matter. As he put it, 'in some cases we have to do with a system of hypotheses which, *if subjected to a too consistent analysis*, would result in a contradiction between themselves or with a certain accepted law, but which we use in a way that is restricted so as not to yield a self-evident falsehood' (our italics, see [67], p.144). It is clear that we can give this at least one reading according to which Jaśkowski seemed already to have been worried about the effects of consistent contradictions!

Now, non-contradictoriness will be a necessary but *no more* a sufficient requirement for us to prove consistency. In the case of explosive logics, of course, the concepts of non-contradictoriness and non-triviality will coincide, so that non-contradictoriness and consistency are also to be identified. Paraconsistent logics are situated exactly in that terra incognita which lays in between non-explosive logics and trivial ones, and they comprehend exactly those logics which are both non-explosive and non-trivial (examples of such logics are provided by the whole of the literature on paraconsistent logics)! So, again, consistency divides the logical space in between consistent (and so, explosive and non-contradictory) logics, and inconsistent ones, and these last ones may, at their turn, be either paraconsistent (and so, non-explosive, and possibly even contradictory), or trivial.

**1.2 Paraconsistent, but not contradictory!** In fact, there is another point that we want to stress here, for it seems that much confusion has been unnecessarily raised around it. In general, paraconsistent logics do *not* validate contradictions or invalidate anything like the 'Principle of Non-Contradiction' (though there are a few that do). Most paraconsistent logics, actually, are just fragments of some other given consistent logic (such as some version of classical logic, or else some normal modal logic), so that they *cannot*, in any case, be contradictory! However, a good way of making this whole point much less ambiguous (even though still open to dispute, but now on a different level) is by considering formal definitions of those so-called (meta)logical *principles*.

Let us say that a logic respects the *Principle of Non-Contradiction*, (PNC), if it is non-contradictory, according to our previous definitions, that is, if it has non-contradictory theories, that is, theories in which no contradictory pair of formulas $A$ and $\neg A$ may be inferred. Let us also say that a logic respects the *Principle of Non-Triviality*, (PNT), (a realization of da Costa's Principle of Tolerance inside of the logical space) if it is non-trivial, thus possessing non-trivial theories, and say that a logic respects the *Principle of Explosion*, or *Pseudo-Scotus*, (PPS), if it is explosive, that is, if all of its theories explode when in contact with a contradiction. It is clear now that all paraconsistent logics, by their very nature, must disrespect (PPS), aiming to retain (PNT), but it is also clear that they cannot disrespect (PNC) as long as they are defined as fragments of other logics that do respect (PPS)! The gist and legacy of paraconsistent logic indeed lies in showing that logics may be constructed in which the Principle of Pseudo-Scotus is controlled in its power, and this has 'in principle' *nothing to do* with the validity or not of the Principle of Non-Contradiction as we understand it. Yet a few logics exist which are not only paraconsistent, but that in fact disrespect (PNC). Such logics are usually put forward in order to formalize some dialectical principles, and are accordingly known as *dialectical logics*. Being able to infer contradictions, however, such dialectical logics cannot be fragments of any consistent logic, and in order to avoid trivialization they should also usually assume, for instance, the failure of Uniform Substitution, at least when applied to some specific formulas, such as the contradictions that those logics can infer (or else any other contradiction, and thus any other formula, would be inferable). Much weaker versions of the Principle of Non-Contradiction have nevertheless been considered in the literature, as for instance the following one, deriving from semantical approaches to the matter: a logic is said to respect the *Principle of Non-Contradiction, second* *form*, (PNC2), if it has non-trivial models for pairs of contradictory formulas. But then, of course, every model for the falsification of (PPS), that is, every model for a paraconsistent logic, would also satisfy (PNC2), and vice-versa, so that not only would (PNC2) be unnecessary as a new principle, but there would also be no principle dealing specifically with the existence of dialectical logics. Too bad! And, of course, there is a BIG difference in having models for some specific contradictions and having *all* models of a given logic validating some contradictory pair of formulas — this amounts, in the end, to the same difference which exists, in classical logic, in between contingent formulas, on the one hand, and (tautological or) contradictory ones, on the other hand…

The above definition of (PNC) will also prevent us from identifying this principle, inside some arbitrary given logic **L**, with the validity in **L** of some particular formula, such as $\neg(A \wedge \neg A)$ (as in da Costa's first requisite for the construction of his paraconsistent calculi —check the subsection **3.8**). But it *is* true that such a formula can, as well as many other formulas, be identified, in some situations, to the expression of consistency inside of some specific logics, such as da Costa's $C_1$! Let us, in general, say that a theory $\Gamma$ is *gently explosive* when there is always a way of expressing the consistency of a given formula $A$ by way of formulas which depend only on $A$, that is, when there is a (set of) formula(s) constructed using $A$ as their sole variable and that cannot be added to $\Gamma$ together with a contradiction $A$ and $\neg A$, unless this leads to triviality. A gently explosive logic, then, is exactly a logic having only gently explosive theories, and we can now formulate (gPPS), a 'gentle version' of (PPS), for a given logic **L**, asserting that this logic must be gently explosive. Gently explosive paraconsistent logics, thus, are precisely those logics that we have above dubbed **LFI**s, the logics of formal inconsistency. In the logics we will be studying in this paper, we will in general assume that the consistency of each formula $A$ can be expressed by operators already at their linguistic level, and in the simplest case this will be written as $\circ A$, where '$\circ$' is the 'consistency connective'. The **C**-systems (in this paper they will be supposed to be based on classical logic), will be particular **LFI**s illustrating some different ways in which one can go on to axiomatize the behavior of this new connective.

There are also some other forms of explosion, as the *partial* one, which does not trivialize the whole logic, but just part of it (for instance, when a contradiction does not prove every other formula, but does prove every other *negated* formula). We will let our paraconsistent logics also reject this kind of explosion. There is *ex falso*, which asserts that at least one element should exist in our logics so that everything follows from it (a kind of *falsum*, or *bottom particle*). There is *controllable* explosion, which states that, if not all, at least some of our formulas should lead to trivialization when taken together with their negations. And, finally, there is *supplementing* explosion, which states that our logics should possess, or be able to define, a *supplementing*, or *strong* negation, to the effect that strongly negated propositions (that we have above called strong contradictions) should explode. (There are also all sorts of combinations of these forms of explosions, and perhaps some other forms still to be uncovered, but these are the ones we will concentrate on, here.) All of these alternative forms of explosion can be turned into logical (meta)principles, and none of these rejects, by their own right, 'full' Pseudo-Scotus —all of them, nonetheless, can still be held even when the Pseudo-Scotus does not hold! The para-consistent logics studied in this paper will, of course, disrespect Pseudo-Scotus, and in addition to that they will also disrespect the principle regarding partial explosion, while, in most cases, they will still respect the principles regarding gentle explosion, *ex falso*, supplementing explosion, and, often, controllable explosion as well. This will be made much clearer in section **2**, where this study will be made more precise, and the interrelations between all of those principles will be more deeply investigated.

**1.3 What do you mean?** Let's now briefly describe the exciting things that await the reader in the next sections (we will skip section **1** in our description —you are reading it—, but do not stop here!).

Section **2** is *General Abstract Nonsense*. No particular systems of paraconsistent logic are studied here (though some are mentioned), but most of the definitions and preprocessed material that you will need to understand the rest are to be found in this section. There is nothing for you to lose your appetite —you can actually intensify it, even if, or especially if, you do not agree with some of our positions. We first make clear what we mean by *logics*, introduced by their *consequence relations*, and what we mean by *theories* based on these logics, and on the way you will also learn what *closed theories* and *monotonicity* mean, and what it means to say that a logic is a *fragment* or an *extension* of another logic. This is just preparatory work. We then introduce the logical notions of *contradictoriness*, *explosiveness* and *triviality*, concerning theories and logics, and pinpoint some immediate connections between these notions. This already takes us to the subsection **2.1**, where the first logical (meta)-principles are introduced, namely the principles of non-contradiction, (PNC), of non-triviality, (PNT), and of explosion, (PPS) (a.k.a. Pseudo-Scotus, or *ex contradictione*). You will even learn a little bit about the (pre-)history of these principles, their interrelations, and some confusions about them which lurk around. Some of their ontological aspects are also lightly touched. The subsection **2.2** brings us to *paraconsistent logics*, formulated in two equivalent (but not necessarily so) presentations, one of them saying that they should allow for contradictory non-trivial theories, the other one saying that they must disrespect (PPS). After you learn what it means to say that two given formulas / theories are *equivalent* inside some given logic, FACT 2.8 will call your attention to the discrimination that paraconsistent logics ought to make between contradictions: they cannot be all equivalent inside such logics. *Dialectical logics*, being those logics disrespecting (PNC), are mentioned to fill the gaps in the general picture, but they will not be studied here. In the subsection **2.3** we start talking about finite trivializability, and look at some remarkable examples of this phenomenon, as for example the one of a logic having *bottom particles* —thus respecting a principle that we call *ex falso*, (ExF)—, and the one of a logic having *strong negations* —and respecting a principle we call *supplementing explosion*, (sPPS). We also consider some properties of adjunction (and so, of conjunction), and in the end we draw a map to show the relationships between (PPS), finite trivializability, (ExF) and (sPPS), noting that no two of these principles are to be necessarily identified (and, in particular, *ex falso* does not coincide with *ex contradictione*). Pay special attention to FACT 2.10(ii), in which all non-trivial logics respecting *ex falso* are shown to have strong negations. Subsection **2.4** considers what happens when one says farewell to (PPS) but still maintains some of the other special forms of explosion exposed be-fore and hints are given as to some disadvantages presented by paraconsistent logics which disrespect all of those principles at once. Some other misunderstandings about the construction of paraconsistent logics are discussed, and the difference between contradictoriness and inconsistency is finally called into scene. Logics respecting a so-called principle of *gentle explosion*, (gPPS), are introduced as the ones in which *consistency* can be expressed, and even a *finite* version of gentle explosion, (fgPPS), is considered, as a particular case of finite trivialization. *Logics of formal inconsistency*, **LFI**s, are then defined to be exactly those respecting (gPPS) while disrespecting (PPS), and the great majority of the **LFI**s that we will be studying in the following will actually also respect (fgPPS). The new definition of consistency that we introduce is shown, for a given logic, in general to coincide simply with the sum of (PPS) and (PNT). Systems of paraconsistent logic known as *discussive* (or *discursive*) *logics* are shown to be representable as **LFI**s. In the subsection **2.5**, the principles of *partial explosion*, (pPPS), and *controllable explosion*, (cPPS), are finally introduced. A *boldly paraconsistent* logic is defined as one in which not only (PPS) but also (pPPS) is disrespected, and we try to concentrate exclusively on such logics. Classical logic respects all of the above principles, but for each of those principles, except (PNT), examples will be explicitly presented, or at least referred to, at some point or another, of logics disrespecting it. Multiple connections between those principles are exhibited not only along these lines but also in the section **3**. In this respect, pay also special attention to FACT 2.19 and the comments around it, in the subsection **2.6**, which show that the **LFI**s are ubiquitous: an enormous subclass of the already known paraconsistent logics can have its members recast as logics of formal inconsistency. **C**-systems are also introduced in this last subsection, and the map of the paraconsistent land presented in the subsection **2.4** gets richer and richer.

Section **3** brings a very careful syntactical study of a large class of **C**-systems based on classical logic. Each new axiom is justified as it is introduced, and its effects and counter-effects are exhibited and discussed. The systems presented are initially linearly ordered by extension, but soon spread out in many directions. The remarkable unifying character of our approach in terms of **LFI**s is made clear while most logics produced by the 'Brazilian school' in the last forty years or so are shown to smoothly fit the general schema and together make up a whole coherent map of **C**-land. Subsection **3.1** presents a kind of minimal paraconsistent logic (for our purposes), called $C_{min}$ and constructed from the positive part of classical logic by the addition of (the axiom which represents the principle of) excluded middle, plus an axiom for double negation elimination. The Deduction Metatheorem holds for this logic and its extensions, and $C_{min}$ is shown to be paraconsistent. Comparisons are drawn between $C_{min}$ and one of its fragments, da Costa's $C_\omega$, and the facts that no strong negation, or bottom particle, or finitely trivializable theory, or negated theorem are to be found in these logics are mentioned. You will also learn that, in these logics, no two different negated formulas are provably equivalent. Of course, as a consequence of these last facts, these logics cannot be **LFI**s, what to say **C**-systems, but the **C**-systems which will be studied in the following subsections are all extensions of them. We make some observations about versions of *proof by cases* provable in these logics, by way of excluded middle, and we adjust some of its axioms to better suit deduction. A way to turn these logics into classical logic, simply by adding back the Pseudo-Scotus to them, is also demonstrated.

In subsection **3.2**, we introduce the *basic logic of (in)consistency*, called **bC**, by adding a new axiom to $C_{min}$, and we show how to immediately extract from this axiom a strong negation and a bottom particle. We now have '∘', the *consistency* connective, at our disposal as a new primitive constructor in our language, realizing the finite gentle explosion. The logic **bC**, which is, in fact, a conservative extension of $C_{min}$, is shown already to have negated theorems and equivalent negated formulas, but on the other hand it does not have any provably consistent formulas. Sufficient and necessary conditions for a **bC**-theory to behave classically are presented. The axiom defining **bC** define a kind of restricted *Pseudo-Scotus*, as obvious. Some related restricted forms of *reductio ad absurdum* which are also present are studied, and the elimination of double negation shows its purpose in THEOREM 3.13, where you will learn that some forms of partial trivialization are avoided by all paraconsistent extensions of **bC**. Restricted forms of *reductio* deduction and inference rules are shown to be present in **bC**, and some other rules relating contradictions and consistency are exhibited. No paraconsistent extensions of **bC** will contain the formula $(A \wedge \neg A)$ as a bottom particle (but some other **LFI**s, such as Jaśkowski's **D2** —at least under some presentations— do have it as a bottom particle). A formula such as $\neg(A \wedge \neg A)$ is also not provable in **bC**, but can be proved in some of its extensions, such as the three-valued maximal paraconsistent logics **LFI1** and **LFI2**.

Subsection **3.3** is mostly composed of negative results. It starts by showing that not many rules making the interdefinition of connectives possible hold in **bC**. The reader will also learn about the obligatory failure of *disjunctive syllogism* in vast extensions of the paraconsistent land, and the failure of 'full' *contraposition* inference rules in **bC** and all of its paraconsistent extensions, though some *restricted* forms of it had already been shown to hold, in the previous subsection. The uses of disjunctive syllogism and of contraposition to derive the Pseudo-Scotus had already been pointed out a long ago, respectively, by C. I. Lewis (and, much before, by the 'Pseudo-Scotus' himself), and by Popper. Some asymmetries related to negation of equivalent formulas are pointed out, and as a result it will not be possible to prove a *replacement* theorem for **bC**, which would establish the validity of *intersubstitutivity of provable equivalents*, (IpE), and the same phenomenon will be observed in most, but not all, of **bC**'s extensions. Reasons for all these failures, and possible solutions for them, are discussed.

In subsection **3.4** the problem of adding an inconsistency connective '•' to **bC**, intended as dual to '∘', making inconsistency coincide with non-consistency, and consistency coincide with non-inconsistency, is shown to be not as easy as it may seem, as a consequence of the last negative results. Some intermediary logics obtained in the strive towards the solution of this problem are exhibited, and hints are given on how this solution should look. A first solution, adopted, in fact, in the whole literature, is to make inconsistency equivalent to contradiction, and this is exactly what the logic **Ci**, introduced in the subsection **3.5**, does. It will *not* be the case, however, that consistency in **Ci** can be identified with the negation of a conjunction of contradictory formulas, that is, the consistency of a formula $A$ will not be equivalent to any formula such as $\neg(A \wedge \neg A)$. New forms of gentle explosion and restricted contraposition deduction rules are shown to hold in **Ci**, and provably consistent formulas in **Ci** are shown to exist. Indeed, the notable FACT 3.32 shows that provably consistent formulas in **Ci** coincide with the formulas causing controllable explosion in this logic and in all of its extensions. Some more restricted forms of contraposition inference rules introduced by **Ci** are also exhibited, and the failure of (IpE) also for this logic is pointed out. In fact, as we already know that full contraposition cannot be added to **Ci** in order to get (IpE), some weaker contraposition deduction rules which would also do the job are tested, and these are also shown to lead to collapse into classical logic (see FACT 3.36). But there still can be a chance of obtaining (IpE) in extensions of **Ci** by the addition of even weaker forms of contraposition deduction rules, as the reader is going to see at the end of subsection **3.7**, where positive results for this are shown for extensions of **bC**. The connectives '∘' and '•' have a good behavior inside of **Ci**, and we show that any formula having them as the main operator is consistently provable, and that consistency propagates through negation (and inconsistency back-propagates through negation). In addition to that, schemas such as $(A \rightarrow \neg\neg A)$, which are shown *not* to hold in the general case, are indeed shown to hold if $A$ has the form $\circ B$, for some $B$. All of this comes either as an effect or a consequence of the fact that a restricted replacement is valid in **Ci**, as proved in the subsection **3.6**, to the effect that inconsistency here can be really introduced, by definition, as non-consistency, or else consistency can be introduced, by definition, as non-inconsistency.

Subsection **3.7** shows how to compare the previously introduced **C**-systems with an extended version of classical logic, **eCPL** (adding innocuous operators for consistency and inconsistency). As a result, we can show that the strong negation that we had defined for **bC** does not have all properties of classical negation, but another strong negation can be defined in **bC**, which *does* have a classical character. In **Ci** these two negations are shown to be equivalent, but the interesting output of a strong classical negation is making it possible for us to conservatively translate classical logic inside of all our **C**-systems (THEOREM 3.46 and comments after THEOREM 3.48), so that any classical inference can be faithfully reproduced, up to a translation, inside of **bC** or of any of its extensions. About **Ci**, we can now prove that it has some redundant axioms, and the remarkable FACT 3.50, showing that only consistent or inconsistent formulas can themselves be consistently provable in this logic, and so these are the only formulas that can cause controllable explosion in **Ci**. An even more remarkable THEOREM 3.51 shows several conditions which *cannot* be fulfilled by paraconsistent systems in order to render the proof of full replacement, (IpE), possible. But paraconsistent extensions of **bC** in which (IpE) holds are indeed shown to exist, the same task remaining open for extensions of **Ci**.

Subsection **3.8** presents the **dC**-*systems*, which are the **C**-systems in which the connectives '∘' and '•' can be dispensed, definable from some combination of the remaining connectives. The particular combinations chosen by da Costa in the construction of his calculi $C_n$ are surveyed, and we start concentrating more and more on general parallels of da Costa's original requisites for the construction of paraconsistent calculi (which does not mean that we shall feel obliged to obey them *ipsis litteris*). Criticisms on the particular choices made by da Costa and some of their consequences are surveyed. Again, in a particular case, that of da Costa's $C_1$, the consistency of a formula $A$ is identified, as we have said before, with the formula $\neg(A \wedge \neg A)$, and the extension of **Ci** which makes this identification is called **Cil**. This system is shown still to suffer some strange asymmetries related to the negations of equivalent formulas (for instance, $\neg(\neg A \wedge A)$ is not equivalent to $\neg(A \wedge \neg A)$), and some partial or full solutions to that are discussed at that point and below. Connections between our **C**-systems and *relevance* or *intuitionistic* logics are touched. In particular, the problem of defining **C**-systems based on intuitionistic, rather than classical, logic is touched, but no interesting solutions are presented (because they seem still not to exist, but perhaps the reader will have the pleasure of finding them in the future).

In subsection **3.9**, **dC**-systems are put aside for a moment and the addition of an axiom for the introduction of double negation is considered, together with its consequences. Arguments, both positive and negative, for the 'proliferation of inconsistencies' that such an axiom could cause are presented, and rejected. Subsection **3.10** surveys various ways in which consistency (and inconsistency) can propagate from simpler to more complex formula and vice-versa. One of these forms, perhaps the most basic one, is illustrated by an extension of **Ci**, the logic **Cia** (or else an extension of **Cia**, the logic **Cila**, which was recorded into history under the name $C_1$), and this logic will be shown to make possible a new and interesting conservative translation from classical logic inside of it, or any of its extensions. If all the reader wants to know about is da Costa's original calculi $C_n$, this subsection is the place (together with some earlier comments in the subsection **3.7**), but in that case be warned: You may miss most of the fun! Properties of the $C_n$ are surveyed, and the problem of finding a real *deductive limit* to this hierarchy is presented together with its solution, so that the reader can forget once and for all any ideas they may have had about the logic $C_\omega$ having its place as part of this hierarchy. In particular, the deductive limit of the $C_n$, the logic $C_{Lim}$, is shown to constitute an **LFI**, though we are not sure if its form of gentle explosion can be made finite. Again, (IpE) is shown not to hold for the calculi $C_n$, so that *Lindenbaum-Tarski*-like algebraizations for these logics can be forgotten, but the situation for them is actually worse, for it has been proven that they just cannot define any non-trivial congruence, putting aside also the possibility of finding *Blok-Pigozzi*-like algebraizations for them. But several extensions of the $C_n$ can indeed fix this last problem, and so we try to concentrate on some stronger forms of propagation of consistency which will help us with this. In particular, the logics $C_1^+$ (later proposed by da Costa and his disciples), as well as five three-valued logics, $\mathbf{P}^1$, $\mathbf{P}^2$, $\mathbf{P}^3$, **LFI1** and **LFI2**, proposed in several studies, are also axiomatized and studied as extensions of **Ci**. An increasingly detailed and clear map of the **C**-systems based on classical logic is being drawn.

In the subsection **3.11**, the last five three-valued logics are shown to constitute part of a much larger family of 8,192 three-valued paraconsistent logics, each of them proven to be axiomatizable as extensions of **Ci** containing suitable axioms for propagation of consistency. Each of these three-valued logics can also be shown to be *distinct* from all of the others, and *maximal* relative to classical logic, **eCPL**, solving one of the main requisites set down by da Costa, to the effect that 'most rules and schemas of classical logic' should be provable in a 'good' paraconsistent logic. We also count how many of those 8,192 logics are in fact **dC**-systems, and not only **C**-systems, and show many *connections* between them and the other logics presented before, all of them fragments of some of these three-valued maximal paraconsistent logics. Interestingly, $\mathbf{P}^i$ is shown to be conservatively translatable inside of any of the other three-valued logics, and all of these are shown to be conservatively translatable inside of **LFI1**. (IpE) is proven not to hold in any of these logics, but there are some other interesting connectives which they can define, as some sort of *'highly' classical* negation, and *congruences* which will make possible, in subsection **3.12**, the definition of non-trivial (Blok-Pigozzi) algebraizations of all of these three-valued logics. Indeed, the subsection **3.12** surveys positive and negative results regarding algebraizations of the **C**-systems.

Section **4** sets some exciting open problems and directions for further research, for the reader's recreation.

**1.4 Standing on the shoulders of each other.** It would be very unwise of us to present this study, which includes a technical survey of its area, without trying to connect it as much as possible to the rest of the related literature. But we pledge to have done our very best to highlight, wherever opportune, some of the relevant papers which come close or very close to our points, or on which we simply base our study at some points! As in the case of the famous legendary caliph who set the books of the library of Alexandria on fire, we could say that the relevant papers which are not cited here at some point or another are, in most cases, either *blasphemous* (meaning that there is no context here for them to be mentioned), or *unnecessary* (meaning that in general you have to go no farther than the bibliography of our bibliography to learn about them). Other options would be our total *ignorance* about such and such papers at the moment of writing this (about which we thank for any enlightenment that we might receive), or else because we felt it was already well represented by another publication on the same matter, or because it integrates our list of *future research* (that was a good excuse, wasn't it?). Or perhaps it was *not* relevant at all! (You wouldn't know that, would you?) Read the text and judge for yourself. We just want to mention in this subsection a few other papers whose structural or methodological similarity (or dissimilarity) to some of our themes is most striking —so that we can better highlight our *own* originality on some topics, whenever it becomes the case.

Our study in section **2** is totally situated at the level of a general theory of consequence relations, a field sometimes referred to as that of *General Abstract Logics* (cf. [111]), or *Universal Logic* (cf. [19]). There are a few (rare) papers dealing with the definitions of the logical principles at a purely logical level. One of them is Restall's [99], where an approach to the matter quite different from ours is tackled. Also starting from the definition of logic as determined by its consequence relation (even though monotonicity is not pressuposed), and assuming from the start that an adjunctive conjunction is available, the author also requires one sort of contraposition deduction rule to be valid for all of the negations he considers (something that, in most logics herein studied, does not hold), and fixes the relevance logic $R$ in the formulation of most of his results. Several versions of the 'law' of non-contradiction are then presented, starting from the outright identification of this principle with the principle of explosion, passing through the identification of the principle of non-contradiction with the principle of excluded middle, or with the validity of the formula $\neg(A \wedge \neg A)$, and going up to some sort of difference in degree between accepting the inference of *all* propositions at once from some given formula $(A \wedge \neg A)$ instead of accepting *each* at a turn. It is clear that the outcome of all this is completely diverse from what we propose here. Some other studies go so far as to also study some of the alternative forms of explosion that we concentrate on here. This is the case, for instance, for Batens's [10], and Urbas's [108]. We are unaware, however, of any study which has taken these alternative forms as far as we do, and have studied them in precise and detailed terms. Such a study is presumed to be essential to help clarify the foundations, the nature and the reach of paraconsistent logics. There has been, for instance, arguments to the effect that the negations of paraconsistent logics are not (or may not be) negation operators after all (cf. Slater's [104] and Béziau's [23]). Béziau's argument amounts to a request for the definition of some minimal 'positive properties' in order to characterize paraconsistent negation really as constituting a *negation* operator, instead of something else. Slater argues for the *inexistence* of paraconsistent logics, given that their negation operator is not a 'contradictory forming functor', but just a 'subcontrary forming one', recovering and extending an earlier argument from Priest & Routley in [93]. Evidently, the same argument about not being a contradictory forming functor applies as well to intuitionistic negation, or in general to any other negation which does not have a classical behavior. Regardless of whether you wish to call such an operator 'negation' or something else, the negations of paraconsistent logics had better be studied under a less biased perspective, by the investigation of general properties that they can or cannot display inside paraconsistent logics. Some good examples of that kind of critical study may be found not only in the present paper but also in Avron's [9], Béziau's [18], and especially Lenzen's [70], among others.

In section **3** we investigate **C**-systems based on classical logic. In this respect, the present study has at least one very important ancestor, namely Batens's [10], where a general investigation of logics extending the positive classical logic (not all of them being **C**-systems!) is presented. This same author has also presented, elsewhere one of the best arguments that may be used to support our approach in terms of logics of formal inconsistency, **LFI**s. Criticizing Priest's logic *LP* (cf. [90]), Batens insists that:

> There simply is no way to express, within this logic, that *B* is *not* false or that *B* behaves consistently. (Cf. [13], p.216)

Asserting that 'paraconsistent negation should not and cannot express rejection' (id., p.223), Batens wants to say that it is not because a negated sentence $\neg B$ is inferred from some non-trivial theory of a paraconsistent logic that we can conclude that *B* is not also to be inferred from that, i.e. $\neg B$ does not express the *rejection* of *B*. From that he will draw several lessons along his article, such as that: (i) the presence of a strong negation (he writes 'classical negation', but this is clearly an extrapolation —see our note 15, in the subsection **3.7**) inside of a paraconsistent logic is not only a sufficient requisite but also a necessary one to express (classical) rejection; (ii) that one needs a controllably explosive paraconsistent logic (he calls it 'non-strictly paraconsistent') to be able to 'fully describe classical logic' (id., p.225); (iii) that the existence of a bottom particle is also sufficient for the above purposes, for it may define a strong negation (this appeared in an addenda to the paper). So, all at once, this author argues for the validity of three of our alternative explosion principles: (sPPS), (cPPS), and (ExF). From these, of course, we already know (see FACT 2.19) that our principle (gPPS) will often follow, so that according to his recommendations we are finally left with **LFI**s, instead of something like Priest's logic, which, again according to Batens, and for the above reasons, 'fails to capture natural thinking' (though it was proposed to such an effect), and does not provide sufficient environment for us to do 'paraconsistent mathematics'. Priest's response to such criticisms seems to us to be somewhat of a cheat, for he proposes to introduce such a strong negation using an ill-defined bottom particle (see note 25, in the subsection **3.10**). The only point where Batens goes too far to be right seems to be on his argument about paraconsistent logics not being adequate to be used on our metalanguage, because we would be in need of strong negations to complete any consistent description of the world. But now we know that an **LFI** would be more than enough to such an end, being able to fully reproduce all the classical inferences (THEOREM 3.46). And do remember that **LFI**s are especially tailored in order to *express* the fact that $B$ behaves consistently, attending, thus, to Batens's requisite above (and also to his praxis, given that he has already been using in his articles, since long, some symbol to express inconsistency in the object language, be it just an abbreviation or some sort of metaconnective —check the symbol '!' in [15]).

The very idea of a paraconsistent logic still has, nowadays, as strong defenders as attackers. Though, as we know, many attacks are but misunderstandings, many defenses are also poor or unsound. We hope here to contribute to this debate, in one way or another, combining as much precision and clarity as we can. At a more fundamental philosophical level, also, paraconsistency has raised diverse excited opinions about the contribution (or damage) it makes to the very notion of *rationality*. An author such as Mario Bunge will on the one hand compare the Pseudo-Scotus with some sort of *cancer* (cf. [32], p.17), and on the other hand observe that 'a refined symbolism can hide a brazen irrationalism' (id., p.23). He asserts that paraconsistent logic is non-rational by definition, because 'it does not include the principle of non-contradiction' (id., p.24). About this same point other authors will concede similar verdicts, and yet arrive at different conclusions. As Gilles-Gaston Granger put it, paraconsistent logics can be seen as a 'provisory recourse to the irrational', for maintaining an indicium of the rational (sic), namely the principle of non-trivialization, while also maintaining an indicium of the irrational, namely the possible presence of contradictions, to be 'philosophically justified' (cf. [65], p.175). Yet some other authors, such as Newton da Costa, defended that, according to some pragmatic principles of reason which 'seem to be present in all processes of systematization of rational knowledge', this same rational knowledge can be said, among other things, to be both intuitive and discursive, to result from the interaction of the spirit and its environs, and not to be identifiable with a particular system of logic. About reason, on its own turn, he maintained that it is tied to its historical evolution, has its range of application determinable only pragmatically, and is always expressible by way of some logic, which, in each case, is supposed to be uniquely determined by each given context, as being precisely the logic that is most adequate to that context. To determine the concept of adequacy, finally, da Costa recurs again to pragmatic factors, such as simplicity, convenience, facility, economy, and so on (cf. [51]). This is in fact a very thought-provoking issue, and several other authors have advanced positions on the relations of paraconsistency and rationality, such as Francisco Miró Quesada (in [81] and [80]), Nicholas Rescher (in [98]), Jean-Yves Béziau (in [23]), and Bobenrieth (in [26]). The reader is invited to read those authors directly, if this is their interest. We will not venture here any further steps in this slippery slope, for our aim is much less ambitious. After reading this comprehensive technical survey, however, we hope that the reader will feel illuminated enough to risk their own rationally-based judgements on the matter.

> Logic is the chosen resort of clear-headed people, severally convinced of
> the complete adequacy of their doctrines. It is such a pity that they cannot
> agree with each other.
> —A. N. Whitehead, "Harvard: The Future", Atlantic Monthly 158, p. 263.

Many a logician will agree that the fundamental notion behind logic is the notion of
'derivation', or rather should we say the notion of 'consequence'. On that account, in
our common heritage it is to be found the Tarskian notion of a *consequence relation*.
As usual, given a set *For* of formulas, we say that $\Vdash \subseteq \wp(For) \times For$ defines a con-
sequence relation on *For* if the following clauses hold, for any formulas $A$ and $B$,
and subsets $\Gamma$ and $\Delta$ of *For:* (formulas and commas at the left-hand side of $\Vdash$ de-
note, as usual, sets and unions of sets of formulas)

(Con1)  $A \in \Gamma \implies \Gamma \Vdash A$                                          (reflexivity)

(Con2)  $(\Delta \Vdash A$ and $\Delta \subseteq \Gamma) \implies \Gamma \Vdash A$                    (monotonicity)

(Con3)  $(\Delta \Vdash A$ and $\Gamma, A \Vdash B) \implies \Delta, \Gamma \Vdash B$          (transitivity)

So, a logic **L** will here be defined simply as a structure of the form $<For, \Vdash>$, con-
taining a set of formulas and a consequence relation defined on this set. We need not
suppose at this point that the set *For* should be endowed with any additional structure,
like the usual algebraic one, but we will hereby suppose, for convenience, that *For*
is built on a denumerable language having $\neg$ as its (primitive or defined) negation
symbol, and we will also suppose the connectives to be constructing operators on the set
of formulas. Any set $\Gamma \subseteq For$ is called a *theory* of **L**. A theory $\Gamma$ is said to be *proper*
if $\Gamma \neq For$, and a theory $\Gamma$ is said to be *closed* if it contains all of its consequences,
i.e. if the converse of (Con1) holds: $\Gamma \Vdash A \implies A \in \Gamma$. Whenever we have, in a
given logic, that $\Gamma \Vdash A$, for a given theory $\Gamma$ and some formula $A$, we will say that $A$
is *inferred* from $\Gamma$ (in this logic); if, for all $\Gamma$, we have that $\Gamma \Vdash A$, that is, if $A$ is
inferred from any given theory, we will say that $A$ is a *thesis* (of this logic).

Not all known logics respect all the above clauses, or only them. For instance, those
logics in which (Con2) is either dropped out or substituted by a form of 'cautious
monotonicity' are called *non-monotonic*, and the logics whose consequence relations
are closed under substitution are called *structural*. Unless explicitly stated to the con-
trary, we will from now on be working with some fixed arbitrary logic **L** = $<For, \Vdash>$,
and with some fixed arbitrary theory $\Gamma$ of **L**. Properties (Con1)–(Con3) will be as-
sumed to hold irrestrictedly, and they will be used in some proofs here and there.
Some interesting and quite immediate consequences from (Con2) and (Con3) which
we shall make use of are the following:

FACT 2.1  The following properties hold for any logic, any given theories $\Gamma$ and $\Delta$,
  and any formulas $A$ and $B$:
  (i)  $\Gamma, \Delta \nVdash A \implies \Gamma \nVdash A$;
  (ii)  $(\Gamma \Vdash A$ and $A \Vdash B) \implies \Gamma \Vdash B$;
  (iii)  $(\Gamma \Vdash A$ and $\Gamma, A \Vdash B) \implies \Gamma \Vdash B$.
**Proof:**   (i) follows from (Con2); (ii) and (iii), from (Con3).                    □

Given two logics **L**1 = $<For_1, \Vdash_1>$ and **L**2 = $<For_2, \Vdash_2>$, we will say that **L**1 is a
*linguistic extension* of **L**2 if $For_2$ is a proper subset of $For_1$, and we will say that **L**1

is a *deductive extension* of **L2** if $\Vdash_2$ is a proper subset of $\Vdash_1$. Finally, if **L1** is both a linguistic and deductive extension of **L2**, and if the restriction of **L1**'s consequence relation $\Vdash_1$ to the set $For_2$ will make it identical to $\Vdash_2$ (that is, if $For_2 \subset For_1$, and for any $\Gamma \cup \{A\} \subseteq For_2$ we have that $\Gamma \Vdash_1 A \Leftrightarrow \Gamma \Vdash_2 A$) then we will say that **L1** is a *conservative extension* of **L2**. In any of the above cases we can more generally say that **L1** is an *extension* of **L2**, or that **L2** is a *fragment* of **L1**. These concepts will be used mainly in the next section, where we will build and compare a number of paraconsistent logics. Just as a guiding note to the reader, however, we could remark that usually, but not obligatorily, linguistic extensions are also deductive ones, but it is quite easy to find in the realm of non-classical logics, on the other hand, deductive fragments which are not linguistic ones (like intuitionistic logic is a deductive fragment of classical logic). Most paraconsistent logics in the literature are also deductive fragments of classical logic themselves, but the ones we shall be working on here, the **C**-*systems*, are in general deductive fragments only of a conservative extension of classical logic —by the addition of (explicitly definable) connectives expressing consistency / inconsistency). A particular case of them, the **dC**-systems, will nevertheless be shown to be characterizable as deductive fragments of good old classical logic, dispensing its mentioned extension. But these assertions will be made much clearer in the near future.

Let $\Gamma$ be a theory of **L**. We say that $\Gamma$ is *contradictory with respect to* $\neg$, or simply *contradictory*, if it is such that, for some formula $A$, we have $\Gamma \Vdash A$ and $\Gamma \Vdash \neg A$. With some abuse of notation, but (hopefully) no risk of misunderstanding, we will from now on write these sort of sentences in the following way:

$$\exists A \,(\Gamma \Vdash A \;\text{ and }\; \Gamma \Vdash \neg A). \tag{D1}$$

For any such formula $A$ we may also say that $\Gamma$ is *$A$-contradictory*, or simply that $A$ is *contradictory* for such a theory $\Gamma$ (and such an underlying logic **L**). It follows that:

**FACT 2.2** For a given theory $\Gamma$: (i) If $\{A, \neg A\} \subseteq \Gamma$ then $\Gamma$ is $A$-contradictory. (ii) If $\Gamma$ is both $A$-contradictory and closed, then $\{A, \neg A\} \subseteq \Gamma$.

**Proof:** Part (i) comes from (Con1), part (ii) from the very definition of a closed theory. □

A theory $\Gamma$ is said to be *trivial* if it is such that:

$$\forall B \,(\Gamma \Vdash B). \tag{D2}$$

Hence, a trivial theory can make no difference between the formulas of a logic — all of them may be inferred from it. Of course, using (Con1) we may notice that the non-proper theory *For* is trivial. We may also immediately conclude that:

**FACT 2.3** Contradictoriness is a necessary condition for triviality in a given theory. (D2) $\Rightarrow$ (D1)

A theory $\Gamma$ is said to be *explosive* if:

$$\forall A \,\forall B \,(\Gamma, A, \neg A \Vdash B). \tag{D3}$$

Thus, a theory is called explosive if it trivializes when exposed to a pair of contradictory formulas. Evidently:

**FACT 2.4** (i) If a theory is trivial, then it is explosive. (ii) If a theory is contradictory and explosive, then it is trivial. (D2) $\Rightarrow$ (D3); (D1) and (D3) $\Rightarrow$ (D2)

**Proof:** Use (Con2) in the first part and FACT 2.1(iii) in the second. □

**2.1 A question of principles.** Now, remember that talking about a logic is talking about the inferential behavior of a set of theories. Accordingly, using the above definitions, we will now say that a given logic **L** is *contradictory* if all of its theories are contradictory. (D4)

In much the same spirit, we will say that **L** is *trivial*, or *explosive*, if, respectively, all of its theories are trivial, or explosive. respect. (D5), (D6)

The empty theory may be here regarded as playing an important role, revealing some intrinsic properties of a given logic, in spite of the behavior of any of its specific non-empty theories (also called 'non-logical axioms'). Indeed:

FACT 2.5 A monotonic logic **L** is contradictory / trivial / explosive if, and only if, its empty theory is contradictory / trivial / explosive.

We can now tackle a formal definition for some of the so-called *logical principles* (relativized for a given logic **L**), namely:

PRINCIPLE OF NON-CONTRADICTION (PNC)
   **L** must be non-contradictory: $\exists\Gamma\,\forall A\,(\Gamma\nVdash A$ or $\Gamma\nVdash\neg A)$.

PRINCIPLE OF NON-TRIVIALITY (PNT)
   **L** must be non-trivial: $\exists\Gamma\,\exists B\,(\Gamma\nVdash B)$.

PRINCIPLE OF EXPLOSION, or PRINCIPLE OF PSEUDO-SCOTUS[3] (PPS)
   **L** must be explosive: $\forall\Gamma\,\forall A\,\forall B\,(\Gamma, A, \neg A \Vdash B)$.

This last principle is also often referred to as *ex contradictione sequitur quodlibet*. The reader will immediately notice that these principles are somewhat interrelated:

FACT 2.6 (i) An explosive logic is contradictory if, and only if, it is trivial. (ii) A trivial logic is both contradictory and explosive. (iii) A logic in which the Principle of Explosion holds is a trivial one if, and only if, the Principle of Non-Contradiction fails. (D6) ⇒ [(D4) ⇔ (D5)]; (D5) ⇒ [(D4) and (D6)];
(PPS) ⇒ [not-(PNT) ⇔ not-(PNC)]

**Proof:** Just consider FACT 2.4, and the definitions above. □

A trivial logic, i.e. a logic in which (PNT) fails, cannot be a very interesting one, for in such a logic anything could be inferred from anything, and any intended capability of modeling 'sensible' reasoning would then collapse. Of course, (PPS) would still hold in such a logic, as well as any other universally quantified sentence dealing with the behavior of its consequence relation, but this time only because they would be unfalsifiable! It is readily comprehensible then that triviality might have been regarded as the mathematician's worst nightmare. Indeed, (PNT) constituted what Hilbert called 'consistency (or compatibility) principle', with which proof his Metamathematical enterprise was crafted to cope. Well aware of the preceding fact, and working inside the environment of an explosive logic such as classical logic, Hilbert transposed the 'problem of consistency' (that is, the problem of non-triviality) to the problem of proving that there were no contradictions among the axioms of arithmetic and their consequences (this was Hilbert's Second Problem, cf. [66]). By the way, this situation would eventually lead Hilbert to the formulation of a curious criterion according to which the non-contradictoriness of a mathematical object is a necessary and sufficient condition for its very *existence*.[4] Perhaps he would have never proposed such a criterion if he had only considered the existence of non-explosive logics, with or without (PNC), logics in which contradictory theories do not necessarily lead to trivialization! The search for such logics would give rise, much later, to the 'paraconsistent enterprise'.[5]

---

[3] Which was made visible by a reedition of a collection of commentaries on Aristotle's *Prior Analytics*, long attributed, in error, to Johannes Duns Scotus (1266-1308), in the twelve books of the *Opera Omnia*, 1639 (reprint 1968). The current most plausible conjecture about the authorship of these books will trace them back to John of Cornwall, around 1350. See [26] for more on its history.

In classical logic, of course, all the three principles above hold, and one could naively speculate from such that they are all 'equivalent', in some sense. Indeed, they have all been now and again confused in the literature and each one of them has, in turn, been identified with the 'Principle of Non-Contradiction' (and these will not exhaust all formulations of this last principle that have been proposed here and there). The emergence of paraconsistent logic, as we shall see, will serve to show that this equivalence is far from being necessary, for an arbitrary logic **L**.

**2.2 The paraconsistency predicament.** Some decades ago, S. Jaśkowski ([67]) and N. C. A. da Costa ([49]), the founders of paraconsistent logic, proposed, independently, the study of logics which could accommodate contradictory yet non-trivial theories. Accordingly, a *paraconsistent logic* (a denomination which would be coined only in the seventies, by Miró Quesada) would be initially defined as a logic such that:

$$\exists \Gamma \exists A \exists B \ (\Gamma \Vdash A \ \text{and} \ \Gamma \Vdash \neg A \ \text{and} \ \Gamma \nVdash B). \tag{PL1}$$

*Attention:* This definition says *not* that (PNC) is not to hold in such a logic, for it says nothing about *all* theories of a paraconsistent logic being contradictory, but only that *some* of them should be contradictory, and yet non-trivial. As a consequence, following our definitions above, the notion of paraconsistent logic has, in principle, nothing to do with the rejection of the Principle of Non-Contradiction, as it is commonly held! On the other hand, it surely has something to do with the rejection of explosiveness. Indeed, consider the following alternative definition of a paraconsistent logic, as a logic in which (PPS) fails:

$$\exists \Gamma \exists A \exists B \ (\Gamma, A, \neg A \nVdash B). \tag{PL2}$$

Now one may easily check that:

FACT 2.7 (PL1) and (PL2) are equivalent ways of defining a paraconsistent logic, if its consequence relation is reflexive and transitive.

$$[(\text{Con1}) \ \text{and} \ (\text{Con3})] \Rightarrow [(\text{PL1}) \Leftrightarrow (\text{PL2})]$$

---

[4] Girolamo Saccheri (1667-1733) had already paved the way much before to set non-contradictoriness, instead of intuitiveness, as a sufficient, other than necessary, criterion for the legitimateness of a mathematical theory (cf. [1]). The so-called 'Hilbert's criterion for existence in mathematics' seems thus to constitute a further step in taking this method to its ultimate consequences.

[5] Assuming intuitively (cf. [46], p.7) that a contradiction could painlessly be admitted in a given theory if only this theory was not to be trivialized by it, even some years before the actual proposal of his first paraconsistent systems, da Costa was eventually led to trace the Metamathematical's problem about the utility of a formal system back to (PNT). At that point da Costa was even to suggest that Hilbert's criterium for existence in mathematics should be changed, and that *existence*, in mathematics, should be equated with non-triviality, rather than with non-contradictoriness (cf. [47], p.18). To be more precise about this point, da Costa has in fact recovered Quine's motto ([96], chapter I): 'to be is to be the value of a variable' —from which follows that the ontological commitment of our theories is to be measured by the domain of its variables—, and then proposed the following modification to it: 'to be is to be the value of a variable, in a given language of a given logic' (cf. [52], and the entry 'Paraconsistency' in [33]). This was meant to open space for the appearance of different ontologies based on different kinds of logic, analogously to what had happened in the XIX century with the appearance of different geometries based on different sets of axioms.

**Proof:** To show that (PL1) implies (PL2), use (Con3), or directly FACT 2.1; to show the converse, use (Con1). □

Say that two formulas $A$ and $B$ are *equivalent* if each one of them can be inferred from the other, that is:

$$(A \Vdash B) \text{ and } (B \Vdash A). \tag{Eq1}$$

In a similar manner, say that two sets of formulas $\Gamma$ and $\Delta$ are equivalent if:

$$\forall A \in \Delta \; (\Gamma \Vdash A) \text{ and } \forall B \in \Gamma \; (\Delta \Vdash B). \tag{Eq2}$$

We will alternatively denote these facts by writing, respectively, $A \dashv\Vdash B$, and $\Gamma \dashv\Vdash \Delta$.

Now, an essential trait of a paraconsistent logic is that it does not see all contradictions at the same light —each one is a different story. Indeed:

FACT 2.8 Given any arbitrary transitive paraconsistent logic, it cannot be the case that all of its contradictions are equivalent.

**Proof:** If, for whatever formulas $A$ and $B$, we have that $\{A, \neg A\} \dashv\Vdash \{B, \neg B\}$, then any $A$-contradictory theory, would also be, by transitivity and definition (Eq2), a $B$-contradictory theory. But if a theory infers every pair of contradictory formulas, it infers, in particular, any given formula at all, and so it is trivial. □

Once again, the reader should note that the existence of a paraconsistent logic **L** presupposes only the existence of *some* non-explosive theories in **L**; this does not mean that *all* theories of **L** should be non-explosive —and how could they all be so? (recall FACT 2.4(i)) Moreover, once more according to our proposed definitions, the reader will soon notice that the great majority of the paraconsistent logics found in the literature, and all the paraconsistent logics studied in this paper, are non-contradictory (i.e. 'consistent', following the usual model-theoretic connotation of the word). In particular, they usually have non-contradictory empty theories, which means, from a proof-theoretical point of view, that they bring no built-in contradiction in their axioms, and that their inference rules do not generate contradictions from these axioms. Even so, because of their paraconsistent character, they can still be used as underlying logics to extract some sensible reasoning of some theories that are contradictory and are still to be kept non-trivial. This phenomenon is no miracle, and certainly no sleight of hand, as the reader will understand below, but is obtained from suitable constraints on the power of explosiveness, (PPS). So, all paraconsistent logics which we will present here are in some sense 'more conservative' than classical logic, in the sense that they will extract less consequences than classical logic would extract from some given classical theory, or at most the same set of consequences, but never more. Our paraconsistent logics then (as most paraconsistent logics in the literature) will not validate any bizarre form of reasoning, and will not extract any contradictory consequence in the cases where classically there were no such consequences. It is nonetheless possible to also build logics which disrespect both (PPS) *and* (PNC), and thus might be said to be 'highly' non-classical, in a certain sense, once they *do* have theses which are not classical theses. Such logics will constitute a particular case of paraconsistent logics that are generally dubbed *dialectical logics*, or *logics of impossible objects*, and some specimens of these may be found, for instance, in [56], [83], [88], and [100]. We will not study these kind of logics here.

We shall, from now on, make use of either one of the above definitions for paraconsistent logic, indistinctly.

**2.3 The trivializing predicament.** Given (PL1), we know that any paraconsistent logic must possess contradictory non-trivial theories, and from (PL2) we know that these must be non-explosive. Evidently, not all theories of a given logic can be such: we already also know that any trivial theory is both contradictory and explosive, and every logic has trivial theories (consider, for instance, the non-proper theory *For*, i.e. the whole set of formulas). It is possible, though, and in fact very interesting, to further explore this no man's land which lays in between plain non-explosiveness and out-right explosiveness, if one considers some paraconsistent logics having some suitable explosive proper theories. That is what we will do in the following subsections.

A logic **L** is said to be *finitely trivializable* when it has finite trivial theories.     (D7)
Evidently:

FACT 2.9 If a logic is explosive, then it is finitely trivializable.     (D6) ⇒ (D7)
**Proof:** All theories of an explosive logic are explosive, in particular the empty one. Thus, for any *A*, the finite theory $\{A, \neg A\}$ is trivial.     □

This same fact does not hold for non-explosive logics. In fact, we will present, in the following, a few paraconsistent logics which are *not* finitely trivializable, although these shall, in general, not concern us in this article, for reasons which will soon be made clear. Let us first state and study some few more simple definitions.

A logic **L** has a *bottom particle* if there is some formula *C* in **L** that can, by itself, trivialize the logic, that is:

$$\exists C \, \forall \Gamma \, \forall B \; (\Gamma, C \Vdash B).$$     (D8)

We will denote any fixed such particle, when it exists, by ⊥. Evidently, no arbitrary monotonic and transitive logic can have a bottom particle as a thesis, under pain of turning this logic into a trivial logic —in which, of course, all formulas turn to be bottom particles.

It is instructive here to remember another formulation of (PPS) which sometimes shows up in the literature:

PRINCIPLE OF 'EX FALSO SEQUITUR QUODLIBET'     (ExF)
    **L** must have a bottom particle.

Now, if we are successful in isolating logics that disrespect (PPS) while still respecting (ExF) we will show that *ex contradictione* (*sequitur quodlibet*) does not need to be identified with *ex falso* (*sequitur quodlibet*), as is quite commonly held.[6]

We say that a logic **L** has a *top particle* if there is some formula *C* in **L** that is a consequence of every one of its theories, no matter what, that is:

$$\exists C \, \forall \Gamma \; (\Gamma \Vdash C).$$     (D9)

We will denote any fixed such particle, when it exists, by ⊤. Evidently, given a monotonic logic, any of its theses will constitute such a top particle (and logics with no theses, like Kleene's three-valued logic, will have no such particles). Also, given transitivity and monotonicity, it is easy to see that the addition of a top particle to a given theory is pretty innocuous, for in that case $(\Gamma, \top \Vdash B)$ if and only if $(\Gamma \Vdash B)$.

Let **L** now be some logic, let σ: *For* → *For* be a mapping (if *For* comes equipped with some additional structure, we will require σ to be an endomorphism), and let this mapping be such that σ(*A*) is to denote a formula which *depends only on A*. By this we shall mean that σ(*A*) is a formula constructed using but *A* itself and some

---

purely logical symbols (such as connectives, quantifiers, constants). In more general terms, given any sequence of formulas $A_1, A_2, ..., A_n$, we will let $\sigma(A_1, A_2, ..., A_n)$ denote a formula which depends only on the formulas of the sequence. Similarly, we will let $\Gamma(A_1, A_2, ..., A_n)$ denote a set of formulas each of which depends only on the sequence $A_1, A_2, ..., A_n$. In some situations it will help to assume this $\sigma$ to be a *schema*, that is, that given any two sequences $A_1, A_2, ..., A_n$ and $B_1, B_2, ..., B_n$, we must have that $\sigma(A_1, A_2, ..., A_n)$ will be made identical to $\sigma(B_1, B_2, ..., B_n)$ if we only change each $A_i$ for $B_i$, in $\sigma(A_1, A_2, ..., A_n)$ (this means, in some sense, that all these $\sigma$-formulas will share some built-in *logical form*). Usually, when saying that we have a formula, or set of formulas, depending only on some given sequence of formulas, we further presuppose that this dependency is schematic, but this supposition will in general be not strictly necessary to our purposes.

We say that a logic **L** has a *strong* (or *supplementing*) *negation* if there is a schema $\sigma(A)$, depending only on $A$, that does not consists, in general, of a bottom particle and that cannot be added to any theory inferring $A$ without causing its trivialization, that is:

(a) $\exists A$ such that $\sigma(A)$ is not a bottom particle, and
(b) $\forall A \, \forall \Gamma \, \forall B \, [\Gamma, A, \sigma(A) \Vdash B]$.  (D10)

We will denote the strong negation of a formula $A$, when it exists, by $\sim A$.

Parallel to the definition of contradictoriness with respect to $\neg$, we might now define a theory $\Gamma$ to be *contradictory with respect to* $\sim$ if it is such that:

$$\exists A \, (\Gamma \Vdash A \ \text{and} \ \Gamma \Vdash \sim A).$$  (D11)

Accordingly, a logic **L** is said to be *contradictory with respect to* $\sim$ if all of its theories are contradictory with respect to $\sim$.  (D12)

Here we may of course introduce yet another version of (PPS):

SUPPLEMENTING PRINCIPLE OF EXPLOSION  (sPPS)
    **L** must have a strong negation.[7]

Some immediate consequences of the last definitions are:

FACT 2.10  (i) If a logic has either a bottom particle or a strong negation, then it is finitely trivializable.  (ii) If a non-trivial logic has a bottom particle, then it admits a strong negation.  (iii) If a logic is explosive and non-trivial, then it is supplementing explosive.
[(D8) or (D10)] $\Rightarrow$ (D7); [not-(D5) and (D8)] $\Rightarrow$ (D10);
[(PNT) and (ExF)] $\Rightarrow$ (sPPS); [(PPS) and (PNT)] $\Rightarrow$ (sPPS)

**Proof:** (i) is obvious. To prove (ii), define the strong negation $\sim A$ of a formula $A$ by stipulating that, for any theories $\Gamma$ and $\Delta$, we have (a) $(\Gamma, \Delta \Vdash \sim A)$ iff $(\Gamma, \Delta, A \Vdash \bot)$, and (b) $(\Gamma, \Delta, \sim A \Vdash \bot)$ iff $(\Gamma, \Delta \Vdash A)$. By (Con1), we have that $(\Gamma, \sim A \Vdash \sim A)$, and so, part (a) will give us $(\Gamma, \sim A, A \Vdash \bot)$, choosing $\Delta = \{\sim A\}$. But $\bot \Vdash B$, for any formula $B$, once $\bot$ is a bottom particle. So, by FACT 2.1(ii), we conclude that $(\Gamma, A, \sim A \Vdash B)$, for any $B$. Now, to check that such a strong negation, thus defined, cannot be always a bottom particle, notice that part (b) will give us $\Vdash A$ iff $\sim A \Vdash \bot$, choosing both $\Gamma$ and $\Delta$ to be empty. So, if $\sim A$ were a bottom particle, $\sim A \Vdash \bot$ would be the case, and hence any $A$ would be a thesis of this logic, which is not the case, once we have supposed it to be non-trivial.[8] To check (iii), just note that a non-trivial ex-plosive logic will come already equipped with a built-in strong negation, coinciding with its own primitive negation. □

---

[7] A strong negation should *not* be confused with a 'classical' one! Take a look at THEOREM 3.42.
[8] In the presence of a convenient implication, for instance, obeying the Deduction Metatheorem (THEOREM 3.1) such an 'implicit' definition of a strong negation from a bottom particle can be internalized by the underlying logic as an 'explicit' definition (as in the case of intuitionistic logic). Check also our remarks about this matter in our discussion of Beth Definability Property, at the end of the subsection **3.12**.

FACT 2.11 Let **L** be a logic with a strong negation ~. (i) Every theory which is contradictory with respect to ~ is explosive. (ii) A logic is contradictory with respect to ~ if, and only if, it is trivial. (D11) ⇒ (D3); (D12) ⇔ (D5)

A logic **L** is said to be *left-adjunctive* if for any two formulas $A$ and $B$ there is a schema $\sigma(A, B)$, depending only on $A$ and $B$, with the following behavior:

(a) $\exists A \, \exists B$ such that $\sigma(A, B)$ is not a bottom particle, and
(b) $\forall A \, \forall B \, \forall \Gamma \, \forall D \, [\Gamma, A, B \Vdash D \;\Rightarrow\; \Gamma, \sigma(A, B) \Vdash D]$.     (D13)

Such a formula, when it exists, will be denoted by $(A \wedge B)$, and the sign $\wedge$ will be called a *left-adjunctive conjunction* (but it will not necessarily have, of course, all properties of a classical conjunction). Similarly, a logic **L** is said to be *left-disadjunctive* if there is a schema $\sigma(A, B)$, depending only on $A$ and $B$, such that (D12) is somewhat inverted, that is:

(a) $\exists A \, \exists B$ such that $\sigma(A, B)$ is not a top particle, and
(b) $\forall A \, \forall B \, \forall \Gamma \, \forall D \, [\Gamma, \sigma(A, B) \Vdash D \;\Rightarrow\; \Gamma, A, B \Vdash D]$.     (D14)

In general, whenever there is no risk of misunderstanding, we might also denote this formula, when it exists, by $A \wedge B$, and we will accordingly call $\wedge$ a *left-disadjunctive conjunction*. Now, one should be aware of the fact that, in principle, a logic can have just one of these conjunctions, or it can have both a left-adjunctive and a left-disadjunctive conjunction without the two of them coinciding.

To convince themselves of the naturalness of these definitions and the comments we made about them, we invite the reader to consider the following two more 'concrete' properties of conjunction:

(a) $\exists A \, \exists B$ such that $A \wedge B$ is not a bottom particle, and
(b) $\forall \Gamma \, \forall A \, \forall B \, (\Gamma, A \wedge B \Vdash A \;\; \text{and} \;\; \Gamma, A \wedge B \Vdash B)$.     (pC1)

(a) $\exists A \, \exists B$ such that $A \wedge B$ is not a top particle, and
(b) $\forall \Gamma \, \forall A \, \forall B \, (\Gamma, A, B \Vdash A \wedge B)$.     (pC2)

Now, it is easy to see that:

FACT 2.12 Let **L** be a logic obeying (Con1)–(Con3). (i) A conjunction in **L** is left-adjunctive iff it respects (pC1). (ii) A conjunction in **L** is left-disadjunctive iff it respects (pC2). [(Con1)–(Con3)] ⇒ {[(D13) ⇔ (pC1)] and [(D14) ⇔ (pC2)]}
**Proof:** To prove that a left-adjunctive conjunction respects (pC1) and that a left-disadjunctive conjunction respects (pC2), use (Con1) and (Con2). For the converses, use (Con3). □

The reader might mind to notice that a conjunction which is both left-adjunctive and left-disadjunctive is sometimes called, in the scope of relevance logic, an *intensional* conjunction, and in the scope of linear logic such a conjunction is said to be a *multiplicative* one (also, in [9], this is what the author calls an *internal* conjunction).

We may now check that:

FACT 2.13 Let **L** be a left-adjunctive logic. (i) If **L** either is finitely trivializable or has a strong negation, than it has a bottom particle. (ii) If **L** is finitely trivializable, then it will be supplementing explosive. (iii) If **L** respects *ex contradictione*, then it will respect *ex falso*. (D13) ⇒ {[(D7) or (D10)] ⇒ (D8)};
[(D13) and (D7)] ⇒ (sPPS); (D13) ⇒ [(PPS) ⇒ (ExF)]

**Proof:** To prove (i), note that if **L** has a finite trivial theory $\Gamma$, one may define a bottom particle from the conjunction of all formulas in $\Gamma$; in case it has a strong negation, any formula in the form $(A \wedge {\sim} A)$, for some formula $A$ of **L**, will suffice. Parts (ii) and (iii) are immediate. □

Consider now the *discussive logic* proposed by Jaśkowski in [67], **D2**, which is such that $\Gamma \vDash_{\mathbf{D2}} A$ iff $\Diamond\Gamma \vDash_{S5} \Diamond A$, where $\Diamond\Gamma = \{\Diamond B\colon$ for all $B \in \Gamma\}$, $\Diamond$ denotes the possibility operator, and $\vDash_{S5}$ denotes the consequence relation defined by the well-known modal logic $S5$. It is easy to see that in **D2** one has that $(A, \neg A \vDash_{\mathbf{D2}} B)$ *does not* hold in general, though $(A \wedge \neg A) \vDash_{\mathbf{D2}} B$ *does* hold, for any formulas $A$ and $B$. This phenomenon can only happen because (pC1) holds while (pC2) does not hold in **D2**, and so its conjunction is left-adjunctive but not left-disadjunctive, while $(A \wedge \neg A)$ defines a bottom particle. Hence, the fact above still holds for **D2**, and this logic indeed displays a quite immediate example of a logic respecting (ExF) but not (PPS).

To sum up with the latest definitions and their consequences, we can picture the situation as follows, for some given logic **L**:

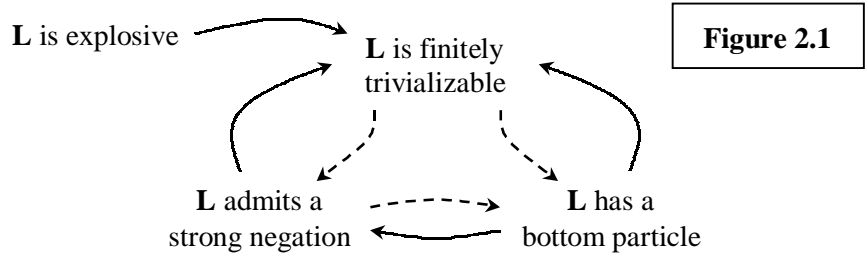

| Figure 2.1 |

Where:
- ✗ ⟶ ✓ means that ✗ entails ✓
- ✗ ⤍ ✓ means that ✗ plus left-adjunctiveness entails ✓

**2.4 Huge tracts of the logical space.** Lo and behold! If now the reader only learns that all properties mentioned in the last subsection *are* compatible with the definition of a paraconsistent logic, they are sure to obtain a wider view of the paraconsistent landscape. Indeed, general non-explosive logics, that is, logics in which not all theories are explosive, can indeed uphold the existence either of finitely trivializable theories, strong negations, or bottom particles! (A rough map of this brave new territory may be found in **Figure 2.2**.) Logics which are paraconsistent but nevertheless have some special explosive theories, such as the ones just mentioned, will constitute the focus of our attention from now on, for, as we shall argue, they may let us explore some fields into which we would not tread in the lack of those theories. Some interesting new concepts can now be studied —this is the case of the notion of *consistency* (and its dual, the notion of *inconsistency*), as we shall argue.

Consider for instance the logic *Pac*, given by the following matrices:

| $\wedge$ | **1** | ½ | **0** |
|---|---|---|---|
| **1** | 1 | ½ | 0 |
| ½ | ½ | ½ | 0 |
| **0** | 0 | 0 | 0 |

| $\vee$ | **1** | ½ | **0** |
|---|---|---|---|
| **1** | 1 | 1 | 1 |
| ½ | 1 | ½ | ½ |
| **0** | 1 | ½ | 0 |

| $\rightarrow$ | **1** | ½ | **0** |
|---|---|---|---|
| **1** | 1 | ½ | 0 |
| ½ | 1 | ½ | 0 |
| **0** | 1 | 1 | 1 |

| | $\neg$ |
|---|---|
| **1** | 0 |
| ½ | ½ |
| **0** | 1 |

where both 1 and ½ are distinguished values. This is the name under which this logic appeared in Avron's [8] (section 3.2.2), though it had previously appeared, for instance, in Avron's [7], under the name $RM_3^\supset$, and, even before than that, in Batens's [10], under the name $PI^s$. It is easy to see that, in such a logic, for no formula $A$ it can be the case that $A, \neg A \vDash_{Pac} B$, for all $B$. So, $Pac$ is a non-explosive, thus paraconsistent, logic. Conjunction, disjunction and implication in $Pac$ are fairly classical connectives: in fact, the whole positive classical logic is validated by its matrices. But the negation in $Pac$ is in some sense strongly non-classical in its surrounding environment, and the immediate consequence of this is that $Pac$ does not have any explosive theory as the ones mentioned above. If such a three-valued logic would define a negation having all properties of classical negation, the table at the right shows how it would look. It is very easy to see that such a negation (in fact, a strong negation with all classical properties) is *not* definable in $Pac$, for any truth-function of this logic having only ½'s as input will also have ½ as output. As a consequence, $Pac$ will not

| | ~ |
|---|---|
| **1** | 0 |
| ½ | 0 |
| **0** | 1 |

respect *ex falso*, having no bottom particle (being unable, thus, as we shall argue before the end of this section, to express the consistency of its formulas), and once it is evidently a left-adjunctive logic as well, it will not even be finitely trivializable at all. One could then criticize such a logic for providing a very weak interpretation for negation, once in this logic all contradictions are admissible. This has some weird consequences and is certainly too light a way of obtaining a paraconsistent logic (this is also the central point of Batens's criticism of Priest's $LP$,[9] see [13]): if some contradictions will give you trouble just assume, then, that no contradiction at all can ever really hurt your logic! Under our present point of view, proposing a logic in which no single contradiction can ever have a harmful effect on their underlying theories is quite an extremist position, and may take us too far away from any classical form of reasoning.[10]

Now, if one endows the language of $Pac$ either with such a strong negation or a *falsum* constant (a bottom particle), with the canonical interpretation, what will result is a well-known conservative extension of it, called $\mathbf{J}_3$, which is still paraconsistent but has all those special explosive theories neglected by $Pac$. This logic $\mathbf{J}_3$ was first introduced by D'Ottaviano and da Costa in 1970 (cf. [60]) as a 'possible solution to the problem of Jaśkowski', and reappeared quite often in the literature after that. The first presentation of $\mathbf{J}_3$ did not bring the strong negation ~ as a primitive connective, but displayed instead a primitive 'possibility connective' $\nabla$ (see its table to the right). In [61] it was once more presented, but this time having also a sort of 'consistency connective' $\circ$ as primitive (table to the right), and in [44] we have explored more deeply the expressive and inferential power of this logic, and the possibility of applying it to the study of inconsistent

| | $\nabla$ | $\circ$ |
|---|---|---|
| **1** | 1 | 1 |
| ½ | 1 | 0 |
| **0** | 0 | 1 |

---

[9] By the way, Priest's logic $LP$ is nothing but the implicationless fragment of $Pac$ (cf. [90]).

[10] This constitutes indeed the kernel of a long controversy between H. Jeffreys and K. Popper. The first author argued in 1938 that contradictions should not be reasonably supposed to imply anything else, to which the second author replied in 1940 saying that contradictions are fatal and should be avoided at all costs, to prevent science from collapsing. Jeffreys aptly reiterated, in 1942, that he was not suggesting that *all* contradictions should be tolerated, but at least *some*. Popper responded to this successively in 1943, 1959 and 1963, saying that he himself had thought about a system in which contradictory sentences were not 'embracing', that is, did not explode, but he abandoned this system because it turned out to be too weak (lacking, for instance, *modus ponens*), and he hastily concluded from that that no useful such a system could ever be attained. See more details and references about this dispute in [26], chapter VI.

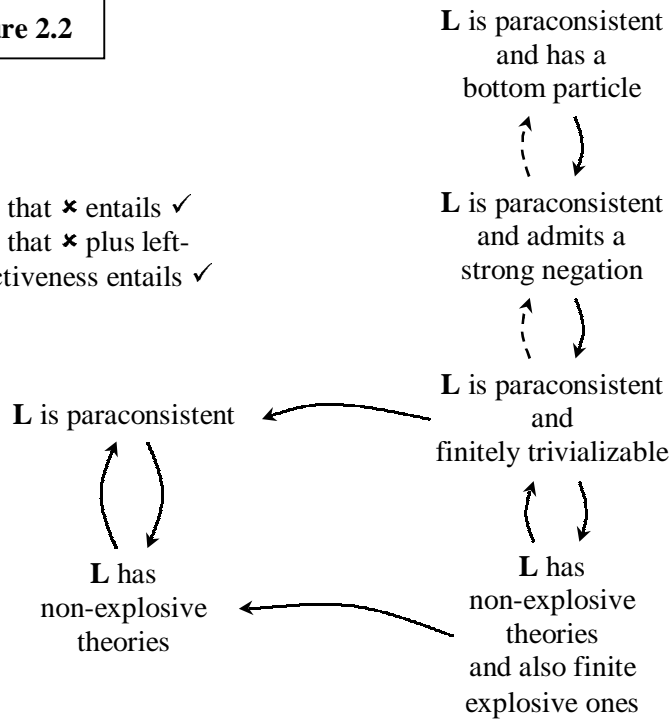

**Figure 2.2**

Where:

✗ ⟶ ✓ means that ✗ entails ✓
✗ ⤍ ✓ means that ✗ plus left-
adjunctiveness entails ✓

**L** is paraconsistent
and has a
bottom particle

**L** is paraconsistent
and admits a
strong negation

**L** is paraconsistent          **L** is paraconsistent
                                    and
                                 finitely trivializable

**L** has                        **L** has
non-explosive                    non-explosive
theories                         theories
                                 and also finite
                                 explosive ones

databases, abandoning ~ and ∇ but still maintaining ◦ as primitive. As a result, we have argued that this logic (now renamed **LFI1**, one of our main 'logics of formal inconsistency') has been shown to be perfectly adequate, among other options, for the task of formalizing the notion of (in)consistency in a very strong and sensible way. But we will have much more to say about this further on.

The reader could now certainly ask himself: If paraconsistency is about non-explosiveness, why are you so interested in having these special explosive theories? Because our interest lies much further than the simple control of the explosive power of contradictions —we want to be able to retain classical reasoning, if only under some suitable interpretation of a fragment of our paraconsistent logics, and we also want to use these paraconsistent logics not only to reason under conditions which do not presuppose consistency, but we want to be able to take hold of the very notion of consistency inside of our logics! From this point of view, the paraconsistent logics which shall interest us are exactly those which permit us to formalize, and get a good grip on, the intricate phenomenon of *inconsistency*, as opposed to mere cut and dried *contradictoriness*.

Whatever inconsistency might mean, by our previous analysis, we might surely suppose a trivial theory to be not only contradictory but inconsistent as well. But yet, a contradiction is certainly one of the many guises of inconsistency! So one may conjecture that *consistency* is exactly what a contradiction might be lacking to become explosive —if it was not explosive from the start. Roughly speaking, we are going to suppose that a 'consistent contradiction' is likely to explode, even if a 'regular' contradiction is not. In logics such as classical logic, consistency is well established, and indeed all theories are explosive; therefore, in any given classical theory, a contradiction turns out to be not only a necessary but also a sufficient condition for triviality.

Now, based on the above considerations, let us suppose in general that a proposition *can* be contradictory and still does not cause much harm, in general, in a paraconsistent logic, if only its consistency is not guaranteed, or cannot be established. Thus, an 'inconsistent' contradiction will be allowed to show up with no big commotion, but still a 'consistent' one should behave classically, and explode! This is how we will put it in formal terms. Let $\Delta(A)$ here denote a (possibly empty) set of schemas depending only on $A$. We will call a theory $\Gamma$ *gently explosive* if:

(a) $\exists A$ such that $\Delta(A) \cup \{A\}$ is not trivial, $\Delta(A) \cup \{\neg A\}$ is not trivial, and
(b) $\forall A\, \forall B\, [\Gamma, \Delta(A), A, \neg A \Vdash B]$.         (D15)

The gently explosive theory $\Gamma$ will be said to be finitely so when $\Delta(A)$ is a finite set, so that a finitely gently explosive theory will be simply one that is finitely trivialized in a very distinctive way.       (D16)

Accordingly, a logic $\mathbf{L}$ will be said to be *[finitely] gently explosive* when all of its theories are [finitely] gently explosive.       [(D17)] (D18)

Thus, in any such a gently explosive logic, given a contradictory theory there is always something 'reasonable' —to wit, consistency— which one can add to it in order to guarantee that it will become trivial. We may now consider the following gentle versions of the Principle of Explosion:

#### [FINITE] GENTLE PRINCIPLE OF EXPLOSION       [(fgPPS)] (gPPS)
   $\mathbf{L}$ must be [finitely] gently explosive.

So, according to the interpretation proposed above, what we are implicitly assuming in the above principles is that, for any given formula $A$, the (finite) set $\Delta(A)$ will express, in a certain sense, the *consistency* of $A$ relative to the logic $\mathbf{L}$.

Based on that, we may define the consistency of a logic in the following way. $\mathbf{L}$ will be said to be *consistent* if:

  (a) $\mathbf{L}$ is gently explosive, and (b) $\forall A\, [\Delta(A) = \varnothing$ or $\forall\Gamma\, (\forall B \in \Delta(A))(\Gamma \Vdash B)]$.  (D19)

It immediately follows, from these definitions and the preceding ones, that:

**FACT 2.14** (i) Any non-trivial explosive theory / logic is finitely gently explosive. (ii) Any transitive logic is consistent if, and only if, it is both explosive and non-trivial. (iii) Any transitive consistent logic is finitely gently explosive. (iv) Any left-adjunctive finitely gently explosive logic is supplementing explosive.

$$[\text{not-(D2) and (D3)}] \Rightarrow (\text{D16}); [(\text{PNT}) \text{ and } (\text{PPS})] \Rightarrow (\text{fgPPS});$$
$$(\text{Con3}) \Rightarrow \{(\text{D19}) \Leftrightarrow [(\text{D6}) \text{ and not-(D5)}]\};$$
$$[(\text{Con3}) \text{ and } (\text{D19})] \Rightarrow (\text{D17});$$
$$[(\text{D13}) \text{ and } (\text{fgPPS})] \Rightarrow (\text{sPPS})$$

**Proof:** To check (i), just let $\Delta(A)$ be empty, for every formula $A$. This result evidently parallels FACT 2.10(iii), about supplementing explosive logics. To see, in (ii), that any given consistent logic is explosive use transitivity whenever you meet a non-empty $\Delta$. Part (iii) follows from (i) and (ii), and part (iv) simply reflects FACT 2.13(ii).   □

So, based on the above definition of a consistent logic and the subsequent fact, if we were to define a so-called *Principle of Consistency*, it would then simply coincide with the sum of (PNT) and (PPS), for logics obeying transitivity. We shall, therefore, not insist in explicitly formulating here such a principle.

We may now finally define what we will mean by a *logic of formal inconsistency* (**LFI**), which will be nothing more than a logic that allows us to 'talk about consis-tency' in a meaningful way. We will consider, of course, an *inconsistent* logic to be simply one that is not consistent. This assumption, together with FACT 2.14(ii), explains why paraconsistent logics were early dubbed, by da Costa, 'inconsistent formal systems', once all paraconsistent logics are certainly inconsistent in the sense of not respecting (D19), even though they are always also non-trivial and quite often they are non-contradictory as well. Those inconsistent logics which went so far as to be trivial, and thus no more paraconsistent at all, were dubbed, by Miró Quesada, *absolutely inconsistent* logics (cf. [80]). Now, an **LFI** will be any non-trivial logic in which consistency does not hold, but can still be expressed, thus being a gently explosive and yet non-explosive logic, that is, a logic in which:

$$\text{(a) (PPS) does not hold, but (b) (gPPS) holds.} \tag{D20}$$

Classical logic, then, will not be an **LFI** just because (PPS) holds in it. *Pac* will also not be an **LFI**, even though it is paraconsistent, for *Pac* is not finitely trivializable. But D'Ottaviano & da Costa's $\mathbf{J}_3$ (and, consequently, our **LFI1**), which conservatively extends *Pac*, will *indeed* be an **LFI**, where consistency is expressed by the connective ∘ (see above), and inconsistency, as usual, is expressed, by the negation of this connective. Also, Jaśkowski's **D2** will constitute an **LFI**, as the reader can easily check, where the consistency of a formula $A$ can be expressed by the formula ( $A \vee \neg A$), written in terms of the necessity operator of *S5*.[11]

'Only' **LFI**s —though these seem to comprise by far the *great majority* of all known paraconsistent logics— will interest us in this study.

**2.5 DEFCON 2: one step short of trivialization.** The distinction between the original formulation of explosiveness, its formulation in terms of *ex falso*, and its supplementing and gentle formulations offered above does not tell you everything you need to know about the ways of exploding. Indeed, there are more things in the realm of explosiveness, dear reader, than are dreamt of in your philosophy! Thus, for instance, a not very interesting scenario seems to unfold if contradictions are still prevented from rendering a given theory trivial but nevertheless are allowed to go half the way, causing some kind of 'partial trivialization'. So, a theory $\Gamma$ will be said to be *partially trivial with respect to* a given schema $\sigma(C_1, …, C_n)$, or $\sigma$-*partially trivial*, if:

$$\text{(a) } \exists C_1 … \exists C_n \text{ such that } \sigma(C_1, …, C_n) \text{ is not a top particle, and}$$
$$\text{(b) } \forall C_1 … \forall C_n \ [\Gamma \Vdash \sigma(C_1, …, C_n)]. \tag{D21}$$

Following this same path, a theory $\Gamma$ will be said to be *partially explosive with respect to* the schema $\sigma(C_1, …, C_n)$, or $\sigma$-*partially explosive*, if:

$$\text{(a) } \exists C_1 … \exists C_n \text{ such that } \sigma(C_1, …, C_n) \text{ is not a top particle, and}$$
$$\text{(b) } \forall C_1 … \forall C_n \forall A \ [\Gamma, A, \neg A \Vdash \sigma(C_1, …, C_n)]. \tag{D22}$$

Of course, a logic **L** will be said to be $\sigma$-partially trivial / $\sigma$-partially explosive if all of its theories are $\sigma$-partially trivial / $\sigma$-partially explosive. respect. (D23), (D24)

More simply, a theory, or a logic, can now be said to be *partially trivial / partially explosive* if this theory, or logic, is $\sigma$-partially trivial / $\sigma$-partially explosive, for some

---

[11] This needs to be qualified. Among the various formulations among which **D2** has appeared in the literature, it is not completely clear if its language has a necessity operator available so as to make this definition possible, or not. If this is not available, it may well be that **D2** is not characterizable as an **LFI** after all (even though a situation for a necessity operator would quite naturally appear, to all practical purposes, in the trivial case in which there is just one person 'discussing', or even more unlikely, a situation in which all contenders just agree with each other).

schema σ. We can now immediately formulate the following new version of the Principle of Explosion:

PRINCIPLE OF PARTIAL EXPLOSION                                    (pPPS)
  **L** must be partially explosive.

One may immediately conclude that:

FACT 2.15  (i) Any partially trivial theory / logic is partially explosive.  (ii) Any explosive logic is partially explosive.        (D21) ⇒ (D22); (D23) ⇒ (D24); (PPS) ⇒ (pPPS)

A well-known example of a logic which is not explosive but is partially explosive even so, is given by Kolmogorov & Johánsson's Minimal Intuitionistic Logic, **MIL**, which is obtained by the addition to the positive part of intuitionistic logic of some forms of *reductio ad absurdum* (cf. [68] and [69]). What happens, in this logic, is that $\forall \Gamma \forall A \forall B$ ($\Gamma, A, \neg A \Vdash B$) is *not* the case, but still it *does* hold that $\forall \Gamma \forall A \forall B$ ($\Gamma, A, \neg A \Vdash \neg B$). This means that **MIL** is paraconsistent in a broad sense, for contradictions do not explode, but still all *negated* propositions can be inferred from any given contradiction!

It is something of a consensus that an interesting paraconsistent logic should not only avoid triviality but also partial triviality. Thus, the following definition now comes in handy.  A logic **L** will be said to be *boldly paraconsistent* if:

$$\text{(pPPS) fails for } \mathbf{L}. \tag{BPL}$$

Evidently:

FACT 2.16  A boldly paraconsistent logic is paraconsistent.        (BPL) ⇒ (PL2)

Now, let's tackle a somewhat inverse approach.  Call a theory $\Gamma$ *controllably explosive in contact with* a given schema $\sigma(C_1, \ldots, C_m)$ if:

(a) $\exists C_1 \ldots \exists C_m$ such that $\sigma(C_1, \ldots, C_m)$ and $\neg\sigma(C_1, \ldots, C_m)$ are not bottom particles, and (b) $\forall C_1 \ldots \forall C_n \forall B$ [$\Gamma, \sigma(C_1, \ldots, C_m), \neg\sigma(C_1, \ldots, C_m) \Vdash B$].        (D25)

Accordingly, a logic **L** will be said to be *controllably explosive in contact with* $\sigma(C_1, \ldots, C_m)$ when all of its theories are controllably explosive in contact with this schema.        (D26)

Some given theory / logic can now more simply be called *controllably explosive* when this theory / logic has some schema in contact with which it is controllably explosive. An immediate new version of the Principle of Explosion that suggests itself then is:

CONTROLLABLE PRINCIPLE OF EXPLOSION                              (cPPS)
  **L** must be controllably explosive.

Similarly to the case of FACT 2.14, parts (i) and (iii), it follows here that:

FACT 2.17  (i) Any non-trivial explosive theory / logic is controllably explosive.  (ii) Any transitive consistent logic is controllably explosive.        [not-(D2) and (D3)] ⇒ (D25)
                    [(PNT) and (PPS)] ⇒ (cPPS); [(Con3) and (D19)] ⇒ (D26)

By the way, we may also now emend FACT 2.9 so as to immediately conclude that:

FACT 2.18  Any finitely-gently / controllably explosive logic is finitely trivializable, and yet non-trivial.        [(D17) or (D26)] ⇒ [(D7) and not-(D5)]

This fact can be used to update and complement the information conveyed in **Figure 2.1**.

Now, there seems to be no good reason to rule out controllably explosive theories, as we did in the case of partially explosive theories by way of the bold definition of paraconsistency, (BPL). In fact, it seems that most, if not all, finitely gently explosive logics *are* controllably explosive, and vice-versa! We will see, later on, many examples of paraconsistent logics —indeed, of **LFI**s— which not only are obviously gently explosive, but are also controllably explosive in contact with schemas such as $(A \wedge \neg A)$, or such as $\circ A$, where $\circ$, we recall, is a connective expressing consistency (Jaśkowski's **D2**, for instance, may already be one of these, but the logic **LFI1**, on the other hand, explodes only in contact with the second of these schemas). There are even logics which controllably explode in contact with large classes of non-atomic propositions (see [78], and ahead, for a number of them). An extreme case of these, as we shall see, is given by Sette's three-valued logic $\mathbf{P}^1$ (cf. [103]), which controllably explodes in contact with *any* complex formula, and so can be said to behave paraconsistently only at the level of its atoms. It is also not uncommon for some paraconsistent logic **L** having a strong negation $\div$ to be controllably explosive. In fact, it suffices that such a logic is transitive and infers $\neg \div A \Vdash A$, and of course it will turn out to be controllably explosive in contact with $\div A$, or at least in contact with $\div \div A$ (see, for instance, FACT 3.76, or THEOREM 3.51(i) and FACT 3.66). Many **LFI**s will moreover be controllably explosive in contact with any consistent formula (see FACT 3.32). And so on, and so forth.

A range of variations on the above versions of the Principle of Explosion can be obtained if we only mix the ones we already have. We shall nevertheless not investigate this theme here any further, but only notice that the multiple relations, hinted above, between (sPPS), (gPPS) and (cPPS), the supplementing, the gentle and the controllable forms of explosion, certainly deserve a closer and more attentive look by the 'paraconsistent community' and sympathizers.

**2.6 C-systems.** Given a logic $\mathbf{L} = <For, \Vdash>$, let $For^+ \subseteq For$ denote the set of all *positive formulas* of **L**, that is, the *negationless* fragment of $For$, or, in still other words, the set of all formulas in which no negation symbol $\neg$ occurs. The logic $\mathbf{L}1 = <For_1, \Vdash_1>$ is said to be *positively preserving relative to* the logic $\mathbf{L}2 = <For_2, \Vdash_2>$ if:

$$\text{(a) } For_1^+ = For_2^+, \text{ and (b) } (\Gamma \Vdash_1 A \iff \Gamma \Vdash_2 A), \text{ for all } \Gamma \cup \{A\} \subseteq For_1^+. \quad \text{(D27)}$$

So, if **L**1 is positively preserving relative to **L**2, then it will in general be a conservative extension of the positive fragment of **L**2. Now, as an example of the ubiquity of **LFI**s inside the realm of paraconsistent logics, just notice that:

FACT 2.19 Any paraconsistent logic that is positively preserving relative to classical logic and has a bottom particle can be characterized as an **LFI**.

**Proof:** Just define $\circ A$ as $(A \to \bot) \vee (\neg A \to \bot)$, and check that, in general, $\circ A$ is not a top particle, $\{\circ A, A\}$ is not always trivial, and $\{\circ A, \neg A\}$ is not always trivial, but that, in any case, $\{\circ A, A, \neg A\}$ is indeed a trivial theory. This result actually holds for any logic having a *left-adjunctive disjunction*, that is, a binary connective $\vee$ such that $(B \vee C)$ is not a bottom particle, for some formulas $B$ and $C$, and such that $\forall B \forall C \forall \Gamma \forall \Delta \forall D \{(\Gamma, B \Vdash D) \text{ and } (\Delta, C \Vdash D)] \Rightarrow [\Gamma, \Delta, (B \vee C) \Vdash D]\}$ (for a particular consequence of this feature, see FACT 3.7), and having *modus ponens*: $\forall \Gamma \forall A \forall B [\Gamma, A, (A \to B) \vdash_{min} B]$. You just have to choose $\Gamma = \{A\}$, $B = (A \to \bot)$, $\Delta = \{\neg A\}$, and $C = (\neg A \to \bot)$, and notice that, in this case, both $(\Gamma, B \Vdash \bot)$ and $(\Delta, C \Vdash \bot)$, by *modus ponens*. □

This last result shows that any paraconsistent logic conservatively extending the positive classical logic and respecting either one of the principles of *ex falso* or of supplementing explosion will be finitely gently explosive as well, throwing some light on some hitherto unsuspected connections between (ExF), (sPPS) and (fgPPS), and consequently any such a logic can be easily recast as an **LFI** (take another look at **Figure 2.2**). Consequently, for all such logics, it amounts to be more or less the same starting either with a consistency operator, or with a strong negation, or with a bottom particle: each of these can be used to define the others. This does not mean, however, that 'only' such logics are **LFI**s (see the case of $C_{Lim}$, in the subsection **3.10**).

To specialize a little bit from this very broad definition of **LFI**s above we will now define the concept of a **C**-system. The logic **L1** will be said to be a **C**-*system based on* **L2** if:

(a) **L1** is an **LFI** in which consistency or inconsistency
are expressed by operators (at the object language level),
(b) **L2** is not paraconsistent, and
(c) **L1** is positively preserving relative to **L2**.                    (D28)

Any logic constructed as a **C**-system based on some other logic will more generally be identified simply as a **C**-*system*. In the next section we will study various logics which are **C**-systems, and pinpoint some which are not.

Jaśkowski's **D2**, as we have already seen in the above subsections, *is* an **LFI** and *can* define an operator expressing consistency —at least under some presentations (see note 11). But, in order for it to be characterized as a **C**-system it would still have to be clarified on which logic it is based, that is, where does its peculiar positive (non-adjunctive) part come from! This same question arises with respect to all other logics that are left-adjunctive but not left-disadjunctive, as well as with respect to many relevance logics.

All **C**-systems we will be studying below are inconsistent, non-contradictory and non-trivial. Furthermore, they are boldly paraconsistent (though the proof of *this* fact will be left for [42]), and often controllably explosive as well, they have strong negations and bottom particles, and are positively preserving relative to classical propositional logic —so, that they will respect (PNC), (PNT), (ExF), (sPPS), (gPPS) and often (cPPS), but they will not respect neither (PPS) nor (pPPS). Let's now jump to them.

## 3   COOKING THE **C**-SYSTEMS ON A LOW FLAME

> Indeed, even at this stage I predict a time when there will be mathematical
> investigations of calculi containing contradictions, and people will actu-
> ally be proud of having emancipated themselves even from consistency.
> —Wittgenstein, Philosophical Remarks, p.332.

Underlying the original approach of da Costa to the concoction of a propositional calculus capable of admitting contradictions, yet remaining sensible to performing reasonable deductions, laid the idea of maintaining the positive fragment of classical logic unaltered. This explains why his approach to paraconsistency has eventually received the inelegant label of 'positive (logic) plus approach' and, more recently, the not much descriptive (and in some cases plainly inadequate) label of 'non-truth-functional approach' (cf., respectively, [92] and [94]). Surely, competitive approaches do exist, like the one stemming from Jaśkowski's or Rescher & Brandom's investi-gations, which rejects left-disadjunction, and is usually referred to as a 'non-adjunctive approach' (cf. [67] and [98]), and which has more recently been tentatively dubbed, by J. Perzanowski, as 'parainconsistent logic'.[12] Another megatrend comes from the 'relevance approach' to paraconsistency, captained by the American-Australian school, whose concern is not so much with negation as with implication, giving rise to 'relevance logics' (cf. [3]). Still another very interesting proposal came from Belgium, under the appellation of 'adaptive logics', which do not worry so much about proving consistency, but assume it instead from the very start, as some kind of default (cf. [11] and [12]). Now, let us make it crystal clear that our concentration in this study on the investigation of **C**-systems, born from the first approach mentioned above, wishes not to diminish the other approaches, nor affirms that they should be held as mutually exclusive. Our intention, indeed, is but to present the **C**-systems under a more general and suggestive background, and from now on we shall draw on the other approaches only when we feel that as a really necessary or instructive step. To the reader particularly interested in them, we prefer simply to redirect them to the competent sources.

### 3.1 Paleontology of C-systems.

All definitions and remarks made above were set forth directing an arbitrary consequence relation $\Vdash$, be it syntactical, semantical or defined in any other mind-boggling way. Once the surfacing of contradictions on a theory involves negation, and nothing but that, it is appealing to consider and explore the intuitive idea that an interesting class of paraconsistent logics is to be given by the ones which are positively preserving relative to classical logic, differing from classical logic only in the behavior of formulas involving negation. This is the idea into which we will henceforth be digging, by axiomatically proposing a series of logics characterized by their syntactical consequence relations, $\vdash$, and containing all rules and schemas which hold in the positive part of classical logic. Thus, let's initially consider $\wedge$, $\vee$, $\rightarrow$, and $\neg$ to be our primitive connectives, and consider the set of formulas *For*, as usual, to be the free algebra generated by these connectives. We will start our journey from the following set of axioms:

(Min1)   $\vdash_{min} (A \rightarrow (B \rightarrow A))$;
(Min2)   $\vdash_{min} ((A \rightarrow B) \rightarrow ((A \rightarrow (B \rightarrow C)) \rightarrow (A \rightarrow C)))$;
(Min3)   $\vdash_{min} (A \rightarrow (B \rightarrow (A \wedge B)))$;
(Min4)   $\vdash_{min} ((A \wedge B) \rightarrow A)$;
(Min5)   $\vdash_{min} ((A \wedge B) \rightarrow B)$;
(Min6)   $\vdash_{min} (A \rightarrow (A \vee B))$;
(Min7)   $\vdash_{min} (B \rightarrow (A \vee B))$;
(Min8)   $\vdash_{min} ((A \rightarrow C) \rightarrow ((B \rightarrow C) \rightarrow ((A \vee B) \rightarrow C)))$;
(Min9)   $\vdash_{min} (A \vee (A \rightarrow B))$;
(Min10)  $\vdash_{min} (A \vee \neg A)$;
(Min11)  $\vdash_{min} (\neg\neg A \rightarrow A)$.

Here, by writing $\vdash_{min} (A \rightarrow (B \rightarrow A))$ we will be abbreviatedly denoting that:

$$\forall \Gamma \, \forall A \, \forall B \, [\Gamma \vdash_{min} (A \rightarrow (B \rightarrow A))],$$

---

and so on, for the other axioms. The only inference rule, as usual, will be *modus ponens*, (MP): $\forall \Gamma \forall A \forall B [\Gamma, A, (A \to B) \vdash_{min} B]$. The logic built using such axioms, plus (MP) and the usual notion of proof from premises (we may now be calling *proofs*, *theorems* and *premises* which we have previously called, respectively, inferences, theses and theories) was called $C_{min} = \langle For, \vdash_{min} \rangle$ and studied by the authors in [39].

First of all, let us observe that the so-called *Deduction Metatheorem* is here valid:

**THEOREM 3.1** $[\Gamma, A \vdash_{min} B \Rightarrow \Gamma \vdash_{min} (A \to B)]$.[13]
**Proof:** It is a familiar and straightforward procedure to show that the Deduction Metatheorem holds for any logic containing (Min1) and (Min2) as provable schemas and having only *modus ponens* as a primitive rule.

Evidently, by monotonicity and transitivity, *modus ponens* already gives us the converse of THEOREM 3.1. This makes it possible for us to introduce all axioms as some sort of axiomatic inference rules, and this is what we shall do from now on. Moreover, using the Deduction Metatheorem and its converse, one could now equivalently represent, in $C_{min}$, the fact that (PPS) (the Principle of Explosion) does not hold by the unprovability of the theorem (tPS): $(A \to (\neg A \to B))$. And indeed:

**THEOREM 3.2** (tPS) is not provable by $C_{min}$.
**Proof:** Use the matrices of *Pac*, in the subsection **2.4**, to check that all axioms above are validated and that (MP) preserves validity, while (tPS) is not always validated. This shows that $C_{min}$ is a fragment of *Pac*, and so it also cannot prove (tPS). In fact, (tPS) is more than non-provable, it is *independent* from $C_{min}$ (and *Pac*) for its negation is not even classically provable, and *Pac* is a deductive fragment of classical logic.

As usual, *bi-implication*, $\leftrightarrow$, will be defined by setting $(A \leftrightarrow B) \overset{\text{def}}{=} ((A \to B) \land (B \to A))$. Note that, by the above considerations, $\vdash_{min} (A \leftrightarrow B)$ if, and only if, $A \vdash_{min} B$, and $B \vdash_{min} A$, which is the same as writing $A \dashv\vdash_{min} B$. So, bi-implication holds between two formulas if, and only if, they are (provably) equivalent (see (Eq1), in the subsection **2.2**). Nevertheless, as the reader shall see below, having two equivalent formulas, in the logics we will be studying here, usually does *not* mean, as in classical logic, that these formulas can be freely intersubstituted everywhere (take a look, ahead, for instance, at results 3.22, 3.35, 3.51, 3.58, 3.65, and 3.74).

Axioms (Min1)–(Min8) are known at least since Gentzen's [62] as providing an axiomatization for the so-called 'positive logic'. Of course, they immediately tell us, among other things, that the conjunction of this logic is both left-adjunctive and left-disadjunctive (just take a look at axioms (Min3)–(Min5)). Nevertheless, (Min9): $(A \lor (A \to B))$, which *is* a positive schema, is *not* provable even if one uses (Min10) and (Min11) in addition to (Min1)–(Min8) and (MP) (i.e. the logic axiomatized as $C_{min}$ minus the axiom (Min9))! Indeed:

**THEOREM 3.3** (Min9) is not provable by $C_{min} \backslash \{(Min9)\}$.
**Proof:** Use the following matrices (cf. [2]) to check that (Min9) is independent from $C_{min} \backslash \{(Min9)\}$:

---

[13] Read this kind of sentence as a universally quantified one —in this case, for example, it would be $\forall \Gamma \forall A \forall B (\Gamma, A \vdash_{min} B \Rightarrow \Gamma \vdash_{min} (A \to B))$

| ∧ | 1 | ½ | 0 | | ∨ | 1 | ½ | 0 | | → | 1 | ½ | 0 | | ¬ | |
|---|---|---|---|---|---|---|---|---|---|---|---|---|---|---|---|
| **1** | 1 | ½ | 0 | | **1** | 1 | 1 | 1 | | **1** | 1 | ½ | 0 | | **1** | 0 |
| **½** | ½ | ½ | 0 | | **½** | 1 | ½ | ½ | | **½** | 1 | 1 | 0 | | **½** | 1 |
| **0** | 0 | 0 | 0 | | **0** | 1 | ½ | 0 | | **0** | 1 | 1 | 1 | | **0** | 1 |

where 1 is the only distinguished value.

So, what is this thing that Gentzen (and Hilbert before him) have dubbed 'positive logic', if even a deductive extension of it is unable to prove all positive theorems of classical logic? Here is the trick: Gentzen referred of course to positive *intuitionistic* logic, and not to the classical logic! So, this logic $C_{min}\backslash\{(\text{Min9})\}$, which was proposed by da Costa (cf. [49]) and called $C_\omega$ by him, turns out to be only positively preserving relative to intuitionistic logic, and not relative to classical logic. In [39] we have proven that its deductive extension $C_{min}$, obtained by adjoining (Min9) to $C_\omega$, is indeed positively preserving relative to classical logic, and moreover:

**THEOREM 3.4** $C_{min}$ does have neither a strong negation nor a bottom particle, and is not finitely trivializable.

**Proof:** PROPOSITION 2.5, in [39], shows that $C_{min}$ does not have a bottom particle, and so, by left-disadjunction and FACT 2.13, it does not have neither a strong negation nor is it finitely trivializable.

Moreover, in [39] we also proved that:

**THEOREM 3.5** $C_{min}$ does not have any negated theorem, i.e. ($\nvdash_{min} \neg A$).

Of course, both results above are valid, *a fortiori*, for $C_\omega$. Indeed, as shown by Urbas (cf. [107]), these logics are very weak with respect to negation, so that the following holds:

**THEOREM 3.6** No two different negated formulas of $C_{min}$ are provably equivalent.

**Proof:** The THEOREM 2, in [107], shows that $\neg A \dashv\vdash \neg B$ is derivable in $C_\omega$ if and only if $A$ and $B$ are the same formula. It is straightforward to adapt this result also to $C_{min}$.

Much more about the provability (or validity) of negated theorems will be seen in the paper [42], which brings semantics to most logics here studied.

THEOREM 3.4 shows that $C_{min}$, or $C_\omega$, *cannot* be **C**-systems based on classical logic, or intuitionistic logic, once they are both *compact* (all proofs are finite) and not finitely gently explosive, so that they cannot be gently explosive at all, and thus cannot formalize 'consistency', in the precise sense formulated in the subsection **2.4**. We had better then make them deductively stronger in order to get what we want.

We make a few more important remarks before closing this subsection. First, note that (Min10): $(A \lor \neg A)$ was added in order to keep $C_{min}$ and $C_\omega$ from being *paracomplete* as well as paraconsistent (let's investigate one deviancy at a time!), and this axiom can indeed be pretty useful in providing us with a form of *proof by cases*:

**FACT 3.7** $[(\Gamma, A \vdash_{min} B)$ and $(\Delta, \neg A \vdash_{min} B)] \Rightarrow (\Gamma, \Delta \vdash_{min} B)$.

**Proof:** From (Min8) and (Min10), by *modus ponens*, monotonicity, and the Deduction Metatheorem (from now on, we will not mention these last three every time we use them anymore).

It will also be practical here and there to use $[(A \rightarrow B), (B \rightarrow C) \vdash_{min} (A \rightarrow C)]$ (a kind of logical version for the transitivity property) as an alternative form of the axiom (Min2). Indeed:

FACT 3.8 (Min2) can be substituted, in $C_{min}$, by $[(A \rightarrow B), (B \rightarrow C) \vdash_{min} (A \rightarrow C)]$.

We shall often make use of both these forms without discriminating which.

In the next subsection (see THEOREM 3.13) we will learn about the utility of (Min11): $(\neg\neg A \rightarrow A)$, which was added by da Costa as a way of rendering the negation of his calculi a bit stronger, using as an argument the intended duality with the logics arising from the formalization of intuitionistic logic, in which usually only the converse of (Min11), i.e. the formula $(A \rightarrow \neg\neg A)$, is valid.

It is quite interesting as well to notice that the addition of the 'Theorem of Pseudo-Scotus', (tPS), to $C_{min}$ as a new axiom schema will not only prevent the resulting logic from being paraconsistent, but it will also provide a complete axiomatization for the *classical propositional logic* (hereby denoted **CPL**). In fact, it is a well-known fact that:

THEOREM 3.9 Axioms (Min1)–(Min11) plus (tPS): $(A \rightarrow (\neg A \rightarrow B))$, and (MP), provide a sound and complete axiomatization for **CPL**.

Actually, the axiom (Min11) can be discharged from the above axiomatization, being proved from the remaining ones. Axiom (Min9) also turns to be redundant (take a look at the FACT 3.45, below).

**3.2 The basic logic of (in)consistency.** Let's consider an extension of our language by the addition of a new unary connective, $\circ$, representing *consistency*. Let's now also add, to $C_{min}$, a new rule, realizing the Finite Gentle Principle of Explosion:

(bc1)     $\circ A, A, \neg A \vdash_{\mathbf{bC}} B$.     'If $A$ is consistent and contradictory, then it explodes'

We will call this new logic, characterized by axioms (Min1)–(Min11) and (bc1), plus (MP), the *basic logic of (in)consistency*, or **bC**. Clearly, thanks to (bc1), we know that **bC** is indeed an **LFI**, i.e. a *logic of formal inconsistency*, and so it is in fact a **C**-*system* based on **CPL**. A strong negation, $\sim$, for a formula $A$ can now be easily defined by setting $\sim A \overset{\text{def}}{=} (\neg A \wedge \circ A)$, and evidently we will have $[A, \sim A \vdash_{\mathbf{bC}} B]$, as expected. A *bottom particle*, of course, is given by $(A \wedge \sim A)$, for any $A$. For alternative ways of formulating **bC**, consider FACT 2.19 and the comments which follow it.

We can already show that THEOREMS 3.5 and 3.6 do *not* hold for **bC**:

THEOREM 3.10 **bC** does have negated theorems, and equivalent negated formulas (but, on the other hand, it has no consistent theorems, that is, theorems of the form $\circ A$).
**Proof:** Consider any bottom particle $\perp$ of **bC**. By definition, it must be such that $(\perp \vdash_{\mathbf{bC}} B)$, for any formula $B$, and so, in particular, $(\perp \vdash_{\mathbf{bC}} \neg\perp)$. But we also have that $(\neg\perp \vdash_{\mathbf{bC}} \neg\perp)$, and proof by cases (FACT 3.7) tells us then that $\vdash_{\mathbf{bC}} \neg\perp$. By the way, this result also transforms THEOREM 3.4 into a corollary of THEOREM 3.5 —if a reflexive logic has proof by cases and no negated theorems, then it cannot contain a bottom particle. Evidently, any bottom particle is equivalent to any other. To check that no formula of the form $\circ A$ is provable, one may just use the classical matrices for $\wedge$, $\vee$, $\rightarrow$ and $\neg$, and pick for $\circ$ a matrix with value constant and equal to 0.

Now, it is easy to see that, in such logic **bC**, if the consistency of the right formulas is guaranteed, than its inferences will behave exactly like in **CPL**. Indeed:

**THEOREM 3.11**  $[\Gamma \vdash_{\mathbf{CPL}} A] \Leftrightarrow [\circ(\Delta), \Gamma \vdash_{\mathbf{bC}} A]$, where $\circ(\Delta) = \{\circ B : B \in \Delta\}$, and $\Delta$ is a finite set of formulas.

**Proof:**  On the one hand, one may just reproduce line by line a **CPL** proof in **bC**, and when it comes to an application of (tPS) —see an axiomatization of **CPL** in the THEOREM 3.9— one will have to use (bc1) instead, and add as a further assumption the consistency of the formula in the antecedent. The converse is immediate.

We know that **bC** is a both a linguistic and a deductive extension of $C_{min}$, once it not only introduces a new connective but has an axiomatic rule telling you what to do with it. But we know more than that:

**THEOREM 3.12  bC** is a conservative extension of $C_{min}$.

**Proof:**  Indeed, if you consider the **bC**-inferences in the language of $C_{min}$, you can no more use (bc1) along a proof, and so you can prove nothing more than you could prove before.

What we have then, in (bc1), is a sort of rough logical clone for the finite gentle rule of explosion. Now, da Costa, in the original presentation of his calculi, which guides us here, has never used a gentle form of explosion but used instead a gentle form of *reductio ad absurdum*:

(RA0)    $\circ B, (A \rightarrow B), (A \rightarrow \neg B) \vdash \neg A$.
     'If supposing $A$ will bring us to a consistent contradiction, then $\neg A$ should be the case'

Notice, by the way, that $((A \rightarrow B) \rightarrow ((A \rightarrow \neg B) \rightarrow \neg A))$ was exactly the form of *reductio* used by Kolmogorov and Johánsson in the proposal of their Minimal Intuitionistic Logic, mentioned above as an example of a logic which is paraconsistent and yet not boldly paraconsistent. Now, the reader might be suspecting that it would really make no difference whether we used (bc1) or (RA0) in the characterization of **bC**. They are right, but this assertion could be made more precise. Indeed, consider the two following alternative versions of these rules:

(bc0)    $\circ A, A, \neg A \vdash \neg B$;
    'If $A$ is consistent and contradictory, then it partially explodes with respect to negated propositions'
(RA1)    $\circ B, (\neg A \rightarrow B), (\neg A \rightarrow \neg B) \vdash A$,
    'If supposing $\neg A$ will bring us to a consistent contradiction, then $A$ should be the case'

and consider the logic *PI* (that is how it was called when it appeared in [10]), characterized simply by (Min1)–(Min10) plus (MP), that is, $C_{min}$ deprived of the schema (Min11): $(\neg\neg A \rightarrow A)$. Then we can prove that:

**THEOREM 3.13  (i)** It does *not* have the same effect adding either (bc1) or (RA0) to *PI*. **(ii)** It *does* have the same effect adding to *PI*: a) (bc0) or (RA0); b) (bc1) or (RA1). **(iii)** It *does* have the same effect adding to $C_{min}$ whichever of the schemas (bc0), (bc1), (RA0) or (RA1). **(iv) bC** cannot be extended into a $\neg$-partially explosive paraconsistent logic.

**Proof:**  To check part (i), use the classical matrices (with values 1 and 0) for $\wedge$, $\vee$ and $\rightarrow$, but let both $\neg$ and $\circ$ have matrices constant and equal to 1 —this way you will see that (bc1) is not provable by the logic obtained from the addition of (RA0) to *PI*. Part (ii) is easy: use FACT 3.8 to prove (bc0) in *PI* plus (RA0), and to prove (bc1) in *PI* plus (RA1); use (Min1) and the proof by cases to prove (RA0) in *PI* plus (bc0), and to prove (RA1) in *PI* plus (bc1). We leave part (iii) as an even easier exercise to the reader (*hint*: use (Min11)). (iv) is an immediate consequence of (iii).

So, this last result gives one reason for us to have our study started from $C_{min}$ rather than from *PI*: we will be avoiding that paraconsistent extensions of our initial logic might turn out to be partially explosive with respect to negated propositions in general, as what occurred with **MIL**, the Minimal Intuitionistic Logic (recall the subsection **2.5**). This feature will help in making many results below more symmetrical. But, to be sure, this does not guarantee that all such extensions will be boldly paraconsistent as well!

The reader should notice that there are, however, some restricted forms of 'reasoning by absurdum' left in **bC**. For example:

FACT 3.14 The following *reductio* deduction rules hold in **bC**:
  (i) $[(\Gamma \vdash_{\textbf{bC}} \circ A)$ and $(\Delta, B \vdash_{\textbf{bC}} A)$ and $(\Lambda, B \vdash_{\textbf{bC}} \neg A)] \Rightarrow (\Gamma, \Delta, \Lambda \vdash_{\textbf{bC}} \neg B)$;
  (ii) $[(\Gamma, B \vdash_{\textbf{bC}} \circ A)$ and $(\Delta, B \vdash_{\textbf{bC}} A)$ and $(\Lambda, B \vdash_{\textbf{bC}} \neg A)] \Rightarrow (\Gamma, \Delta, \Lambda \vdash_{\textbf{bC}} \neg B)$;
  (iii) $[(\Gamma, \neg B \vdash_{\textbf{bC}} \circ A)$ and $(\Delta, \neg B \vdash_{\textbf{bC}} A)$ and $(\Lambda, \neg B \vdash_{\textbf{bC}} \neg A)] \Rightarrow (\Gamma, \Delta, \Lambda \vdash_{\textbf{bC}} B)$.
**Proof:** Part (i) comes immediately from (RA0), part (ii) comes from part (i) using reflexivity and proof by cases, part (iii) comes as a variation of (ii), if you use (Min11).

But we still have not mentioned some of the most decisive features of **bC**! We are now ready for this. Consider, to start with, the following result:

THEOREM 3.15 (i) $(A \wedge \neg A)$ is not a bottom particle in any paraconsistent extension of **bC**. (ii) $\neg(A \wedge \neg A)$ and $\neg(\neg A \wedge A)$ are not top particles in **bC**.
**Proof:** For part (i), just use left-disjunction and THEOREM 3.2 (but the reader might recall from the subsection **2.3** that this formula *is* a bottom particle in some non-left-disjunctive paraconsistent logics such as Jaśkowski's **D2**). To check part (ii) use the following matrices to confirm that neither $\neg(A \wedge \neg A)$ nor $\neg(\neg A \wedge A)$ are provable by **bC**:

| $\wedge$ | **1** | ½ | **0** |
|---|---|---|---|
| **1** | 1 | 1 | 0 |
| ½ | 1 | 1 | 0 |
| **0** | 0 | 0 | 0 |

| $\vee$ | **1** | ½ | **0** |
|---|---|---|---|
| **1** | 1 | 1 | 1 |
| ½ | 1 | 1 | 1 |
| **0** | 1 | 1 | 0 |

| $\to$ | **1** | ½ | **0** |
|---|---|---|---|
| **1** | 1 | 1 | 0 |
| ½ | 1 | 1 | 0 |
| **0** | 1 | 1 | 1 |

| | $\neg$ | $\circ$ |
|---|---|---|
| **1** | 0 | 1 |
| ½ | 1 | 0 |
| **0** | 1 | 1 |

where 1 and ½ are the distinguished values. By the way, the matrices of $\wedge$, $\vee$, $\to$, and $\neg$ are exactly the same matrices which originally defined the maximal three-valued logic $\textbf{P}^1$, proposed in [103], and mentioned in the subsection **2.5** as a logic which is paraconsistent and yet controllably explosive when in contact with any non-atomic formula.

As to the relations between contradictions and inconsistencies what we will find here are some variations on the intuitive idea that a contradiction should not be consistent (but not necessarily the other way around):

FACT 3.16 These are some special rules of **bC**, relating contradiction and consistency:
  (i) $A, \neg A \vdash_{\textbf{bC}} \neg \circ A$;
  (ii) $(A \wedge \neg A) \vdash_{\textbf{bC}} \neg \circ A$;
  (iii) $\circ A \vdash_{\textbf{bC}} \neg(A \wedge \neg A)$;
  (iv) $\circ A \vdash_{\textbf{bC}} \neg(\neg A \wedge A)$.
The converses of these rules do *not* hold in **bC**.

**Proof:** Use FACT 3.14 to prove (i), and left-adjunction to jump from this fact to (ii); play similarly to prove (iii) and (iv). To show that none of the converses of (ii)–(iv) are provable by **bC**, use the same matrices as in THEOREM 3.15(ii), substituting only the matrix for negation by this one to the right.

| | $\neg$ |
|---|---|
| **1** | 0 |
| ½ | 0 |
| **0** | 1 |

The significance of stating both (iii) and (iv) is to draw attention to the fact that, in what follows, logics will be shown in which, due to some unexpected asymmetry, only one of their converses hold. This is the case, for instance, for $C_1$, the first logic of the pioneering hierarchy of paraconsistent logics, $C_n$, $1 \leq n < \omega$, proposed by da Costa (cf. [49] or [50]). As we shall see, the converse of (iii) holds in $C_1$, while the converse of (iv) fails, so that $\neg(A \wedge \neg A)$ and $\neg(\neg A \wedge A)$ are *not* equivalent formulas in this logic (in this respect, see also THEOREM 3.21(iii)).

As the reader will learn in the next subsection (THEOREM 3.20), the regular forms of 'reasoning by contraposition' cannot be valid in any logic which is, as **bC** and its extensions (cf. THEOREM 3.13(iv)), both positively preserving with respect to classical logic and not partially explosive with respect to negation. But there are some restricted forms of it that hold already in **bC**:

FACT 3.17 These are some restricted forms of contraposition that hold in **bC**:
  (i)  $\circ B, (A \rightarrow B) \vdash_{\mathbf{bC}} (\neg B \rightarrow \neg A)$;
 (ii)  $\circ B, (A \rightarrow \neg B) \vdash_{\mathbf{bC}} (B \rightarrow \neg A)$;
(iii)  $\circ B, (\neg A \rightarrow B) \vdash_{\mathbf{bC}} (\neg B \rightarrow A)$;
 (iv)  $\circ B, (\neg A \rightarrow \neg B) \vdash_{\mathbf{bC}} (B \rightarrow A)$.

**Proof:** To check (i), let $\Gamma = \Delta = \Lambda = \{ \circ B, (A \rightarrow B), \neg B \}$ and apply FACT 3.14(ii) to $\Gamma \cup \{A\}$, so as to obtain $\Gamma \vdash_{\mathbf{bC}} \neg A$. From this it follows that $[\circ B, (A \rightarrow B) \vdash_{\mathbf{bC}} (\neg B \rightarrow \neg A)]$. Part (ii) is similar to (i). For parts (iii) and (iv) apply FACT 3.14(iii).

Now, may the reader be aware that rules such as $[\circ A, (A \rightarrow B) \vdash_{\mathbf{bC}} (\neg B \rightarrow \neg A)]$ *do not* hold in this logic!

**3.3 On what one cannot get.** If 'logic is about trade-offs', as Patrick Blackburn likes to put it, let us now start counting the dead bodies to see what we have irremediably lost, up to now. The connectives $\wedge$, $\vee$ and $\rightarrow$ of **bC**, for example, show up as quite independent from one another, and cannot be interdefined as in the classical case:

THEOREM 3.18 The following rule holds in **bC**:
    (i)  $(\neg A \rightarrow B) \vdash_{\mathbf{bC}} (A \vee B)$,
  but none of the following rules hold in **bC**:
  (ii)  $(A \vee B) \vdash_{\mathbf{bC}} (\neg A \rightarrow B)$;
 (iii)  $\neg(\neg A \rightarrow B) \vdash_{\mathbf{bC}} \neg(A \vee B)$;
  (iv)  $\neg(A \vee B) \vdash_{\mathbf{bC}} \neg(\neg A \rightarrow B)$;
   (v)  $(A \rightarrow B) \vdash_{\mathbf{bC}} \neg(A \wedge \neg B)$;
  (vi)  $\neg(A \wedge \neg B) \vdash_{\mathbf{bC}} (A \rightarrow B)$;
 (vii)  $\neg(A \rightarrow B) \vdash_{\mathbf{bC}} (A \wedge \neg B)$;
(viii)  $(A \wedge \neg B) \vdash_{\mathbf{bC}} \neg(A \rightarrow B)$;
  (ix)  $\neg(A \wedge B) \vdash_{\mathbf{bC}} (\neg A \vee \neg B)$;
   (x)  $(\neg A \vee \neg B) \vdash_{\mathbf{bC}} \neg(A \wedge B)$;
  (xi)  $\neg(\neg A \vee \neg B) \vdash_{\mathbf{bC}} (A \wedge B)$;
 (xii)  $(A \wedge B) \vdash_{\mathbf{bC}} \neg(\neg A \vee \neg B)$.

**Proof:** This is much easier to directly check after you take a look at the semantics and decision procedure of **bC**, in the paper [42]. But it also comes as a consequence from the fact that this is already valid for $C_{min}$, as we have proved in [39], and that **bC** is a conservative extension of it (THEOREM 3.12).

Notice that any uniform substitution of a component formula $C$ for its negation $\neg C$, or vice-versa, will not alter the fact that the above rules hold or not in **bC**. That is to say, for instance, that $(A \rightarrow \neg B) \vdash_{\mathbf{bC}} (\neg A \vee \neg B)$ does hold but $(\neg A \vee \neg B) \vdash_{\mathbf{bC}} (A \rightarrow B)$ does not. Of course, the failure of a rule such as $(A \vee \neg B) \vdash_{\mathbf{bC}} \neg(A \wedge \neg B)$ was already to be expected from the fact that $(A \vee \neg A)$ is provable (it is (Min10)) but $\neg(A \wedge \neg A)$ is not (see THEOREM 3.15(ii)).

Now, it should be crystal-clear that the above fact is only about **bC**, and that it does not necessarily carry on to stronger logics. In fact, it is not hard at all to check, for instance, that the three-valued maximal logic **LFI1**, whose matrices were presented in the subsection **2.4**, both extends **bC** and validates all the rules above, except for (ii) and (vi). Once more, the non-validity of (vi) is barely circumstantial, for there are logics extending **bC** in which it holds, such as the above mentioned $\mathbf{P}^i$ (see also the results 3.68 and 3.70, below). Still and all, there *is* a very good reason for the failure of (ii)! Indeed, this is a consequence of the following fact:

**THEOREM 3.19** The *disjunctive syllogism*, $[A, (\neg A \vee B) \vdash B]$, cannot hold in any paraconsistent extension of positive (classical or intuitionistic) logic.
**Proof:** Assume that it held. From (Min6), we would have that $[\neg A \vdash (\neg A \vee B)]$ and so, ultimately, we would conclude, by the transitivity of $\vdash$, that $[A, \neg A \vdash B]$.

Finally, as we have already advanced above, 'full' contraposition is lost (cf. [54]):

**THEOREM 3.20** The regular forms of *contraposition*, such as $[(A \rightarrow B) \vdash (\neg B \rightarrow \neg A)]$, cannot hold irrestrictedly in any paraconsistent extension of **bC**. Furthermore, they cannot hold in any extension of the positive classical logic which happens to be not $\neg$-partially explosive.
**Proof:** If the above rule held in a logic **L** that extends the positive classical logic, from (Min1) we would obtain $[B \vdash (A \rightarrow B)]$, and from (MP) we obtain $[(A \rightarrow B), \neg B \vdash \neg A]$. These two rules would ultimately lead to $[B, \neg B \vdash \neg A]$, and so **L** would be partially explosive with respect to negated propositions. If we assume **L** to be **bC**, then a particular case of $[B, \neg B \vdash \neg A]$ would be $[B, \neg B \vdash \neg \neg C]$, taking $A$ as $\neg C$, and (Min11) would then give $[B, \neg B \vdash C]$, and so it would not be paraconsistent at all. Indeed, this addition of contraposition to **bC** would simply cause the collapse of the resulting logic into classical logic (by THEOREM 3.9). Still some other forms of this contraposition rule, such as $[(\neg A \rightarrow B) \vdash (\neg B \rightarrow A)]$, could be ruled out even without recurring to (Min11), or to partial explosion.

The use of the disjunctive syllogism (THEOREM 3.19) constitutes indeed the kernel of the well-known argument laid down by C. I. Lewis for the derivation of (PPS) in classical logic (cf. [73], pp.250ff), and this was, in fact, a rediscovery of an argument used by the Pseudo-Scotus, much before.[14] The use of contraposition (THEOREM 3.20) to the same purpose was pointed out in an argument by Popper (cf. [89], pp.320ff). Of course, in a logic where both the disjunctive syllogism and contraposition are invalid derivations, these arguments do not apply as such.

---

[14] See Duns Scotus's *Opera Omnia*, pp.288ff. Cf. also note 3.

The failure of contraposition gives us a good reason for having doubts also about the validity of the *intersubstitutivity of provable equivalents*, which states that, given a schema $\sigma(A_1, \ldots, A_n)$:

$$\forall B_1 \ldots \forall B_n \left[ (A_1 \dashv\vdash B_1) \text{ and } \ldots \text{ and } (A_n \dashv\vdash B_n) \right] \Rightarrow \qquad \text{(IpE)}$$
$$\left[ \sigma(A_1, \ldots, A_n) \dashv\vdash \sigma(B_1, \ldots, B_n) \right].$$

Now, as a particular example, if we had (IpE), from $A \dashv\vdash B$ we would immediately derive, for instance, $\neg A \dashv\vdash \neg B$. But this is not the case here. Indeed, in what follows we exhibit some samples of that failure in **bC**:

**THEOREM 3.21** In **bC**:

   (i)  $(A \wedge B) \dashv\vdash_{\mathbf{bC}} (B \wedge A)$ holds, but $\neg(A \wedge B) \dashv\vdash_{\mathbf{bC}} \neg(B \wedge A)$ does not;

   (ii)  $(A \vee B) \dashv\vdash_{\mathbf{bC}} (B \vee A)$ holds, but $\neg(A \vee B) \dashv\vdash_{\mathbf{bC}} \neg(B \vee A)$ does not;

   (iii)  $(A \wedge \neg A) \dashv\vdash_{\mathbf{bC}} (\neg A \wedge A)$ holds, but $\neg(A \wedge \neg A) \dashv\vdash_{\mathbf{bC}} \neg(\neg A \wedge A)$ does not.

**Proof:** The parts which hold are easy, using positive classical logic. Now, to check that none of the other parts hold, even if axioms and rules of **bC** are taken into consideration, use the same matrices and distinguished values as in THEOREM 3.15(ii), changing only the values of $(1 \wedge \frac{1}{2})$ and $(1 \vee \frac{1}{2})$ from 1 to $\frac{1}{2}$ (but leaving the values of $(\frac{1}{2} \wedge 1)$ and $(\frac{1}{2} \vee 1)$ as they are, equal to 1).

**COROLLARY 3.22** (IpE) does not hold for **bC**.

The reader should keep in mind that this last result is, initially, only about **bC**, and that some deductive extensions of it may fix some or even all the counter-examples to intersubstitutivity. Now, given that (IpE) holds for classical logic, it will obviously hold for the positive (classical) fragment of **bC** as well, that is, for the set of formulas in which neither $\neg$ nor $\circ$ occur. Adding contraposition as a new inference rule, it is easy to see, by the transitivity of the consequence operator and the Deduction Metatheorem, that one could extend (IpE) from positive logic to include also the fragment of **bC** containing negation. But then (bold) paraconsistency would be lost, as we learn from THEOREM 3.20! What happens, though, is that the contraposition inference rule is much more than one needs in order to obtain intersubstitutivity for the consistencyless fragments of our logics. In fact, any of the following 'contraposition' deduction rules would of course do the job equally well (cf. [107] and [105]):

$$\forall A \, \forall B \, [(A \vdash B) \Rightarrow (\neg B \vdash \neg A)]; \qquad \text{(RC)}$$
$$\forall A \, \forall B \, [(A \dashv\vdash B) \Rightarrow (\neg B \vdash \neg A)]. \qquad \text{(EC)}$$

It is obvious that (EC) can be inferred from (RC), and Urbas has shown in [107] that the paraconsistent logic obtained by adding (EC) to $C_\omega$ is extended by the paraconsistent logic obtained by the addition of (RC) to $C_\omega$ (and both, of course, are extended by classical logic). So, it *is* possible to obtain paraconsistent extensions of $C_\omega$ (and also of $C_{min}$, for Urbas's proof of non-collapse into classical logic by the addition of (EC) also applies to this logic), but then these new logics can all still be shown to lack a bottom particle (as in THEOREM 3.4), constituting thus no **LFI**s! The question then would be if (IpE) could be obtained for *real* **LFI**s. The closest we will get to this here is showing, in THEOREM 3.53, that there are fragments of classical logic extending **bC** for which (IpE) holds, but then these specific fragments turn out not to be paracon-sistent in our sense. At any rate, for various other classes of **LFI**s we will show that such intersubstitutivity results are just unattainable, as shown in THEOREM 3.51 (see also, for instance, FACT 3.74).

To be sure, one does not need to blame *paraconsistency* for these last few negative results. As the reader will see below, the eccentricities in THEOREM 3.21 can be fixed by some extensions of **bC**. As for THEOREM 3.20, one could always throw away some piece of the positive classical logic in an extreme effort to avoid its consequences. This is what is done, for instance, by some logics of relevance. This could, however, have the effect of throwing the baby out with the bath water —most such logics, if not all, will also dismiss the useful Deduction Metatheorem or, regrettably enough, *modus ponens*. Now, suppose that, driven by itches of relevance, one was taken to consider logics such that $(A, B \nvdash A)$. This would definitely mean, thus, that their consequence relations would be no more than 'cautiously reflexive'. If one still insisted that $(A \Vdash A)$ should hold, then the logics produced would be non-monotonic as well. This would mean, of course, that many of the results that we attained in the last section would not be immediately adaptable to such logics (and this remark also applies to adaptive logics, once they are also non-monotonic, even if for other reasons). These are not problems of actual relevance logics, nevertheless, as they are usually relevant only at the level of theoremhood (always invalidating $(A \to (B \to A))$, while in some cases still validating $(A \to A)$), but still not at the level of their consequence relations, as conjectured above (see, for instance, [3] or [9]) —and of course, in all such cases, the Deduction Metatheorem cannot hold. But yes, we had better push our exposition on, instead of scrubbing this matter here any further.

### 3.4 Letting bC talk about (dual) inconsistency.
The reader may find it a bit awkward, indeed, that we would be calling **bC** a logic of formal *inconsistency*, since it only has a connective expressing *consistency*, but not its opposed concept. So, for *us* to be more consistent, let's now consider a further extension of our language, this time adding a new unary connective, •, to represent inconsistency. The intended interpretation about the dual relation between consistency and inconsistency would require exactly that each of these concepts should be opposed to the other. But how do we formalize this? Consider the following additional axiomatic rule:

(bc2)     $\neg \bullet A \vdash_{\mathbf{bbC}} \circ A$.          'If $A$ is not inconsistent, then it is consistent'

This is surely a must, but in fact it does not represent much of an addition. Indeed, consider its contrapositional variation:

(bc3)     $\neg \circ A \vdash_{\mathbf{bbC}} \bullet A$.          'If $A$ is not consistent, then it is inconsistent'

The lack of contraposition (see THEOREM 3.20), despite the presence of some restricted forms of it (such as in FACT 3.17) can be partly blamed for the fact that **bC** plus (bc2) can still not prove (bc3). Indeed:

THEOREM 3.23  (bc3) is not provable by **bC** plus (bc2).
**Proof:** Just consider three-valued matrices such that: $v(A \land B) = 0$ if $v(A) = 0$ or $v(B) = 0$, and $v(A \land B) = 1$, otherwise; $v(A \lor B) = 0$ if $v(A) = 0$ and $v(B) = 0$, and $v(A \lor B) = 1$, otherwise; $v(A \to B) = 0$ if $v(A) \neq 0$ and $v(B) = 0$, and $v(A \to B) = 1$, otherwise; $v(\neg A) = 1 - v(A)$; and the matrices for the non-classical connectives are the ones demonstrated on the right. 0 is the only non-distinguished value.

| | $\circ$ | $\bullet$ |
|---|---|---|
| **1** | 0 | 1 |
| ½ | 0 | 1 |
| **0** | ½ | 0 |

So, let us now, for the sake of symmetry, define the logic **bbC** as given by the addition of both (bc2) and (bc3) to the basic logic of (in)consistency, **bC**. This is still not much… for consider now the converses of these rules:

(bc4)     $\bullet A \vdash_{\mathbf{bbbC}} \neg \circ A$;                      'If $A$ is inconsistent, then it is not consistent'

(bc5)     $\circ A \vdash_{\mathbf{bbbC}} \neg \bullet A$.                      'If $A$ is consistent, then it is not inconsistent'

Will these hold in **bbC**? The answer is once more in the negative:

**THEOREM 3.24** Neither (bc4) nor (bc5) are provable by **bbC**.
**Proof:** Consider the same three-valued matrices for the binary connectives as in THEOREM 3.23, but let now negation be such that $v(\neg A) = 0$ if $v(A) \neq 0$, and $v(\neg A) = 1$, otherwise. The non-classical connectives will now be defined by the new matrices to the right. Once more, 0 is the only non-distinguished value.

|       | ∘   | •   |
|-------|-----|-----|
| **1** | 0   | 1   |
| ½     | ½   | 1   |
| **0** | 0   | 1   |

In reality, the situation is even worse than it may appear at first sight, though predictable. It happens that, once more, it is not enough to add just one of (bc4) or (bc5) to **bbC** —the other one would still not be provable. Indeed:

**THEOREM 3.25** (i) (bc4) is not provable by **bbC** plus (bc5); (ii) (bc5) is not provable by **bbC** plus (bc4).
**Proof:** Consider now the four-valued matrices where $\wedge$, $\vee$, $\rightarrow$ and $\neg$ are once more defined as in THEOREM 3.23 (only that now they have a wider domain, with four values). For part (i), let ∘ and • be given by the matrices to the right. For part (ii), just modify ∘ so that $\circ(\tfrac{2}{3}) = \tfrac{1}{3}$ (and no more $\tfrac{2}{3}$); modify also • in the contrary sense, so that $\bullet(\tfrac{2}{3}) = \tfrac{2}{3}$ (and no more $\tfrac{1}{3}$). In both cases, only 0 should be taken to be a non-distinguished value.

|       | ∘   | •   |
|-------|-----|-----|
| **1** | 1   | 0   |
| ⅔     | ⅔   | ⅓   |
| ⅓     | 0   | ⅔   |
| **0** | 1   | 0   |

Taking the above results into account, we will now define the logic **bbbC** to be given by the addition of both (bc4) and (bc5) to the preceding **bbC**.

It is important to note that the last theorem above also shows that it is ineffective trying to introduce the inconsistency connective in the logic **bC** simply by setting, by definition, $\bullet A \overset{\text{def}}{=} \neg \circ A$. The reason is that, even though this would automatically guarantee that $\bullet A \dashv\vdash \neg \circ A$, and that $\neg \bullet A \dashv\vdash \neg \neg \circ A$, and so on, just by definition and reflexivity, this would *not* guarantee as well that, for instance, we would have $\circ A \vdash \neg \bullet A$. Indeed, to check this you may here just reconsider THEOREM 3.25(ii). So, the relation between ∘ and • cannot, in the cases of **bC** and **bbC**, be characterized by a simple definition. Despite this, one may now establish new presentations for some previous facts and theorems, just slightly different from before:

**THEOREM 3.26** The results 3.11, 3.14, 3.15, 3.16, and 3.17 are all valid for **bbbC**, and are still valid if one substitutes any occurrence of ∘ for ¬•, and ¬∘ for •.
**Proof:** This is routine, just using (bc2)–(bc5). For 3.15 and 3.16 remember to add a matrix for •, just negating the matrix for ∘ presented in THEOREM 3.15(ii).

So, could the relation between ∘ and • be characterized by a definition, now that we have **bbbC**? Another NO is the answer. For if a definition such as $\bullet A \overset{\text{def}}{=} \neg \circ A$ were feasible, this would mean, given (bc5): $\circ A \vdash_{\mathbf{bbbC}} \neg \bullet A$, that $\circ A \vdash_{\mathbf{bbbC}} \neg \neg \circ A$ should hold just by straightforward substitution. But, as it happens, this last rule does *not* hold in **bbbC**:

**THEOREM 3.27**  Neither $\circ A \rightarrow \neg\neg\circ A$ nor $\bullet A \rightarrow \neg\neg\bullet A$ are provable by **bbbC**.

**Proof:**  Consider once more the same three-valued matrices for the binary connectives given in THEOREM 3.23, 0 as the only non-distinguished value, but now let the unary connectives be those pictured to the right.

| | ¬ | ∘ | • |
|---|---|---|---|
| **1** | 0 | ½ | ½ |
| **½** | 1 | 0 | 1 |
| **0** | 1 | 0 | 1 |

Evidently, the above matrices must also display the non-provability by **bbbC** of the schema $A \rightarrow \neg\neg A$, the converse of (Min11). But if the validity of $A \rightarrow \neg\neg A$ clearly implies the validity of the two schemas in THEOREM 3.27, the validity of those schemas certainly *does not* imply the validity of $A \rightarrow \neg\neg A$. Indeed:

**THEOREM 3.28**  $A \rightarrow \neg\neg A$ is not provable by **bbbC** plus $\circ A \rightarrow \neg\neg\circ A$ and $\bullet A \rightarrow \neg\neg\bullet A$.

**Proof:**  Consider the same matrices and distinguished values as in THEOREM 3.27, only that now $\circ$ is constant and equal to 0, and $\bullet$ is constant and equal to 1.

Now, if we added to **bbbC** the axioms $\circ A \vdash \neg\neg\circ A$ and $\bullet A \vdash \neg\neg\bullet A$ this would only shift our problem to proving that $\neg\bullet A \vdash \neg\neg\neg\bullet A$ and $\neg\circ A \vdash \neg\neg\neg\circ A$ hold, and so on, and so forth. Of course, these would be all guaranteed if we now defined **bbbbC** by the addition to **bbbC** of an infinite number of axiomatic rules, to the effect that $\neg^n\circ A \vdash_{\textbf{bbbbC}} \neg^{n+2}\circ A$ and $\neg^n\bullet A \vdash_{\textbf{bbbbC}} \neg^{n+2}\bullet A$, where $\neg^m$ denotes $m$ occurrences of negation in a row. We could also solve all of this at once by fixing $A \vdash \neg\neg A$ as a new axiomatic rule, but we argue that it is a bit too early for this last solution —indeed, there is a gamut of interesting **C**-systems in which this axiom does *not* hold, and we would rather explore them first. So, let us study first, in what follows, some other forms of obtaining the intended duality between $\circ$ and $\bullet$ using a finite set of schemas, and without yet incorporating $A \vdash \neg\neg A$ as a rule.

**3.5  The logic Ci, where contradiction and inconsistency meet.**  While strengthening **bC**, we have been trying to keep up with the intended duality between consistency and inconsistency. But, given the new version of the FACT 3.16 obtained in THEOREM 3.26 (which also applies to **bbbbC**), we know that in any of the logics **b(b(b(b)))C** a contradiction implies an inconsistency, but not the other way around —so, *this* situation has still not been changed. Now, the distinction between contradiction and inconsistency is a contribution of the present study, and we are unaware of any other formal attempts to do so in the same way as we do here. What will happen then if we now introduce new axioms in order to finally obtain the identification of contradiction and inconsistency, getting closer this way to the other paraconsistent logics in the literature? Let's do it. Consider the two following axiomatic rules:

| | | |
|---|---|---|
| (ci1) | $\bullet A \vdash_{\textbf{Ci}} A$; | 'If $A$ is inconsistent, then $A$ should be the case' |
| (ci2) | $\bullet A \vdash_{\textbf{Ci}} \neg A$. | 'If $A$ is inconsistent, then $\neg A$ should be the case' |

Given the classical properties of conjunction, these two rules will evidently have the same effect as the following single one:

| | | |
|---|---|---|
| (ci) | $\bullet A \vdash_{\textbf{Ci}} (A \wedge \neg A)$. | 'An inconsistency implies a contradiction' |

So, let's call **Ci** the logic obtained by the addition of (ci1) and (ci2) (or, equivalently, the addition of (ci)) to **bbbC**, that is, the logic axiomatized by (Min1)–(Min11), (bc1)–(bc5), (ci), and (MP). In **Ci** we finally have that $\bullet A$ and $(A \wedge \neg A)$ are equivalent for-

(bc5), (ci), and (MP). In **Ci** we finally have that $\bullet A$ and $(A \wedge \neg A)$ are equivalent formulas, and we shall see that this will make a BIG difference on **Ci**'s deductive strength.

First, let us note that, even though we now have, in **Ci**, the converse of parts (i) and (ii) of FACT 3.16, the converses of parts (iii) and (iv) still do *not* hold. Indeed:

**FACT 3.29** This rule does hold in **Ci**:

(i) $\neg \circ A \vdash_{\textbf{Ci}} (A \wedge \neg A)$,

but the following rules do not:

(ii) $\neg(A \wedge \neg A) \vdash_{\textbf{Ci}} \circ A$;

(iii) $\neg(\neg A \wedge A) \vdash_{\textbf{Ci}} \circ A$.

**Proof:** The first part is obvious. For the following ones, consider, for instance, the three-valued matrices such that:

$v(A \wedge B) = min(v(A), v(B))$;

$v(A \vee B) = max(v(A), v(B))$;

$v(A \to B) = v(B)$, if $v(A) \neq 0$, and $v(A \to B) = 1$, otherwise;

$v(\neg A) = 1 - v(A)$;

$v(\circ A) = 0$, if $v(A) = v(\neg A)$, and $v(\circ A) = 1$, otherwise;

$v(\bullet A) = 1$, if $v(A) = v(\neg A)$, and $v(\bullet A) = 0$, otherwise,

where 0 is the only non-distinguished value. The attentive reader might have noticed that these are exactly the matrices defining the already mentioned **LFI1**, in the subsection **2.4**.

So, this last theorem reminds us that, even though in **Ci** we *do* have an equivalent way of referring to inconsistency using just the classical language, this does not mean that we should also have an immediate **CPL**-linguistic equivalent manner of referring to *consistency* as well (but confront this with what happens in the case of the **dC**-systems, in the subsection **3.8**)! There are, however, many other things that we *do* have. For instance, in **Ci** the THEOREM 3.15 is still entirely valid. Indeed:

**THEOREM 3.30** $\neg(A \wedge \neg A)$ and $\neg(\neg A \wedge A)$ are not top particles in **Ci** (also the formula $(A \to \neg \neg A)$ is still not provable).

**Proof:** Use again the matrices of $\mathbf{P}^1$ (in THEOREM 3.15(ii)), adding a matrix for '$\bullet$' by negating the matrix for '$\circ$'.

We also have in **Ci** some new ways of formulating gentle explosion and the *reductio* deduction rules:

**FACT 3.31** The following rules hold in **Ci**:

(i) $\circ A, \bullet A \vdash_{\textbf{Ci}} B$;

(ii) $\circ A, \neg \circ A \vdash_{\textbf{Ci}} B$;

(iii) $\bullet A, \neg \bullet A \vdash_{\textbf{Ci}} B$;

(iv) $[(\Gamma, B \vdash_{\textbf{Ci}} \circ A)$ and $(\Delta, B \vdash_{\textbf{Ci}} \bullet A)] \Rightarrow (\Gamma, \Delta \vdash_{\textbf{Ci}} \neg B)$;

(v) $[(\Gamma, B \vdash_{\textbf{Ci}} \circ A)$ and $(\Delta, B \vdash_{\textbf{Ci}} \neg \circ A)] \Rightarrow (\Gamma, \Delta \vdash_{\textbf{Ci}} \neg B)$;

(vi) $[(\Gamma, B \vdash_{\textbf{Ci}} \bullet A)$ and $(\Delta, B \vdash_{\textbf{Ci}} \neg \bullet A)] \Rightarrow (\Gamma, \Delta \vdash_{\textbf{Ci}} \neg B)$.

**Proof:** Part (i) comes from (ci) and (bc1), parts (ii) and (iii) come from part (i) if you use (bc2)–(bc5). Rules (iv), (v) and (vi) are variations on FACT 3.14(ii), using the previous rules.

Parts (ii) and (iii) of FACT 3.31 simply show **Ci** to be controllably explosive in contact either with a consistent or with an inconsistent formula. In fact, in **Ci** one can go on to prove a much more intimate connection between consistency and control-lable explosion, and this will reveal some even stronger consequences of the new axiomatic rule, (ci), that we now consider:

FACT 3.32 A particular given schema in **Ci** (or in any extension of this logic) is consistent if, and only if, **Ci** is controllably explosive in contact with this schema.

**Proof:** To show that $[(\Gamma \vdash_{\mathbf{Ci}} \circ A) \Rightarrow (\Gamma, A, \neg A \vdash_{\mathbf{Ci}} B)]$ just invoke axiom (bc1) and the transitivity of $\vdash$. For the converse, note that, from (ci) and (bc3), one may obtain $[\neg \circ A \vdash_{\mathbf{Ci}} (A \wedge \neg A)]$, and so, from the supposition that $(\Gamma, A, \neg A \vdash_{\mathbf{Ci}} B)$ it follows that $\neg \circ A$ is a bottom particle. One may then conclude, as in THEOREM 3.10, that $\vdash_{\mathbf{Ci}} \neg \neg \circ A$, and, by (Min11), that $\vdash_{\mathbf{Ci}} \circ A$.

FACT 3.33 These are some special theses of **Ci**:

(i) $\vdash_{\mathbf{Ci}} \circ \circ A$;

(ii) $\vdash_{\mathbf{Ci}} \neg \bullet A$;

(iii) $\vdash_{\mathbf{Ci}} \circ \bullet A$;

(iv) $\vdash_{\mathbf{Ci}} \neg \bullet \bullet A$.

**Proof:** Parts (i) and (iii) come directly from FACT 3.32 and from parts (ii) and (iii) of FACT 3.31. For (ii) and (iv), use (bc2) and the previous parts.

This last result (check also [54]) implies that **Ci** will not have consistency or inconsistency appearing at different levels: both consistent and inconsistent formulas are consistent (in contrast to what happened in the case of **bC** —see THEOREM 3.10—, where no formula was provably consistent), and none of them is inconsistent (check also FACT 3.50 for a much stronger version of the last fact in **Ci**).

The reader will recall from the subsection **3.3** that contraposition inference rules not only did not hold in **bC** but could not even be added to any paraconsistent extension of it (THEOREM 3.20). **Ci** can be shown to count, nevertheless, with more restricted forms of contraposition than **bC** (compare the following with FACT 3.17):

FACT 3.34 These are some restricted forms of contraposition introduced by **Ci**:

(i) $(A \rightarrow \circ B) \vdash_{\mathbf{Ci}} (\neg \circ B \rightarrow \neg A)$;

(ii) $(A \rightarrow \neg \circ B) \vdash_{\mathbf{Ci}} (\circ B \rightarrow \neg A)$;

(iii) $(\neg A \rightarrow \circ B) \vdash_{\mathbf{Ci}} (\neg \circ B \rightarrow A)$;

(iv) $(\neg A \rightarrow \neg \circ B) \vdash_{\mathbf{Ci}} (\circ B \rightarrow A)$.

**Proof:** To check (i), let $\Gamma = \Delta = \{(A \rightarrow \circ B), \neg \circ B\}$ and apply FACT 3.31(v) to $\Gamma \cup \{A\}$, so as to obtain $\Gamma \vdash_{\mathbf{Ci}} \neg A$. This will give the desired result. Alternatively, one could use directly FACT 3.17(i) and note that $\circ \circ B$ is a theorem of **Ci** (this is FACT 3.33(i)). The other parts are similar, and we leave them as easy exercises to the reader.

Note that all rules in the last result continue to be valid if one substitutes any '$\circ$' for '$\neg \bullet$', and any '$\neg \circ$' for '$\bullet$'. On the other hand, rules such as $[(\circ A \rightarrow B) \vdash_{\mathbf{Ci}} (\neg B \rightarrow \neg \circ A)]$ *do not* hold in this logic!

Now, we have learned from COROLLARY 3.22 that the intersubstitutivity of provable equivalents, (IpE), does not hold for **bC**. The same result is true for **Ci**, and still the same counter-examples mentioned before can be presented here:

THEOREM 3.35 (IpE) does not hold for **Ci**.

**Proof:** Add to the matrices on THEOREM 3.21 one matrix for $\bullet$ such that $v(\bullet A) = 1 - v(\circ A)$, and check that all the new axioms, defining **Ci** from **bC**, still hold.

Now, in order to go one step further from the actual absence of contraposition in **Ci**, let us recall that in the subsection **3.3** it has been pointed out that the addition of some of the deduction 'contraposition' rules (EC) or (RC) would have been equally sufficient for obtaining (IpE) for consistencyless fragments of our paraconsistent logics. It seems, nevertheless, that obtaining (IpE) will not be an easy task, after all:

**FACT 3.36** The addition of (RC): $[(A \vdash B) \Rightarrow (\neg B \vdash \neg A)]$ to **Ci** causes its collapse into classical logic.

**Proof:** From (ci1) and (ci2), plus (bc3), one obtains, respectively, that $\neg \circ A \vdash_{\mathbf{Ci}} A$, and $\neg \circ A \vdash_{\mathbf{Ci}} \neg A$. Applying (RC) and (Min11) one would have then $\neg A \vdash_{\mathbf{Ci}} \circ A$ and $\neg \neg A \vdash_{\mathbf{Ci}} \circ A$. But then, using the proof by cases, one would conclude that $\vdash_{\mathbf{Ci}} \circ A$, that is, all formulas would be consistent. Looking at THEOREM 3.9 and (bc1), one sees that this was exactly what was lacking in order for classical logic to be characterized.

So, (RC) must be ruled out as an alternative in order to obtain (IpE), in the case of **Ci**. As for (EC), its possible addition to **Ci** will be discussed below, in THEOREM 3.51, FACT 3.52, and the subsequent commentaries on these results.

The new restricted forms of contraposition in FACT 3.34 are, in any case, strong enough for us to show that **Ci** has some redundant axioms as it is. Indeed:

**FACT 3.37** In **Ci**: (i) (bc2) proves (bc3), and vice-versa. (ii) (bc4) proves (bc5), and vice-versa.

Other interesting consequences of (ci) are those that we shall call 'Guillaume's Theses', which regulate the propagation of consistency and the back-propagation of inconsistency through negation:

**FACT 3.38** **Ci** also proves the following:
   (i) $\circ A \vdash_{\mathbf{Ci}} \circ \neg A$;
   (ii) $\bullet \neg A \vdash_{\mathbf{Ci}} \bullet A$.

**Proof:** From (ci) and (bc4), we have that $[\neg \circ \neg A \vdash_{\mathbf{Ci}} (\neg A \wedge \neg \neg A)]$, from $C_{min}$ we have that $[(\neg A \wedge \neg \neg A) \vdash_{\mathbf{Ci}} (A \wedge \neg A)]$, and from FACT 3.16(ii) we know that $[(A \wedge \neg A) \vdash_{\mathbf{Ci}} \neg \circ A]$. So, ultimately, we have the rule $[\neg \circ \neg A \vdash_{\mathbf{Ci}} \neg \circ A]$. By (bc3) and (bc4) we prove part (ii) of our fact. Part (i) comes from this same rule, by an application of FACT 3.34(iv).

This last result will provide us with some other forms for the theses in FACT 3.33, such as:

**FACT 3.39** These are also some special theses of **Ci**:
   (i) $\vdash_{\mathbf{Ci}} \circ \neg \circ A$;
   (ii) $\vdash_{\mathbf{Ci}} \neg \bullet \neg \circ A$;
   (iii) $\vdash_{\mathbf{Ci}} \circ \neg \bullet A$;
   (iv) $\vdash_{\mathbf{Ci}} \neg \bullet \neg \bullet A$.

It will also be useful to note that here we have (contrasting with THEOREM 3.30, which informed us, among other things, that $[A \nvdash_{\mathbf{Ci}} \neg \neg A]$):

**FACT 3.40** Here are some more special theses of **Ci**:
   (i) $\circ A \vdash_{\mathbf{Ci}} \neg \neg \circ A$;
   (ii) $\bullet A \vdash_{\mathbf{Ci}} \neg \neg \bullet A$.

**Proof:** These will follow directly if you apply FACT 3.34 twice. The reader might remember that we lacked these forms in **bbbC** (this was THEOREM 3.27).

Now, do we obtain in **Ci** that intended duality between consistency and inconsistency? The answer is YES. This is the topic for our next subsection.

### 3.6 On a simpler presentation for Ci.

The logic **Ci** provides us with a sufficient environment to prove a kind of restricted *intersubstitutivity* or *replacement* theorem. While we know from THEOREM 3.35 that full replacement for the formulas of **Ci** does not obtain, our present restricted forms of contraposition, nevertheless, will help us to show that intersubstitutivity *does* hold if only we are talking only about substituting some formula whose outmost operator is '∘' by this same formula, but now having '¬•' in the place of that '∘', or if we will substitute some formula whose outmost operator is '•' by this same formula, but now having '¬∘' in the place of that '•'. In simpler terms, what we are saying is that we can now take just one of the operators '∘' and '•' as primitive, and define the other in terms of the negation of that first one. So, we will now show that:

THEOREM 3.41 An equivalent axiomatization for **Ci** is obtained if we consider only axioms (Min1)–(Min11), (bc1), (ci), and (MP), and set one of these two definitions:

   (i) $\bullet A \overset{\text{def}}{=} \neg \circ A$;

   (ii) $\circ A \overset{\text{def}}{=} \neg \bullet A$.

**Proof:** Consider part (i) to be the case. This means that we can take (bc3): $\neg \circ A \vdash_{\mathbf{Ci}} \bullet A$ and (bc4): $\bullet A \vdash_{\mathbf{Ci}} \neg \circ A$, for granted, simply by definition. Now, (bc2): $\neg \bullet A \vdash_{\mathbf{Ci}} \circ A$, will be the case if, and only if, given the definition of '•', $\neg \neg \circ A \vdash_{\mathbf{Ci}} \circ A$ is the case —and it is, because of (Min11). As to (bc5): $\circ A \vdash \neg \bullet A$, it will be the case if, and only if, $\circ A \vdash_{\mathbf{Ci}} \neg \neg \circ A$ is the case —and it is, this time thanks to FACT 3.40(i). An alternative, and much simpler way, of checking that (bc2) and (bc5) should hold here is by taking FACT 3.37 into consideration. The axiomatic rule (bc1): $\circ A$, $A$, $\neg A \vdash_{\mathbf{Ci}} B$ is already in the 'standard form' (we are here eliminating all occurrences of '•'s and leaving only '∘'s), and the rule (ci): $\bullet A \vdash_{\mathbf{Ci}} (A \wedge \neg A)$ can be exchanged, by the definition of '•', that is, by (bc3) and (bc4), for $\neg \circ A \vdash_{\mathbf{Ci}} (A \wedge \neg A)$. Now we have shown that all occurrences of '•' in the axioms of **Ci** can be substituted by an occurrence of '¬∘', and all occurrences of '¬•' in the axioms of **Ci** can be substituted by an occurrence of '∘'. So, if you would have proven a formula in which, respectively, an inconsistency connective '•' or its negated form '¬•' appears at some point, you can now rewrite the proof using the new versions of the axioms above and what will appear in the end will be, respectively, a negated consistency connective '¬∘', or simply the connective '∘'. For part (ii) the procedure is entirely analogous, but now use FACT 3.40(ii), or FACT 3.37 again, when necessary.

So, this last result provides us with a restricted form of replacement theorem for consistent formulas, and guarantees the intended duality between ∘ and •, which could not be obtained in the subsection **3.4**, within **bbbbC** or its fragments. With such a result in hand we need make no big effort to verify that formulas such as $\neg(\circ A \wedge \neg \circ A)$, $\neg(\circ A \wedge \bullet A)$, $\neg(\neg \bullet A \wedge \neg \circ A)$ and $\neg(\neg \bullet A \wedge \bullet A)$ are all equivalent, which could, otherwise, be quite a non-trivial task!

The reason why we can obtain this new axiomatization, as the reader will make out after he is introduced to the semantics of **Ci**, in [42], is that, truth-functionally based on the non-classical behavior of negation, both the consistency and the inconsistency operators of this logic will work quite 'classically'.

**3.7 Using LFIs to talk about classical logic.** At this point, working with **Ci**, perhaps the question would arise as to how far we are from classical propositional logic, **CPL**. The answer is: a lot —and just a little bit. As we are not presupposing any kind of doublethinking, let us then reformulate a few things for the question, and its answer, really to make sense.

To start with, it is hard to compare two logics if they 'talk about different things', and are so disjoint that none of them is an extension of the other. For $C_{min}$ was a conservative extension of positive classical logic, but it was a fragment, in the same language, of 'full' **CPL**, as we know, and **bC** was a conservative extension of $C_{min}$. Thus, **Ci**, which is a deductive extension of **bC**, happens to be written in a richer language than that of **CPL**, but it does not contain all classical inferences, and so these two logics are hardly comparable. Now, this is easy to fix. Let us also conservatively extend **CPL** by the addition of connectives for consistency and inconsistency, whose matrices will be such that ∘ takes always the distinguished value 1, and •, on the contrary, is constant and equal to 0. We will designate this 'new' logic, obtained by such an extension of **CPL**, *extended classical logic*, or **eCPL**. Of course, **eCPL** can be easily axiomatized by the addition of an axiomatic schema such as:

(ext)    $\vdash_{\textbf{eCPL}} \circ A$    'Every $A$ is consistent'

to any axiomatization of **CPL**, like the one mentioned in THEOREM 3.9. The inconsistency connective, •, can be here introduced as a definition: $\bullet A \overset{\text{def}}{=} \neg \circ A$, just as in THEOREM 3.41(i). In this way we obtain an extension of classical logic which looks as a logic of formal inconsistency (see (D20)), having an operator expressing consistency, and of course an axiomatic rule such as $[\circ A, A, \neg A \vdash_{\textbf{eCPL}} B]$, expressing finite gentle explosion, will hold in **eCPL**. But, as it happens, given axiom (ext), we know that **eCPL** is not only finitely gently explosive (and so, non-trivial), but explosive as well. It is, in fact, a *consistent* logic (see (D19)), instead of an **LFI**.

Well and good, but is **Ci** now to be characterized as a deductive fragment of **eCPL**? Indeed! Just check that all axioms of **Ci** are validated by the matrices of **eCPL**, and that's it. So, **Ci** is in fact a fragment of an alternative formulation of classical logic, and this of course will guarantee that **Ci** is a non-contradictory logic (once **Ci** is not explosive, but it is a fragment of **eCPL**, and **eCPL** is still at least as explosive and non-trivial as **CPL** was, and consequently it cannot prove a contradiction). Is that all to it? No, because we will now see that we can still use **Ci** to reproduce in a very faithful way every inference of **CPL** (or of **eCPL**)!

How can this be done? Remember that **Ci** has a strong (or supplementing) negation, which can be defined, as in the case of **bC**, by setting $\sim A \overset{\text{def}}{=} (\neg A \wedge \circ A)$. But, in **bC**, even though this negation had the power of producing (supplementing) explosions, it could not still be said to have all properties of a *classical negation*. Indeed:

THEOREM 3.42 The strong negation ~, in **b(b(b(b)))C**, is not classical.

**Proof:** Just consider once more the classical matrices for the classical connectives, as in the above definition of **eCPL**, but now exchange the matrices of ∘ and •, letting ∘ be constant and equal to 0 (and not to 1, as before), and letting • be constant and equal to 1 (and not to 0, as before). It is easy to see that all axioms and rules of

**b(b(b(b)))C** are validated by such matrices, but (ci) and (ext) are not, and consequently formulas such as $(A \vee \sim A)$ and $(A \to \sim \sim A)$ (recall the definition of $\sim A$) are not validated as well, being independent from all logics we have exposed previous to **Ci**.

This is an interesting result that shows that being explosive is not enough to make a negation classical.[15] But what would be enough? Well, given the axiomatization of classical logic in THEOREM 3.9, we know that any connective $\div$ added in an axiomatic environment where (Min1)–(Min9) hold and which is such that:

(Alt10)    $\vdash_{\textbf{Alt}} (A \vee \div A)$;
(Alt11)    $\vdash_{\textbf{Alt}} (\div \div A \to A)$;
(Alt12)    $\vdash_{\textbf{Alt}} (A \to (\div A \to B))$,

also hold, should behave as the classical negation. So, all we have to do now is to show that (Alt10)–(Alt12) hold in **Ci** if one substitutes $\div$ for the strong negation $\sim$. We could here make use of an auxiliary lemma:

LEMMA 3.43  These are some theorems of **Ci**:
  (i) $\vdash_{\textbf{Ci}} (A \vee \circ A)$;
  (ii) $\vdash_{\textbf{Ci}} (\neg A \vee \circ A)$.

**Proof:** For part (i), observe that, from (Min6), $[\circ A \vdash_{\textbf{Ci}} (A \vee \circ A)]$, and, from (ci1), $[\neg \circ A \vdash_{\textbf{Ci}} A]$, so, once more by (Min6), and transitivity, $[\neg \circ A \vdash_{\textbf{Ci}} (A \vee \circ A)]$. Using the proof by cases one finally concludes that $[\vdash_{\textbf{Ci}} (A \vee \circ A)]$. Part (ii) is similar to (i), but you should now use (ci2).

THEOREM 3.44  The strong negation $\sim$, in **Ci**, is classical.
**Proof:** To check that (Alt10) holds for $\sim$, that is, that $[\vdash_{\textbf{Ci}} (A \vee (\neg A \wedge \circ A))]$, notice that this last schema is equivalent to $[\vdash_{\textbf{Ci}} (A \vee \neg A) \wedge (A \vee \circ A))]$, by positive classical logic, and the latter is provable from (Min10) and LEMMA 3.43(i), using (Min3). Now, (Alt12) is immediate, by the very definition of $\sim$, and to check (Alt11) you might just notice that by reflexivity we have $[\sim \sim A, A \vdash_{\textbf{Ci}} A]$, and from (Alt12) we have $[\sim \sim A, \sim A \vdash_{\textbf{Ci}} A]$; so, using a new form of proof by cases obtained from (Alt10) and (Min8) (as in FACT 3.7), we conclude that $[\sim \sim A \vdash_{\textbf{Ci}} A]$.

So, **Ci** is strong enough to endow its strong negation with all properties of a classical negation. This result has some immediate consequences. For instance, we could use it to show that (Min9) is redundant in **Ci** (and all other logics extending it). Notice, of course, that the two last results did not really need to use the whole positive *classical* logic, but that its *intuitionistic* fragment (which does not contain (Min9)) would have been enough. Confront the following fact with THEOREM 3.3:

FACT 3.45  The schema (Min9): $(A \vee (A \to B))$ is redundant in the axiomatization of **Ci**.
**Proof:** From reflexivity and (Alt12) we have that $[A \vdash_{\textbf{Ci}} A]$ and $[\sim A \vdash_{\textbf{Ci}} (A \to B)]$. But, of course, either $A$ or $(A \to B)$, by (Min6) and (Min7), imply the above schema, $(A \vee (A \to B))$. So, using (Alt10) once more to provide a proof by cases, we are done.

---

[15] There seems to be, at any rate, a widespread mistaken assumption in the literature to that effect (despite the example of intuitionistic negation, strong but not classical). Yet in some other studies, as for instance Batens's [13], note 11, a 'classical' negation, $\div$, in a paraconsistent logic is assumed to be one which is not only strong but it should also be the case that $[\div A \vdash \neg A]$ holds (as in axiom (bun), in the subsection **3.8**). This *is* the case, however, for **bC**'s strong negation $\sim$, but now we know that it is still *not* classical (it just has some kind of intuitionistic behavior). Ten years before, nevertheless, this same author (see [10], page 224) had put things more precisely, and required for that definition that $[\vdash (A \vee \div A)]$ should also be the case.

We shall not list the properties of ~ in **Ci** at this point, but only mention the fact that 'it is a classical negation' when necessary, and then use any property that derives from this fact.

Now, this strong (classical) negation will give us a very interesting result. We already knew that the other binary connectives worked as their classical counterparts, and we were informed above that **Ci** comes also equipped with a negation which works like the classical one; so why don't we use **Ci** to 'talk about classical logic', that is, use **Ci**'s own stuff to reproduce any classical inference? One intuitive procedure to bring forth such an effect would be to pick any classical inference and just substitute any occurrence of a classical negation by an occurrence of a strong negation, and leave the rest as it is. And this indeed works:

**THEOREM 3.46** The following mapping conservatively translates **CPL** inside of **Ci**:
  (t1.1)    $t_1(p) = p$, if $p$ is an atomic formula;
  (t1.2)    $t_1(A \# B) = t_1(A) \# t_1(B)$, if # is any binary connective;
  (t1.3)    $t_1(\neg A) = {\sim} t_1(A)$.
    So, it is the case that $[\Gamma \vdash_{\textbf{CPL}} A] \Leftrightarrow [t_1[\Gamma] \vdash_{\textbf{Ci}} t_1(A)]$.
**Proof:** Given THEOREM 3.44, we know that, by way of the above transformation, a counterpart to **CPL**'s axiomatization can be obtained inside of **Ci**.

**COROLLARY 3.47** We also have a conservative translation of **eCPL** inside of **Ci**.
    Just extend the above mapping by adding:
  (t1.4)    $t_1({\circ}A) = {\circ} {\circ} t_1(A)$.
**Proof:** This comes from the above theorem, **eCPL**'s axiom (ext) and FACT 3.33(i).

The above recursive translation just substitutes one negation for another, thus giving rise to a *grammatically faithful* (cf. [61], chapter X) way of reproducing classical inferences inside of **Ci**, and inside of any other logic deductively stronger than it, as the ones we will be studying below. Of course, other logics may provide yet some other sensible ways of translating classical logic inside of them (see, for instance, COROLLARY 3.62).

To be sure, we already had, in **bC**, a way of reproducing classical inferences (recall THEOREM 3.11), but at that point we had to introduce further premises in our theories —to wit, the premises that some of our propositions were consistent). A natural question which may arise then is whether this was really necessary, given that from THEOREM 3.42 we know that the 'canonical' strong negation of **bC** was not a classical one, and would then not allow the above translations to be performed inside of **bC**, or could it perhaps be the case that all strong negations are indeed strong, but some are stronger than others? This last option is indeed what occurs, for it can be easily shown, if we just recall FACT 2.10(ii), how one can define a classical negation inside of **bC**, despite the weakness of this logic, thus being able to talk about classical logic already inside of the most basic **C**-system we here present:

**THEOREM 3.48** The logic **bC** does have a classical negation.
**Proof:** From (bc1), the axiom that realizes finite gentle explosiveness, and from left-adjunctiveness, we know that $(A \wedge (\neg A \wedge {\circ}A))$ is a bottom particle, for any formula $A$ —let's choose any of these conjunctions and denote it by $\bot$, as usual. Inspired by FACT 2.10(ii) and using the Deduction Metatheorem we then define a new strong negation, $\div$, on **bC** as $\div A \overset{\text{def}}{=} (A \to \bot)$.[16] To check that *this* negation is classical, we just need to prove that $(\div \div A \to A)$ is a theorem of **bC**. To such an end, first note that $[\vdash_{\mathbf{Ci}} (A \lor \div A)]$, given that this is $[\vdash_{\mathbf{Ci}} (A \lor (A \to \bot))]$, a form of axiom (Min9), and this gives us a new form of proof by cases, as in THEOREM 3.44. Next, notice that $[((A \to \bot) \to \bot), (A \to \bot) \vdash_{\mathbf{Ci}} \bot]$, by *modus ponens*, and $[\bot \vdash_{\mathbf{Ci}} A]$ by definition of the bottom particle, but also $[((A \to \bot) \to \bot), A \vdash_{\mathbf{Ci}} A]$. Thus, the new form of proof by cases will immediately give us $[((A \to \bot) \to \bot) \vdash_{\mathbf{Ci}} A]$.

---

[16] This was indeed one of the many 'negations' set forth by Bunder in [30], though this author seems not to have completely understood their properties (see below the subsection **3.8**).

It is easy then to transform THEOREM 3.46 into a grammatically faithful translation of **CPL** already inside of **bC**, but it would be less easy to find a non-trivial analogue of COROLLARY 3.47, the translation of **eCPL**, given that **bC** is already known to have no consistent theorems (recall THEOREM 3.10) —that is, no theorems of the form $\circ A$. As to the status of the two different strong negations presented above inside of the stronger logic **Ci**, one can easily go on to show that:

FACT 3.49  In **Ci** the two strong negations above, $\sim$ and $\div$, are both classical, and are in fact equivalent, in a sense (but not all strong negations are classical in **Ci**).
**Proof:**  That they are both classical is an obvious consequence from THEOREM 3.44 and THEOREM 3.48. To see that they are equivalent in **Ci**, remember, on the one hand, that $[A, \sim\!A \vdash_{\mathbf{Ci}} \bot]$, by definition, and so $[\sim\!A \vdash_{\mathbf{Ci}} (A \to \bot)]$, that is, $[\sim\!A \vdash_{\mathbf{Ci}} \div A]$, by the Deduction Metatheorem. On the other hand, we have both that $[(A \to \bot), A \vdash \bot]$, and thus $[(A \to \bot), A \vdash_{\mathbf{Ci}} \sim\!A]$, and that $[(A \to \bot), \sim\!A \vdash_{\mathbf{Ci}} \sim\!A]$, so the form of proof by cases offered by THEOREM 3.44 will allow us to conclude that $[\div A \vdash_{\mathbf{Ci}} \sim\!A]$. The reader will be right in thinking that all classical negations extending a positive classical basis are equivalent, but it is still the case, nonetheless, that **Ci** can define other strong negations that do not have a classical character, as for instance $(\neg\neg\sim\!A)$ or $(\neg\neg\div A)$. Take a look at [42], our paper on semantics, in the section on **Ci**, to check this claim.

Now, the reader may perhaps think that classically negated propositions in **bC** and **Ci** (especially given the reconstruction of classical inferences inside of these logics that such negations support by way of the above mentioned conservative translations) would be classical enough so as to be consistent propositions themselves, that is, that $\circ \div A$ would be a theorem, for instance, of **Ci**, for some classical negation $\div$. We will now show that this can hardly be the case:

FACT 3.50  Only consistent or inconsistent formulas can themselves be provably consistent in **Ci**. Thus, $(\circ A)$ is a theorem of **Ci** if, and only if, $A$ is of the form $\circ B$, $\bullet B$, $\neg\circ B$ or $\neg\bullet B$, for some $B$.
**Proof:**  On the one hand, we already know from FACT 3.33 and FACT 3.39 that formulas such as $\circ B$ or $\bullet B$, and their variations, are all provably consistent in **Ci**. To see that the converse is also true, consider the following three-valued matrices, such that 0 is the only non-distinguished value and $v(A \land B) = \frac{1}{2}$ if $v(A) \neq 0$ and $v(B) \neq 0$, and $v(A \land B) = 0$, otherwise; $v(A \lor B) = \frac{1}{2}$ if $v(A) \neq 0$ or $v(B) \neq 0$, and $v(A \lor B) = 0$, otherwise; $v(A \to B) = \frac{1}{2}$ if $v(A) = 0$ or $v(B) \neq 0$, and $v(A \to B) = 0$, otherwise; $v(\neg A) = 1 - v(A)$; $v(\circ A) = 1$ if $v(A) \neq \frac{1}{2}$, and $v(\circ A) = 0$, otherwise.

As a consequence of the last result, in particular, formulas of the form $\circ \sim\!A$ and $\circ \div A$ will not be provable in **Ci**, and, from FACT 3.32, we conclude that **Ci** is *not* controllably explosive in contact with (at least some) classically negated propositions. As we shall see, on the other hand, there are many extensions of this logic that do have this property, at least for some particular $A$'s (see FACT 3.66, or FACT 3.76). But what would have happened if we had indeed theorems such as tho ones ruled out above? Let us here allow ourselves some counterfactual reasoning, and ask ourselves about the possible validity of (IpE) in some specific paraconsistent logics, like the extensions of **Ci**, or of some of its fragments (given that we know, from THEOREM 3.35, that (IpE) still does not hold in **Ci**, anyway):

**THEOREM 3.51** (IpE) cannot hold in any paraconsistent extension of **Ci** in which:

(i) $(\circ \div A)$ holds, for some given classical negation $\div$; *or*

(ii) $\neg(A \wedge \neg A)$ or $\neg(\neg A \wedge A)$ hold; *or*

(iii) $[(\neg A \vee \neg B) \vdash \neg(A \wedge B)]$ hold; *or*

(iv) $[\neg(A \wedge B) \vdash (\neg A \vee \neg B)]$ hold.

  (IpE) cannot hold in any paraconsistent extension of **bC** in which:

(v) $[\neg(A \rightarrow B) \vdash (A \wedge \neg B)]$ hold.

  (IpE) cannot hold in any adjunctive paraconsistent extension of $C_{min}$ in which:

(vi) both $[(A \wedge B) \vdash \neg(\neg A \vee \neg B)]$ and $[\neg(\neg A \vee \neg B) \vdash (A \wedge B)]$ hold.

  (IpE) cannot hold in any adjunctive paraconsistent logic in which:

(vii) both $\neg(A \wedge \neg A)$ and $[(A \wedge \neg A) \dashv\vdash \neg\neg(A \wedge \neg A)]$ hold.

**Proof:** For part (i), given that $\div$ is a classical negation we can then assume $[A \dashv\vdash \div\div A]$ to hold. Now, if (IpE) were valid one could conclude, in particular, that $[\circ A \dashv\vdash \circ \div\div A]$, and given that $(\circ \div\div A)$ is a theorem of this logic extending **Ci**, by hypothesis, one would infer $\circ A$ as a theorem, but this is (ext), exactly the axiom that is lacking to make **Ci** collapse into **eCPL**. This generalizes a similar argument to be found in [107], Theorem 9. To check part (ii), recall that, in **Ci**, $[\bullet A \dashv\vdash (A \wedge \neg A)]$, and (IpE) would then give $[\circ A \dashv\vdash \neg(A \wedge \neg A)]$, and we are again left with the theorem $\circ A$, as in part (i). For part (iii), recall that $(\neg A \vee \neg\neg A)$ is a theorem already of $C_{min}$, and the problem reduces then to part (ii). For parts (iv) and (v), we will just show that $[(\neg \div A) \vdash A]$ is obtained, and so we may conclude that controllable explosion occurs in contact with $\div A$. Given that $[\neg(A \wedge B) \vdash (\neg A \vee \neg B)]$ holds, consider the strong negation $\sim A \overset{\text{def}}{=} (\neg A \wedge \circ A)$, for which one would immediately obtain $[\neg\sim A \vdash (\neg\neg A \vee \neg \circ A)]$, and so, from (Min11), (ci1) and (Min8), we get $[\neg\sim A \vdash A]$. Given that $[\neg(A \rightarrow B) \vdash (A \wedge \neg B)]$ holds, pick up $\div A \overset{\text{def}}{=} (A \rightarrow \bot)$, and, from (Min4), we have that $[\neg \div A \vdash A]$. For part (vi), given once more that $(\neg A \vee \neg\neg A)$ is a theorem of $C_{min}$, (IpE) would give us $[\neg(\neg A \vee \neg\neg A) \dashv\vdash \neg(\neg B \vee \neg\neg B)]$, and the rules that we here assume give us $[(A \wedge \neg A) \dashv\vdash (B \wedge \neg B)]$, so, by adjunction, we conclude in particular that $[A, \neg A \vdash B]$. This is the main result in Béziau's [21]. Finally, for part (vii), (IpE) would give us $[\neg\neg(A \wedge \neg A) \dashv\vdash \neg\neg(B \wedge \neg B)]$, and so $[(A \wedge \neg A) \dashv\vdash (B \wedge \neg B)]$, and we are in the same situation as in (vi). This is a stronger version of the main result in Béziau's [22], where actually $[A \dashv\vdash \neg\neg A]$ was assumed, instead of $[(A \wedge \neg A) \dashv\vdash \neg\neg(A \wedge \neg A)]$. Of course, a similar version of this last result arises if one just uniformly substitutes $(A \wedge \neg A)$ for $(\neg A \wedge A)$ in its statement. Notice that the rules mentioned in parts (iii) to (vi) had already shown up as the items (x), (ix), (vii), (xii) and (xi) of THEOREM 3.18.

So far we have some negative results about the validity of (IpE) in some possible paraconsistent extensions of **bC** or **Ci**, but are there paraconsistent extensions of these logics in which (IpE) *does* hold? In the search for an answer, one could start by testing the compatibility of the addition, to those logics, of at least one of the following rules of deduction, (RC): $[(A \vdash B) \Rightarrow (\neg B \vdash \neg A)]$ or (EC): $[(A \dashv\vdash B) \Rightarrow (\neg B \vdash \neg A)]$ (see the subsection **3.3**, where these were argued to be enough for the consistency-less fragment of our language), and also of at least one of the following:

$$\forall A \ \forall B \ [(A \vdash B) \Rightarrow (\circ A \vdash \circ B)]; \tag{RO}$$

$$\forall A \ \forall B \ [(A \dashv\vdash B) \Rightarrow (\circ A \vdash \circ B)]. \tag{EO}$$

We have already shown, in FACT 3.36, that (RC) cannot be added to **Ci** without collapsing into classical logic. We can now actually show more:

**FACT 3.52**  In extensions of **Ci**, the validity of (EC) also guarantees (EO).

**Proof:** From $[A \dashv\vdash B]$ we conclude, by (EC), that $[\neg A \dashv\vdash \neg B]$. From these two sentences, by positive logic, we conclude that $[(A \wedge \neg A) \dashv\vdash (B \wedge \neg B)]$, but from FACT 3.16(ii) and FACT 3.29(i) we know that $[\neg \circ C \dashv\vdash (C \wedge \neg C)]$, and so we have that $[\neg \circ A \dashv\vdash \neg \circ B]$. Finally, from FACT 3.34(iv), we have that $[\circ A \dashv\vdash \circ B]$.

The problem of finding paraconsistent extensions of **Ci** in which (IpE) holds reduces then to the problem of finding out if (EC) can be added to this logic without losing the paraconsistent character. We suspect this can be done, but shall leave it as an open problem at this point. As to extensions of **bC**, on the other hand, we can already present a (very partial) result:

**THEOREM 3.53**  There are fragments of **eCPL** extending **bC** in which (IpE) holds.

**Proof:** We already know, from COROLLARY 3.22, that (IpE) does not hold for **bC** as it is. It suffices now to show that the addition of the rules (EC) and (EO) to **bC** may still originate a paraconsistent fragment of (extended) classical logic, once these rules are evidently enough to ensure (IpE). To such an end, one may simply make use of the following matrices by Urbas ([107], Theorem 8):

| ∧ | **1** | 6/7 | 5/7 | 4/7 | 3/7 | 2/7 | 1/7 | **0** |
|---|---|---|---|---|---|---|---|---|
| **1** | 1 | 6/7 | 5/7 | 4/7 | 3/7 | 2/7 | 1/7 | 0 |
| 6/7 | 6/7 | 6/7 | 3/7 | 3/7 | 3/7 | 3/7 | 2/7 | 0 |
| 5/7 | 5/7 | 3/7 | 3/7 | 5/7 | 1/7 | 3/7 | 0 | 1/7 |
| 4/7 | 4/7 | 3/7 | 5/7 | 1/7 | 4/7 | 0 | 2/7 | 1/7 |
| 3/7 | 3/7 | 3/7 | 1/7 | 4/7 | 3/7 | 0 | 3/7 | 0 |
| 2/7 | 2/7 | 3/7 | 3/7 | 0 | 0 | 0 | 2/7 | 0 |
| 1/7 | 1/7 | 0 | 1/7 | 1/7 | 1/7 | 0 | 0 | 1/7 |
| **0** | 0 | 0 | 0 | 0 | 0 | 0 | 0 | 0 |

| ∨ | **1** | 6/7 | 5/7 | 4/7 | 3/7 | 2/7 | 1/7 | **0** |
|---|---|---|---|---|---|---|---|---|
| **1** | 1 | 1 | 1 | 1 | 1 | 1 | 1 | 1 |
| 6/7 | 1 | 6/7 | 1 | 1 | 6/7 | 6/7 | 1 | 6/7 |
| 5/7 | 1 | 1 | 5/7 | 1 | 5/7 | 1 | 5/7 | 5/7 |
| 4/7 | 1 | 1 | 1 | 4/7 | 1 | 4/7 | 4/7 | 4/7 |
| 3/7 | 1 | 6/7 | 5/7 | 1 | 3/7 | 3/7 | 6/7 | 3/7 |
| 2/7 | 1 | 6/7 | 1 | 4/7 | 6/7 | 2/7 | 4/7 | 2/7 |
| 1/7 | 1 | 1 | 5/7 | 4/7 | 5/7 | 4/7 | 1/7 | 1/7 |
| **0** | 1 | 6/7 | 5/7 | 4/7 | 3/7 | 2/7 | 1/7 | 0 |

| → | **1** | 6/7 | 5/7 | 4/7 | 3/7 | 2/7 | 1/7 | **0** |
|---|---|---|---|---|---|---|---|---|
| **1** | 1 | 6/7 | 5/7 | 4/7 | 3/7 | 2/7 | 1/7 | 0 |
| 6/7 | 1 | 1 | 5/7 | 4/7 | 5/7 | 4/7 | 1/7 | 1/7 |
| 5/7 | 1 | 6/7 | 1 | 4/7 | 6/7 | 2/7 | 4/7 | 2/7 |
| 4/7 | 1 | 6/7 | 5/7 | 1 | 3/7 | 3/7 | 6/7 | 5/7 |
| 3/7 | 1 | 1 | 1 | 1 | 4/7 | 1 | 4/7 | 4/7 |
| 2/7 | 1 | 1 | 5/7 | 1 | 5/7 | 1 | 5/7 | 5/7 |
| 1/7 | 1 | 6/7 | 1 | 1 | 6/7 | 6/7 | 1 | 6/7 |
| **0** | 1 | 1 | 1 | 1 | 1 | 1 | 1 | 1 |

| | ¬ |
|---|---|
| **1** | 0 |
| 6/7 | 5/7 |
| 5/7 | 2/7 |
| 4/7 | 3/7 |
| 3/7 | 4/7 |
| 2/7 | 5/7 |
| 1/7 | 1 |
| **0** | 1 |

where 1 is the only distinguished value, and add to these a matrix for $\circ$ that is constant and equal to 0. It is straightforward to check that the above matrices validate all axioms and rules of **bC** plus the two rules above, while formulas such as $(A \rightarrow (\neg A \rightarrow B))$ and $\neg(A \wedge \neg A)$ are still not validated by them.[17]

### 3.8 Beyond Ci: The dC-systems.

We have now come closer to the more orthodox approach to paraconsistent logics that the reader will find in the field, which *does* identify contradictoriness and inconsistency. All logics that we will be studying from here on, and, we argue, all logics of formal inconsistency presented in the literature so far, do not distinguish between these two notions. But this does *not* mean, the reader should be aware, that one can simply dispense with the new operators that have allowed us, so far, to talk about a formula being consistent or inconsistent, for, we remember from FACT 3.16, even though $[\bullet A \dashv \vdash_{\mathbf{Ci}} (A \wedge \neg A)]$ holds, $[\circ A \dashv \vdash_{\mathbf{Ci}} \neg(A \wedge \neg A)]$, for instance, does not hold (in fact, $[\nvdash_{\mathbf{Ci}} \neg(A \wedge \neg A)]$)! Suppose then that we construct now the logic **Cil** exactly by adding to **Ci** (that is, (Min1)–(Min11), (bc1), (ci), (MP), plus the definition of $\circ$ in terms of $\bullet$) the following 'missing' axiomatic rule:

(cl)      $\neg(A \wedge \neg A) \vdash \circ A$.          'If $\neg(A \wedge \neg A)$ is the case, then $A$ is consistent'

Much confusion has been raised around this particular formula, $\neg(A \wedge \neg A)$. Recall, from THEOREM 3.15(ii) and THEOREM 3.30, that it was not a theorem of our previous logics, **bC** or **Ci**. Now, of course, we can immediately conclude even more:

**FACT 3.54** No paraconsistent extension of **Cil** can have $\neg(A \wedge \neg A)$ as a theorem.
**Proof:** If so, axiom (cl) would give us consistency, thus ruining paraconsistency.

Nonetheless, some other paraconsistent extensions of **Ci**, such as **LFI1** (see its matrices in the subsection **2.4** or at FACT 3.29, and see its axiomatization in THEOREM 3.69), *do* have this as theorem (and, as a consequence of THEOREM 3.51(ii), they must lack a full replacement theorem).

The attribution of a privileged status to the formula $\neg(A \wedge \neg A)$, using it to express consistency inside some paraconsistent logics, stems from the early requisites put forward by da Costa on the construction of his famed calculi $C_n$:

**dC[i]**   in these calculi the principle of non-contradiction [sic], in the form $\neg(A \wedge \neg A)$, should not be a valid schema;
**dC[ii]**   from two contradictory formulas, $A$ and $\neg A$, it would not in general be possible to deduce an arbitrary formula $B$;
**dC[iii]**   it should be simple to extend these calculi to corresponding predicate calculi (with or without equality);
**dC[iv]**   they should contain the most part of the schemas and rules of the classical propositional calculus which do not interfere with the first conditions.

While the requisite **dC[ii]** is nothing but the very definition of a paraconsistent logic (recall the subsection **2.2**), **dC[iii]** is simply a claim for extensions of these logics to higher-order calculi (that we will *not* explore here, reiterating instead the popular and still powerful argument of paraconsistentists about the fact that most, if not all, innovations of paraconsistent logic can already be met at the propositional level), and **dC[iv]** is indeed somewhat vague, having received much attention and many di-verse interpretations from several researchers (our own proposal on its interpretation will be found in the subsection **3.11**), the requisite **dC[i]** is in fact the one to blame for the confusion we were talking about. First of all, under our present perspective, to call the formula $\neg(A \wedge \neg A)$ 'principle of non-contradiction' is quite misleading, not only because this would bring us, as a side effect, to commit to very particular interpretations for the negation of a proposition, the conjunction of contradictory propositions, and the negation of this conjunction, but ascribe to us as well a very particular interpretation for the consistency connective (however, if —and only if— you are working in the context of some specific consistent logics, such as classical logic itself, or intuitionistic logic, we admit that this designation can indeed make sense). Evidently, FACT 3.54 is then just a consequence of such a contract.

[17] Notice, nevertheless, that all that is proved here is that there certainly exist fragments of classical logic extending **bC** for which (IpE) holds. But the matrices above do not fulfill, of course, our requisite for defining a paraconsistent logic (namely, disrespecting (PPS)), so that the question is still left open as to whether there are *paraconsistent* such extensions of **bC**! (With thanks to Dirk Batens for calling our attention to that.)

Now, if, on the one hand, some authors have questioned the validity of $\neg(A \wedge \neg A)$ in the context of a paraconsistent logic,[18] on the other hand the construction of paraconsistent logics in which this formula does not hold has also been rather criticized since then, and for various reasons. Some of these criticisms unfold from or link to, by and large, still that same understandable and widespread confusion between the 'principle of non-contradiction' and the aims of paraconsistent logic, namely to avoid the 'principle of explosion' instead (see **dC[ii]**, and our subsections **2.1** and **2.2**) —but these more or less loose arguments can hardly be recast under our present formal definitions of those principles. A slightly more elaborate argumentation appears in Routley & Meyer's [100], where the authors are looking for some formalization of dialectical logic, and they claim to that effect not only that $\neg(A \wedge \neg A)$ is 'usually' a theorem of the 'entailment systems' that they have examined, but also that this does not conflict with other logical truths of those dialectical systems (in their words, this does not generate any 'intolerable tensions which destroy any prospect of a coherent logic'), even though these systems do have contradictory theorems (that are to be understood as 'synthetic a priori'), and validate adjunction. Moreover, they maintain that 'the orthodox Soviet position appears to retain $\neg(A \wedge \neg A)$ as a thesis', and they want to deal with it.[19] Now, none of the logics we study here are dialectical, in the sense of disrespecting the Principle of Non-Contradiction and actually proving contradictory formulas (recall the subsection **2.2**), and so the last critique above, in any case, simply falls idle.

---

[18] For instance, Béziau's [23], section 2.3, argues that, from a *philosophical* standpoint, it is hard to reconcile the validity of $\neg(A \wedge \neg A)$ and an intuitive interpretation for the negation symbol of a paraconsistent logic; as to the *technical* aspect of his criticism, it seems to consists basically of a consequence of THEOREM 3.51(vii) and the wish to obtain both adjunctiveness and the validity of (IpE).

[19] They also say some other things which seem a bit weird. First, that the 'non-orthodox' systems not containing $\neg(A \wedge \neg A)$ as a theorem are all *weaker* dialectical logics (as if the logics were all linearly ordered by strength!). Secondly, they insist on calling the formula $\neg(A \wedge \neg A)$ 'Aristotle's principle of non-contradiction', and after formally presenting their dialectical logics they argue that this formula is 'correct, both in syntactical and semantical formulations': *syntactically* correct because 'it is a theorem, hence valid, hence true' —despite the seemingly naive *petitio principii* brought therein; and *semantically* correct because '[one of its historical formulations] asserts that no statement is both true and false', and this feature, in the case of their logics, is supposed to be 'guaranteed by the bivalent features of the semantics' —now this is surely a mistake, for what guarantees this fact can only be the functional (rather than relational) character of their proposed interpretation, but in any case this last argument by these authors, even if they were kind enough to clear up the somewhat obscure relation of it with the first one, should hardly be accepted as a justification, given that any associated semantics provided to a consequence relation of a given logic is barely *circumstantial*, in a sense, and can often be recast in many apparently non-equivalent ways (if you're not happy with a particular semantics, you can always *look for another one*). A similar criticism of these points has been made before in Batens's [10], section 9.

But let us first explore some consequences of the new axiom (cl), before really questioning it any deeper, or looking for substitutes. The main and most far-reaching consequence is the following:

**THEOREM 3.55** In **Cil** we can define the inconsistency operator as $\bullet A \stackrel{\text{def}}{=} (A \wedge \neg A)$ (from which the consistency operator will be defined as $\circ A \stackrel{\text{def}}{=} \neg(A \wedge \neg A)$).

**Proof:** It can immediately be seen, from the above definitions, that the axioms (bc1), (ci) and (cl) will still hold if one just substitutes all occurrences of the operators $\bullet$ and $\circ$ by their new definitions.

**COROLLARY 3.56** Given a theorem $B$ of **Cil** we can substitute all occurrences of $\bullet$ and of $\circ$ in its subformulas according to the above definitions.

**Proof:** Recall from THEOREM 3.41 that all axioms of **Ci** can be written just with the use of $\bullet$, substituting $\circ$ for $\neg\bullet$, or just with the use of $\circ$, substituting $\bullet$ for $\neg\circ$. In the first case, where we have only '$\bullet$'s, the above theorem permits us to rewrite in **Cil** the proof of $B$ using $(A \wedge \neg A)$ in the place of each formula $\bullet A$ that appears, and using $\neg(A \wedge \neg A)$ in the place of each formula $\neg\bullet A$ that occurs in the proof. In the second case, and for the same reason, we may rewrite the proof of $B$ using $\neg(A \wedge \neg A)$ in the place of each formula $\circ A$ that appears, and using $(A \wedge \neg A)$ in the place of each formula $\neg\circ A$ that occurs in the proof (you may in this last part also wish to take FACT 3.40 and (Min11) once more into consideration).

The above results are structurally similar to those of THEOREM 3.41 and its consequences, where a restricted form of replacement was obtained for the operators $\bullet$ and $\circ$, and we have seen that each one of them could be substituted in **Ci** by the negation of the other. But now we know more, we know that we can simply dispense with the operators $\bullet$ and $\circ$, substituting each formula $\bullet A$ and each formula $\neg\circ A$ for the formula $(A \wedge \neg A)$, each formula $\circ A$ and each formula $\neg\bullet A$ for the formula $\neg(A \wedge \neg A)$. This brings us to the definition of a particular subclass of the **C**-systems that we will call **dC**-*systems*, such as the **C**-systems in which $\bullet$ and $\circ$ can be defined in terms of the other connectives. **Cil** is the first example of a **dC**-system that we here consider; before presenting other examples let us point out some consequences of this last result for **Cil**. It is easy to see, for instance, that the restricted forms of contraposition presented in FACT 3.17 for **bC** and in FACT 3.34 for **Ci**, as well as the forms of *reductio* presented in FACT 3.14 for **bC** and in FACT 3.31 for **Ci**, and the forms of controllable explosion presented in this last fact, together with the fundamental fact relating controllable explosion and consistency in FACT 3.32, all have new versions in **Cil**, if we just change each occurrence of $\circ$ and $\bullet$ for their definitions in THEOREM 3.55. We can also update THEOREM 3.11 with yet another way of reproducing classical inferences inside of **Cil** by the addition of the appropriate premises to its theories (namely the addition of a finite number of formulas of the form $\neg(A \wedge \neg A)$, the formula that in the present circumstances represents the consistency of $A$). Analogously to what we did in FACT 3.40, we can now prove that $[(A \wedge \neg A) \vdash_{\mathbf{Cil}} \neg\neg(A \wedge \neg A)]$, even though $[A \vdash_{\mathbf{Cil}} \neg\neg A]$ still does not hold.

Now, if FACT 3.50 has provided us with a very precise characterization of the consistent theorems in **Ci**, which turned out to be only consistent or inconsistent formulas themselves, namely the ones appearing in FACT 3.33 and FACT 3.39, we may now use those same results to conclude that formulas such as $\circ(A \wedge \neg A)$, and, by FACT 3.38(i), also $\circ\neg(A\wedge\neg A)$, are theorems of **Cil**. These last theorems have raised yet some other protests in the literature. For instance, Sylvan (cf. [105]) claims that the fact that such a logic validates some 'unjustifiable' intuitionistically invalid theorems, together with the validity of $\circ(A\wedge\neg A)$ 'defeats certain paraconsistent objectives' (too bad that he did not proceed to clear up which objectives were these…). This echoes, in one way or another, to a common, and entirely well-founded, criticism, which has been raised by various authors, both to the fact that **C**-systems such as those we have been studying do maintain the whole of positive classical logic and to the fact that many of them (but not all!) are in fact **dC**-systems, and come up with rather particular definitions for the consistency operator. However, both these aspects can be easily varied and experimented. We have here, by a matter of simplicity, set up an investigation of **C**-systems based on classical propositional logic, but it is clear that other approaches may be tackled, by the investigation of **C**-systems based on relevance logic, or intuitionistic logic, as soon as some paradoxes of relevance, or paradoxes raising from some non-constructive assumption, are decided to be avoided. This is clearly not, however, a problem of *paraconsistency* as we have it, but a further (interesting) problem which can be added to it.

Some **dC**-systems based on intuitionistic logic have in fact been defined and studied, for instance, in Bunder's [29]. **B₁**, the stronger logic of the main hierarchy of calculi proposed by the author of that paper, is obtainable simply by dropping the axioms (Min9), (Min10) and (Min11) out of **Cil**, while adding to the resulting logic the axiom:

(bun)     $(A\rightarrow(\circ B\wedge(B\wedge\neg B)))\vdash\neg A$.    'If $A$ implies a bottom particle, then $\neg A$ is the case'

It is more or less clear that the deletion of the above axioms from **Cil** will give the resulting logic a kind of intuitionistic behavior, and that the addition of (bun) cannot recover any of the classical properties which were lost. Despite of this, most, if not all, other claims that the author advances about this logic seem to be mistaken. He conjectures, for instance, that the strong negation defined by the antecedent of (bun), by setting $\dot{\div}A \overset{\text{def}}{=} (A\rightarrow(\circ B\wedge(B\wedge\neg B)))$, for some formula $B$, *is not* a classical negation, and *continues not to be* a classical negation even if one adds back to **B₁** the axioms (Min10) and (Min11). This is wrong, for we know from FACT 3.45 that in this case the axiom (Min9) turns to be provable, and so we obtain a logic at least as strong as **Ci** (plus (bun), if this would make any difference), but than we remember from FACT 3.49 that $\dot{\div}$ is *indeed* a classical negation in **Ci** (and even in **bC**, as we saw in THEOREM 3.48). The author then claims, and purports to prove, that his **B₁** is *not* a subsystem of da Costa's logic $C_1$, which, as we will see in the subsection **3.10**, is simply an extension of **Cil**. Once more he is wrong, and for the very same reason —FACT 3.49 shows us once more that (bun) is evidently provable in **Ci**, and so already **Ci** (and consequently $C_1$) extends **B₁** and all the other weaker calculi proposed by Bunder. From that point on, all of the remaining remarks made by this author on the comparison of his calculi with the ones proposed by da Costa falls apart. The only point remaining from those calculi, therefore, is that of constituting **dC**-systems based on intuitionistic, rather than classical, logic. But the author did not even try to study them any deeper, looking for instance for interpretations for these calculi!

In another paper, Bunder unwittingly produced an even bigger mistake. Starting from the reasonable idea of looking for other formulations for da Costa's version of *reductio* (that is, (RA0): $[\circ B,\ (A\rightarrow B),\ (A\rightarrow\neg B)\vdash\neg A]$, in the subsection **3.2**), the author simply proposes (in [31]) to change $\circ B$ for $\circ A$ in that formula, asserting that 'there seems to be no particular reason why, in (RA0), the $B$ has a restriction, rather than the $A$'. In this case, however, we can use our THEOREM 3.13 again to see that this proposal would be equivalent to the addition of $[\circ B,\ A,\ \neg A \vdash\ C]$ to $C_{min}$, which is clearly absurd, for it would be trying to express the consistency of $A$ by way of some foreign formula $B$! The author then claims that the 'paraconsistent' calculus $D_1$ (again, the 'strongest' one of a hierarchy $D_n$) that he obtains by adding this last formula as a new axiom to $C_1$ is 'strictly stronger than $C_n$', the calculi of da Costa's hierarchy, and purports to prove some facts about them. Once more these facts turn out to be mistaken, and this is easy to see if one remembers that $\circ\circ D$ is a theorem of $\mathbf{Ci}$ (see our FACT 3.33(i), or his Theorem 5), and so we are left with $[A, \neg A \vdash\ C]$, for any $A$ and $C$, and explosiveness is back. So, the author was actually right about his calculi being extensions of the calculi $C_n$, but only because they all *collapse* into classical logic, after all…[20]

As advanced above (and we shall confirm this below, in the subsection **3.10**), that the identification of consistency with the formula $\neg(A\wedge\neg A)$ was exactly what was done in da Costa's calculus $C_1$, which in fact just adds to $\mathbf{Cil}$ some more axioms to deal with the 'propagation of consistency' from simpler to more complex formulas. Now, many authors have criticized this identification —'there is nothing sacrosanct about the original definition of this schema as $\neg(A\wedge\neg A)$', says Urbas in [107]—, or else its consequences, as we have mentioned above (as the 'anomalies' described, for instance, in Sylvan's [105]). One of the most unexpected consequences of this identification, in fact, has already been pointed out in Urbas's [107], Theorem 4, and was hinted above in our subsection **3.2**:

THEOREM 3.57  In **Cil** the consistency of the formula $A$ can be expressed by the formula $\neg(A\wedge\neg A)$, but *not* by the formula $\neg(\neg A\wedge A)$. In fact, one can even add $\neg(\neg A\wedge A)$ to **Cil**, but not $\neg(A\wedge\neg A)$, without this logic losing its paraconsistent character.

**Proof:** That consistency is so expressed in **Cil** and that $\neg(A\wedge\neg A)$ cannot be added to it are simply consequences, respectively, of COROLLARY 3.56 and FACT 3.54. But while it is easy to see that $[\neg(A\wedge\neg A)\vdash_{\mathbf{Cil}} \neg(\neg A\wedge A)]$, the converse of this does *not* hold, as we see from the matrices in THEOREM 3.21 (those three-valued matrices are in fact a much simpler way of checking the same result for which Urbas has used six-valued ones). Notice, in particular, that these matrices in fact also validate the formula $\neg(\neg A\wedge A)$.

The above phenomenon is a bit tricky, and has actually fooled people working with the calculi $C_n$ for perhaps too long a time (see, for instance, [76], note 6, ch.2, p.49, or else [42]). Let us then consider the following alternatives to the 'levo-'axiom (cl):

(cd)  $\neg(\neg A\wedge A)\vdash \circ A$;

(cb)  $(\neg(A\wedge\neg A)\vee\neg(\neg A\wedge A))\vdash \circ A$.

---

[20] By the way, given the above considerations, one of Bunder's main results about these systems (besides the supposed proof about all the calculi $D_n$ constituting different systems), objected to show that the calculi $D_n$ do not satisfy (IpE) (namely the Theorem 10 that closes his [31]) evidently must fail. It is easy, in fact, to find counter-examples for the validity of the formula $[\neg(A\wedge\neg A),\ (A\to B),\ (A\to\neg B)\vdash \neg A]$ in the matrices that he proposes (picking them up from Urbas's [107], who have used them correctly, in the case of the $C_n$ systems): just choose $v(A)\in\{0,\ 3\}$ and $v(B)=1$.

Evidently, the addition to **Ci** of the 'dextro-'axiom (cd), instead of the axiom (cl), would give us this logic **Cid** which has exactly the same qualities and defects as **Cil**, but which would singularize the formula $\neg(\neg A \wedge A)$ as much as the formula $\neg(A \wedge \neg A)$ has been previously singularized by **Cil**. The addition of (cb), instead, defining the logic **Cib**, would assure to both $\neg(A \wedge \neg A)$ and $\neg(\neg A \wedge A)$ the same status, and this would fix, for instance, the famed asymmetry in THEOREM 3.21(iii) (but not the ones in parts (i) and (ii)). Logics having (cb) instead of (cl) have already been studied (see [36] or [76]), but it should be noted that these still suffer from some anomalies related to the definition of consistency in terms of some operation over a conjunction of contradictory formulas. In fact, if the logics having (cb) as an axiom do identify the two formulas above, they do not necessarily identify these with some other formulas such as $\neg(A \wedge (A \wedge \neg A))$, or $\neg((A \wedge \neg A) \wedge A)$, for instance, even though all of the formulas $(A \wedge \neg A)$, $(\neg A \wedge A)$, $(A \wedge (A \wedge \neg A))$ and $((A \wedge \neg A) \wedge A)$ are equivalent on any **C**-system based on classical logic. All of this will have, of course, deep consequences when we go on to provide semantics to these logics, as the reader will see in [42]). Perhaps a good way of fixing all of this at once is by the addition of a new 'global' axiomatic rule to such **dC**-systems, such as:

(cg)     $(B \leftrightarrow (A \wedge \neg A)) \vdash (\neg B \leftrightarrow \neg(A \wedge \neg A))$,

or else the weaker deduction rule:

(RG)     $[B \dashv\vdash (A \wedge \neg A)] \Rightarrow [\neg B \dashv\vdash \neg(A \wedge \neg A)]$.

Logics having such rules are yet to be more deeply investigated. In one way or another, it is clear that the mere addition of such rules is not enough to remedy the whole of THEOREM 3.21 (but compare the following result to the proposal by Mortensen, in the subsection **3.12**). Indeed, similarly to what had happened in COROLLARY 3.22 and in THEOREM 3.35:

**THEOREM 3.58**  (IpE) does not hold for **Cib** plus (cg) or (RG).
**Proof:**  Check for instance that THEOREM 3.21(iii) still holds, that is, that formulas such as $\neg(A \vee B)$ and $\neg(B \vee A)$ are still not equivalent, once more by way of the same matrices and distinguished values as in THEOREM 3.15(ii), but now changing only the value of $(1 \vee \frac{1}{2})$ from 1 to $\frac{1}{2}$ (and leaving the value of $(\frac{1}{2} \vee 1)$ as it is, equal to 1).

It is also noteworthy that da Costa in fact proposed not just one definition of consistency (the one above, from the calculus $C_1$), but considered instead the possibility of having weaker and weaker logics (see [49] or [50]), modifying the requirement for consistency in such a way as to produce an infinite number of logics at once (he acted more with an illustrative than with a practical purpose, but that manoeuvre has, in one way or another, produced some permanent impression). The idea is simple, namely that of having, for each $C_n$, $0 \leq n < \omega$, more and more premises to be fulfilled in order to guarantee consistency. In the case of $n = 1$ we already know that $\circ A$ (da Costa denoted it $A^\circ$) abbreviated the formula $\neg(A \wedge \neg A)$, for $1 < n < \omega$ it was taken to be $A^{(n)}$, where this abbreviation was recursively defined by first setting $A^n$, $0 \leq n < \omega$, as $A^0 \stackrel{\text{def}}{=} A$ and $A^{n+1} \stackrel{\text{def}}{=} (A^n)^\circ$, and then setting $A^{(n)}$, $1 \leq n < \omega$, as $A^{(1)} \stackrel{\text{def}}{=} A^1$ and $A^{(n+1)} \stackrel{\text{def}}{=} A^{(n)} \wedge A^{n+1}$. In other words, each of da Costa's **dC**-systems was defined by exactly the same axioms, changing only the definition of $\circ A$ in each case for $A^{(n)}$, for each given $n$.[21] It is clear in this way, if one really feels inclined to do it for some reason or another, that for each **dC**-system one can go on to multiply it into an infinite number of (in principle, distinct) **dC**-systems, applying the same strategy above. Of course, the same asymmetries already pointed out in each case of **Cil**, **Cid** and **Cib** (THEOREM 3.57), may still apply in each case in appropriate forms, and the theorems obtained from these systems must also be modified, in each case, according to the specific definition of consistency brought therein.

**3.9 The opposite of the opposite.** Having been introduced to the **dC**-systems, a particular class of **C**-systems that can dispense with the use of the operators ∘ and •, we shall not, nevertheless, dedicate ourselves in what follows exclusively to the study of **dC**-systems. All the logics we will present from this point on are still bound to be **C**-systems extending **Ci**, but only by chance will they turn out to be **dC**-systems as well. What we will consider in this subsection is the addition of the following axiomatic rule for 'expansion' of negations, converse to (Min11): ($\neg\neg A \to A$):

(ce)    $A \vdash \neg\neg A$.

Let **Cie** be the logic obtained by the addition of (ce) to **Ci** (recall, from THEOREM 3.30, that this addition is *not* redundant). In the subsections **3.1** and **3.2** we have learned about the role played by (Min11) in **bC** and the logics which extend it (recall THEOREM 3.13), and suggestions were made as to the reasons why da Costa has introduced (Min11) in his first paraconsistent calculi as a dual substitute to (ce), present in intuitionistic logic, as much as (Min10): ($A \lor \neg A$) was intended to be the dual substitute to (PPS), the explosiveness (or *reductio*) that is lost by all paraconsistent calculi. Despite this, qualified forms of both (PPS) and (ce) are retained by **bC**, in the form of the rule (bc1): [∘$A$, $A$, $\neg A \vdash B$], and the rule [∘$A$, $A \vdash \neg\neg A$], this last one being a rule of **bC** that comes immediately from (bc1) and FACT 3.14(iii). Now, it happens that only (PPS), but not (ce), is a problem of paraconsistent logic, as we put it, and, as far as we know, (ce) was only avoided by da Costa in his first calculi (see [49] and [50]), in spite of his manifest intention, on his requisite **dC[iv]** (see the last subsection), to maintain 'most rules and schemas of classical logic not conflicting with the other requisites', because there seemed to be some apprehension about the addition of (ce) leading us back to classical logic, **CPL**, or perhaps making us just lose the paraconsistency character of our logics, after all. It is, however, very easy to see that this is not the case. Indeed:

THEOREM 3.59 (tPS): ($A \to (\neg A \to B)$) is not provable by **Cie**.

---

[21] One should observe, however, that the definition of the schema $A^{(n)}$ proposed in da Costa's foundational work, [49], in reality does *not* coincide with the definitions to be found in other studies in the literature (such as the well-known da Costa's [50], or da Costa & Alves's [53]), and those that we adopt here. Indeed, on page 16 of [49] the reader will find the following definition, setting $A^{(1)} \overset{\text{def}}{=} A^\circ$ and $A^{(n+1)} \overset{\text{def}}{=} A^{(n)} \land (A^{(n)})^\circ$. If one follows this last definition, one ought to conclude, for instance, that $A^{(3)}$ is to denote the formula $A^\circ \land A^{\circ\circ} \land (A^\circ \land A^{\circ\circ})^\circ$, while the definition we have presented above would give instead $A^\circ \land A^{\circ\circ} \land A^{\circ\circ\circ}$. It is easy to see, however, if one just makes use of any of the semantics and decision procedures that have been associated to the calculi $C_n$ that these two formulas are *not* equivalent in each $C_n$ (see [36] for the semantics of a slightly stronger version of these calculi —in the case of $n=1$, for instance, axiom (cb) is used instead of (cl), and axiom (ce), which appears below, was also added; or else go to [76] and [74], for the original versions).

**Proof:** Use the matrices of **LFI1** again, as in FACT 3.29, or else the matrices of $\mathbf{P}^1$, as in THEOREM 3.30 —but in this last case you must change the matrix of negation, setting the value of $\neg\frac{1}{2}$ as $\frac{1}{2}$, instead of 1. In fact, it is to be remarked that this modification on the matrices of $\mathbf{P}^1$ in fact originates a new and interesting maximal three-valued paraconsistent logic, $\mathbf{P}^2$. These three logics, **LFI1**, $\mathbf{P}^1$ and $\mathbf{P}^2$, will be studied in their own right in the subsection **3.11**, as members of a larger family of similar logics.

Some of the immediate and main syntactical results obtained by the logic **Cie** are:

FACT 3.60 **Cie** proves the following:

  (i) $\circ\neg A \vdash_{\mathbf{Cie}} \circ A$;

  (ii) $\bullet A \vdash_{\mathbf{Cie}} \bullet\neg A$.

**Proof:** Just turn the FACT 3.38 upside-down.

Now, some people felt unease by the presence of FACT 3.60(ii), understanding that this would mean a 'proliferation of inconsistencies' —given that any formula $A$ proved to be inconsistent would, by way of that rule, generate infinitely many 'other' inconsistencies (the negation of $A$, the negation of the negation of $A$, and so on). Be that as it may, it is still clear that in **Cie** there are no *new* inconsistencies added by way of this procedure, in a sense, given that the converse of FACT 3.60(ii) is also valid here (it is the FACT 3.38(ii)), and so, in fact, $\bullet A$ and $\bullet\neg A$ are *equivalent* formulas! In **Cile**, the logic obtained by the addition of (ce) to **Cil**, instead of **Ci** (to see that this logic is also paraconsistent, use again the matrices of $\mathbf{P}^2$ in THEOREM 3.59) we evidently obtain a new version of the above result, and FACT 3.60(ii) converts itself into $[(A\wedge\neg A) \vdash_{\mathbf{Cile}} (\neg A\wedge\neg\neg A)]$, leading, so it seems, into a 'proliferation of contradictions'. Once again, this would perhaps not be said to be the case if the formulas at the right and the left hand side are again remembered to be equivalent.

Routley & Meyer, in [100], on their attempt to define a 'dialectical logic', *DL*, meeting the standards of 'Soviet logic' and to recover the 'orthodox Marxist view of negation' have pondered the possibility of criticism coming from some dialecticians to the effect that their 'negative logic is excessively classical', and considered the constitution of a 'weaker dialectical logic', *DM*, having only (Min11), but not (ce), as an axiom. But, even in the case of their *DL*, they have met inferences such as the ones above, acknowledging the possibility of generation of an infinite number of 'distinct' contradictions from any given one, and still defended that this would be all right — one just has to remember that it is still not the case that *any* contradiction is derivable, indeed, just a very specific set of contradictions, 'forming a chain', are derivable, but this, of course, 'does not result in total system disorganization'. But how could it be that these contradictions that they obtain in *DL* are all *distinct*, as they asserted? This is a bit tricky. Let $A_0$ and $\neg A_0$ be two theorems of *DL* (it's a *dialectical* logic, after all —see the subsection **2.2**), and let $A_n$ abbreviate the formula $(A_{n-1}\wedge\neg A_{n-1})$. We have already learned in the last subsection about Routley & Meyer's argument for the validity of the schema $\neg(A\wedge\neg A)$, and from this we may conclude that each $A_n$ will be a theorem, starting from the mere fact that both $A_0$ and $\neg A_0$ hold. But if, on the one hand, it is clear that $(A_n\rightarrow A_{n-1})$ will be a theorem of *DL*, on the other hand it is equally clear that $[A_n \dashv\vdash A_{n-1}]$ holds, as above, similarly also to what had occurred in the cases of **Cie** and **Ciel** (FACT 3.60(ii) and its variations). So, again, how could it be that 'all these contradictions are distinct', as asserted by these authors? The point is that the propositional bases of both *DL* and *DM* are relevance logics, and so it may occur that $A_n$ and $A_{n-1}$ are equivalent formulas, but it is still the case that $(A_{n-1} \to A_n)$ is not a theorem of *DL*. So, after all, we see that *this* is the sense of 'distinctness' employed by those authors, determined exclusively by the validity or not of a bi-implication, and *not* by the sets of consequences of the formulas under examination.

**3.10 Consistency may be contagious!** Supposing we can really trust the consistency of some formulas in our theories, what can we say about the more complex formulas that one can build using the last ones as components: will *these* be also consistent? From FACT 3.38(i): [$\circ A \vdash \circ \neg A$] we know that already in **Ci** the consistency 'propagates' through negation, that is, the consistency of $A$ is sufficient information for us to be sure about the consistency of $\neg A$. This is essentially a consequence of (Min11): $(\neg \neg A \to A)$ and the identification of inconsistency and contradiction guaranteed by the axiom (ci). Now, what do we know about the propagation of consistency through other connectives besides negation? Not much, so far.

The idea behind the construction of the original calculi $C_n$ by da Costa (see [49] and [50]) was that of requiring each component to be consistent as a sufficient reason to count on the consistency of the more complex formula. Bluntly speaking, da Costa's $C_1$ was built by the addition to **Cil** (see the beginning of the subsection **3.8**) of the following axiomatic rules:

(ca1)  $(\circ A \wedge \circ B) \vdash \circ(A \wedge B)$;
(ca2)  $(\circ A \wedge \circ B) \vdash \circ(A \vee B)$;
(ca3)  $(\circ A \wedge \circ B) \vdash \circ(A \to B)$.

Let's call **Cila** the logic obtained by the addition of (ca1)–(ca3) to **Cil**. The difference from **Cila** and the original formulation of $C_1$ is only one: that the connective $\circ$ in $C_1$ was not taken as primitive, but $\circ A$ was instead denoted as $A^\circ$ and was taken more directly as an abbreviation of the formula $\neg(A \wedge \neg A)$ (recall the THEOREM 3.55). As for the other calculi in the hierarchy $C_n$, $1 \le n < \omega$, they were built using the simple trick of letting $\circ A$ abbreviate more and more complex formulas (as we saw at the end of the subsection **3.8**).

As an immediate consequence of the above definitions, one can easily prove in **Cila** —and in each calculus $C_n$— the following 'translating' results (compare these with the less specific THEOREM 3.11 and with the generally applicable COROLLARY 3.47):

**THEOREM 3.61** [$\Gamma \vdash_{\mathbf{CPL}} A$] $\Leftrightarrow$ [$\circ(\Pi), \Gamma \vdash_{\mathbf{Cia}} A$], where $\circ(\Pi) = \{\circ p : p$ is an atomic formula occurring as a subformula in $\Gamma \cup \{A\}\}$.[22]
**Proof:** Immediate, using (ca1)–(ca3). Note that the axiom (cl) plays no role here.

COROLLARY 3.62 The following mapping conservatively translates **CPL** inside of **Cia**:

(t2.1)  $t_2(p) = \circ p$, if $p$ is an atomic formula;
(t2.2)  $t_2(A \# B) = t_2(A) \# t_1(B)$, if # is any binary connective;

----

[22] Might the reader observe that the first formulations of this result, on da Costa's [49], Theorem 9, page 16, and on da Costa's [50], Theorem 4, page 500, the general case in which an infinite number of atomic formulas occur in $\Gamma \cup \{A\}$ is not considered.

(t2.3)  $t_2(\neg A) = \neg t_2(A)$;

(t2.4)  $t_2(\circ A) = \circ t_2(A)$.

So, working with **Cila**, all we need to rely on in order to go on making 'classical inferences' is on the consistency of the atomic constituents of our formulas. As a particular consequence of that, one can now substitute each new axiomatic rule of **Cila** by an alternative version in terms of '•'s instead of '∘'s. Thus, the axiom (ca3), for instance, can be rewritten as $[•(A \rightarrow B) \vdash (•A \vee •B)]$ (use FACT 3.34(i) and COROLLARY 3.62). And so on.

Proposing an infinite number of calculi, instead of one, as in the case of the $C_n$, only starts to make sense after we prove that we are not just repeating the same tune:

**THEOREM 3.63** Each $C_n$ deductively extends each $C_{n+1}$, for $1 \leq n < \omega$.[23]

**Proof:** We will not here give it a try by usual 'syntactical means'. For sure, this will be much easier to check if one just considers the semantics associated with these calculi, for instance in [53] (*corrected* in [74]) and in [36] (or [76]).

Evidently, all these calculi $C_n$ extend also the calculus $C_\omega$ —they even extend $C_{min}$, the stronger logic on which we based **bC**, our first **LFI** (recall the subsection **3.1** for the definition of these logics). This $C_\omega$, we argue, was indeed a very bad choice as a kind of 'limit' to the hierarchy $C_n$, $1 \leq n < \omega$. Consider, for instance, the following result:

**FACT 3.64** The only addition made by $C_n$ (in fact, by **Cia**, for the axiom (cl) has no use in this result) to the rules provable by **bC** about the interdefinability of the binary connectives (see THEOREM 3.18) is the rule (ix): $[\neg(A \wedge B) \vdash_{\mathbf{Cia}} (\neg A \vee \neg B)]$, and its variants.

**Proof:** We just show that that rule holds already in **Cia**, and point the reader again to the semantical studies of the calculi $C_n$ to check that the other formulas are still not provable in **Cila**. First of all, setting $\Gamma = \{\circ(A \wedge B), \neg(A \wedge B), A\}$ it can immediate be seen that $[\Gamma, B \vdash_{\mathbf{Cia}} \circ(A \wedge B)]$, $[\Gamma, B \vdash_{\mathbf{Cia}} (A \wedge B)]$ and $[\Gamma, B \vdash_{\mathbf{Cia}} \neg(A \wedge B)]$, so we apply FACT 3.14(ii) to obtain $[\Gamma \vdash_{\mathbf{Cia}} \neg B]$, and consequently $[\Gamma \vdash_{\mathbf{Cia}} (\neg A \vee \neg B)]$. But then, also $[\neg A \vdash_{\mathbf{Cia}} (\neg A \vee \neg B)]$, and so the proof by cases will give us $[\circ(A \wedge B), \neg(A \wedge B) \vdash_{\mathbf{Cia}} (\neg A \vee \neg B)]$. By (ca1) we then conclude that $[\circ A, \circ B, \neg(A \wedge B) \vdash_{\mathbf{Cia}} (\neg A \vee \neg B)]$, but we also have, from LEMMA 3.43(ii), that $[\vdash_{\mathbf{Cia}} (\neg A \vee \circ A)]$, and from that we obtain $[\circ B, \neg(A \wedge B) \vdash_{\mathbf{Cia}} (\neg A \vee \neg B)]$. By a similar reasoning, from $[\neg B \vdash_{\mathbf{Cia}} (\neg A \vee \neg B)]$, we finally arrive at our goal, $[\neg(A \wedge B) \vdash_{\mathbf{Cia}} (\neg A \vee \neg B)]$.

**COROLLARY 3.65** (IpE) cannot hold in the calculi $C_n$, or in any extension of them.
**Proof:** Just recall THEOREM 3.51(iv).

The FACT 3.64 also suggests some further information about the plausibility of calling either $C_\omega$ or $C_{min}$ 'limits' for the hierarchy $C_n$, $1 \leq n < \omega$. For, as we have seen in THEOREM 3.18, these new forms of De Morgan rules that we now have in each $C_n$

---

[23] A supposedly general proof of this fact, dating still from the 'syntactical period', when no semantics had yet been presented to those calculi, appears for instance in da Costa's [49], pp.17–9, and once again in Alves's [2], pp.17–9, and is credited to Ayda Arruda. There is surely some mistake, however, in their attempt to prove the independence of each axiom $[A^{(n)}, A, \neg A \vdash B]$ with respect to the axioms of $C_{n+1}$, given that this very axiom assumes non-distinguished values in all matrices $\mathbf{T}_n$ thereby presented, if one only picks 1 as the value of $A$ and picks for $B$ any value in between 1 and $n+2$.

are *not* present even in $C_{min}$. Also, by a combination of FACT 3.40(i) and COROL-LARY 3.56, we know that $((A \wedge \neg A) \rightarrow \neg\neg(A \wedge \neg A))$ is valid in **Cia**, and so in each $C_n$, while we also know, from the matrices of $\mathbf{P}^1$, as in THEOREM 3.15(ii), that this formula cannot be a theorem of neither $C_{min}$ nor $C_\omega$. Now, it is only compelling to think of a *deductive limit* for an infinite hierarchy of increasingly weaker calculi as the logic having as inferences exactly all sets of inferences common to the whole hierarchy![24] As we have shown in [39], it *is* possible to define such a logic, for each hierarchy of **dC**-systems as the one given above, by way of the useful tool of possible-translations semantics, obtaining, as a byproduct, some clear-cut and effective decision procedures, even though some other very interesting questions, such as how to finitely axiomatize this limit-calculi, or how to define a strong negation in them, if this is possible at all, were still left open (check also the paper [42]). Indeed, notice that when we go from each $C_n$ to the following $C_{n+1}$ we need in fact to add a further requirement in order to express the consistency of a formula $A$ —while in $C_n$ this was expressible by way of $A^{(n)}$, or, equivalently, by way of the set $\{A^1, A^2, ..., A^n\}$, in $C_{n+1}$ that same set must be incremented by the formula $A^{n+1}$. So, ultimately, in $C_{Lim}$, the deductive limit of the hierarchy $C_n$, the consistency of $A$ can evidently be expressed by an infinite number of formulas, and again we obtain a logic which is gently explosive, being thus an **LFI**. What we still do not know is if logics such as $C_{Lim}$ can be alternatively characterized in such a way so as to also reveal themselves as *finitely* gently explosive, like all the other logics presented up to this point, based on the axiom (bc1). If this characterization is not possible, it is hard to see how a strong negation or a bottom particle could then be *defined* in such a logic,[25] so that this would make a

---

[24] This logic $C_\omega$ has puzzled people for too long as 'part' of the original hierarchy $C_n$. The existence of a logic as a real deductive limit to this hierarchy (see [39]) shows that it was clearly just a matter of coincidence that $C_\omega$ would appear as a kind of 'syntactical limit' to the original axiomatic formulation of the $C_n$, by the deletion of all the axioms and definitions involving the connectives $\circ$ and $\bullet$ ((bc1), (ci) and (cl)), and the 'deletion' also of (Min9) (actually, this last axiom was not in the original formulation of these calculi, and could not even be proved from the other axioms of $C_\omega$ —see THEOREM 3.3). To put the matter in clear terms, $C_\omega$ *can* of course be studied in its own right, as a very weak paraconsistent (non-**LFI**) logic based on positive intuitionistic logic, but it *should not* be seen as part of the hierarchy $C_n$, $n \geq 1$, for it has no more right to occupy that position than $C_{min}$ or many other logics that could substitute it would have!

This coincidence had also some harmful effects on the philosophical appreciation of the logics produced by da Costa. As da Costa himself has put it in his original piece on these systems (cf. [49], p.21), 'roughly speaking, we could say that human reason seems to attain the peak of its power the more it approaches the danger of trivialization'. This statement has been inspiring people to naively defend stances according to which, for instance, 'the more a theory is useful to found mathematics, the more easily it results to trivialize it; and the more difficult it is to trivialize it, the less it is useful to found mathematics' (see [26], p.243). There are good and bad points about these somewhat hasty conclusions. First, as a general technical assertion about paraconsistent logics in general, da Costa's motto is certainly misleading, given the existence of maximal logics such as the three-valued *Pac* (subsection **2.4**), which is both as strong as a fragment of classical logic as it could be, and at the same time is not finitely trivializable at all. One could, then, restrict their attention to **LFI**s and repeat that motto in an environment in which it seems to make sense. In that case, of course, the second statement above would be affirming that no non-**LFI** could be useful to found mathematics —and *this* statement would be very likely to find its defensors (cf., for instance, Batens's attack [13] on Priest's [90]).

[25] Here, we really mean *defined* inside the logic, as a real *formula* of this logic. For instance, in Priest's [90] a logic containing no bottom particle is presented, but the author argues that such a propositional constant $\perp$ could be 'thought of informally as the conjunction of all formulas' (p.146), so that, for instance, a strong negation ~ would be obtainable from that in the usual way, by letting ~$A$ be defined

case in which the Gentle Principle of Explosion does not coincide with the Supplementing one, or with *ex falso* (compare this with FACT 2.19, and the comments which follow that result).

Some other interesting theses of **Cila** are the following (see Urbas's [107]):

**FACT 3.66** In **Cila** the schemas $\circ \bot$ and $\circ \div \div A$ are provable.

**Proof:** Recall first that $\div A$ was the classical negation defined inside of **bC** (in THEOREM 3.48) as $(A \rightarrow \bot)$, and $\bot$ is a bottom particle that can be defined, for instance, as $(\circ B \wedge (B \wedge \neg B))$, for some $B$. Now, by FACT 3.33 and COROLLARY 3.56, we know that both $\circ \circ B$ and $\circ (B \wedge \neg B)$ are theorems of **Cil**, and then we conclude, by the axiom (ca1), that $\circ (\circ B \wedge (B \wedge \neg B))$, and so $\circ \bot$ is a theorem of **Cila**. Recall also from FACT 3.49 that $\div A$ is equivalent, in **Ci**, to $\sim A$, and this last strong negation was defined as $(\neg A \wedge \circ A)$, and so we have in particular that $[(A \rightarrow \bot) \vdash_{\textbf{Cila}} \circ A]$. So, from this last inference and from the fact that $\circ \bot$ is a theorem of **Cila**, as proved above, we use (ca3) and conclude that $[(A \rightarrow \bot) \vdash_{\textbf{Cila}} \circ (A \rightarrow \bot)]$. As particular cases of the last two inferences, substituting $A$ for $((A \rightarrow \bot) \rightarrow \bot)$ in the first case, and for $(A \rightarrow \bot)$ in the second case, we obtain, respectively, $[\div \div \div A \vdash_{\textbf{Cila}} \circ \div \div A]$ and $[\div \div A \vdash_{\textbf{Cila}} \circ \div \div A]$, and the version of proof by cases obtained for $\div$ leaves us at last with $[\vdash_{\textbf{Cila}} \circ \div \div A]$.

This last result gives us yet another reason for the failure of (IpE) in the $C_n$ and in their extensions (COROLLARY 3.65) —now, it is THEOREM 3.51(i) that applies. As we shall see in subsection **3.12**, the failure of (IpE) makes it impossible for us to find a Lindenbaum-Tarski-like algebraization for these logics. In the case of the $C_n$ the situation is actually worse: as Mortensen ([82]) has shown, no non-trivial congruence is definable for these logics, making these logics non-algebraizable even in a much more general sense, the one of Blok-Pigozzi (cf. THEOREM 3.83). There are, nevertheless, several extensions of the $C_n$ in which non-trivial congruences *can* be defined, being thus much more receptive to algebraic treatments. We will be seeing many examples of these below.

Let us now investigate another way of propagating consistency, by liberalizing a little bit the conditions required by **Cila** (that is, $C_1$). Da Costa, Béziau & Bueno proposed, in [57], to substitute the above axioms, (ca1)–(ca3), by the following:

    (co1)    $(\circ A \vee \circ B) \vdash \circ (A \wedge B)$;
    (co2)    $(\circ A \vee \circ B) \vdash \circ (A \vee B)$;
    (co3)    $(\circ A \vee \circ B) \vdash \circ (A \rightarrow B)$.

We will call **Cilo** the logic obtained by the addition of (co1)–(co3) to **Cil**. It is very easy to see, using the positive axioms, that this logic, christened $C_1^+$ in [57], is a deductive extension of $C_1$. Requiring less assumptions in order to obtain consistency of a complex formula in terms of the consistency of its components, **Cilo** (or even **Cio**, already, without recourse to the axiom (cl)) gives us some interesting results such as:

**THEOREM 3.67** $[\Gamma \vdash_{\textbf{Cio}} \circ A]$ whenever $[\Gamma \vdash_{\textbf{Cio}} \circ B]$, for some subformula $B$ of $A$.
**Proof:** Immediate, using (co1)–(co3).

---

as $(A \rightarrow \bot)$. But such an 'informal' bottom particle is simply no formula of our language! (And if it were, then we would have been done: all paraconsistent logics would turn out to be **LFI**s, in one way or another). This idea of Priest has already been criticized in Batens's [13].

FACT 3.68 **Cio** makes some new additions to FACT 3.64 and to the rules displayed in THEOREM 3.18 about the interdefinability of the binary connectives, namely the following provable schematic rules:

(vi) $\neg(A \wedge \neg B) \vdash_{\mathbf{Cio}} (A \rightarrow B)$;

(vii) $\neg(A \rightarrow B) \vdash_{\mathbf{Cio}} (A \wedge \neg B)$;

(xi) $\neg(\neg A \vee \neg B) \vdash_{\mathbf{Cio}} (A \wedge B)$.

**Proof:** Go to [57] or [76], or else the section on semantics for **Cibo** (and **Ciboe**) in [42], to check this. The actual syntactical proofs are, in any case, structurally similar to the one presented in FACT 3.64.

Of course, given THEOREM 3.51(v), we know that FACT 3.68(vii) gives us yet another reason for the failure of (IpE) from this point on. But that we already knew, from the case of **Cila**, in FACT 3.64. What is new in this case is only that **Cilo** *can*, differently from **Cila**, define non-trivial congruences, making it possible to algebraize it *à la* Blok-Pigozzi (see FACT 3.81).

Once again, the axioms (co1)–(co3) have equivalent versions in terms of •, instead of ∘, and it is an easy exercise to try to find them. The reader should remember that both **Cila** and **Cilo** have not only associated decreasing hierarchies, and can evidently define different calculi as their deductive limits, but they also can be structurally varied in terms of their inner definition of consistency, if we only change the axiom (cl) for axiom (cd), or (cb), or (cg), as we did in the subsection **3.8**, or if we add to them the axiom for expansion of negations, (ce), as in the subsection **3.9** (defining the logics **Cido**, **Cibo**, **Cito**, **Ciloe**, **Cidoe**, **Ciboe** and **Citoe**). In all these cases, we can show that the resulting logics are extended by the three-valued paraconsistent logic $\mathbf{P}^2$, as in the THEOREM 3.59.

There are, actually, an unlimited number of ways of propagating consistency.[26] Before proceeding to a general investigation of the 'extreme cases', in the next subsection, let us just briefly survey some propagation axioms which have already shown up in the literature so far, and some of their consequences. Consider, for instance, the following axioms, converse to (co1)–(co3):

(cr1) $\circ(A \wedge B) \vdash (\circ A \vee \circ B)$;

(cr2) $\circ(A \vee B) \vdash (\circ A \vee \circ B)$;

(cr3) $\circ(A \rightarrow B) \vdash (\circ A \vee \circ B)$.

Adding these to **Cibo** and to **Cio** (that is, **Cibo** minus the axiom (cb)) we build, respectively, the logics **Cibor** and **Cior** (and so on, *mutatis mutandis*, for **Cilo** and **Cido**). These give us yet more perspectives on the propagation (and back-propagation) of consistency, and some of the possible meanings which can be assigned to it). Suppose, on the other hand, that we, more simply, consider the following axioms in order to automatically guarantee the consistency of some complex propositions:

(cv1) $\vdash \circ(A \wedge B)$;

(cv2) $\vdash \circ(A \vee B)$;

---

[26] Working with monotonic logics, the addition of consistency-propagation to some logic always means a gain in deductive strength. However, in a non-monotonic environment in which consistency is pressuposed by default, given that propagating consistency in one direction can mean propagating inconsistency the other way, the addition of such axioms for propagation can either be innocuous or in some cases even have a weakening effect on the resulting system (see [79], and [15]).

(cv3)    $\vdash \circ(A \to B)$;

(cw)    $\vdash \circ(\neg A)$.

Let's add **v** to the name of a logic that contains the axioms (cv1)–(cv3), and add **w** to the name of a logic containing (cw). Let's also recall the axiom (ce): $[A \vdash \neg\neg\neg A]$, from which we obtained the FACT 3.60 (backward propagation of consistency through negation). There are now several new possible combinations to be considered. Evidently, any logic having (cv1)–(cv3) proves not only (co1)–(co3) but also (ca1)–(ca3). So we have, for instance, the logic **Cibv**, and have at least two immediate ways of enriching it with respect to the behavior of negation, obtaining the logics **Cibve** and **Cibvw**. As it happens, **Cibvw** axiomatizes the three-valued maximal paraconsistent logic $\mathbf{P}^1$ (matrices in THEOREM 3.15(ii)), that has the peculiarity of admitting inconsistency only at the atomic level, and **Cibve** axiomatizes $\mathbf{P}^2$ (matrices in THEOREM 3.59), a logic that admits inconsistency only at the level of atomic propositions, or of propositions of the form $(\neg^n p)$, where $p$ is atomic and $\neg^n$ denotes $n$ applications of negation. If, on the other hand, one considers the logic **Ciborw**, once more a three-valued paraconsistent logic pops up, namely the one given by the following matrices:

| $\wedge$ | **1** | **½** | **0** |
|---|---|---|---|
| **1** | 1 | 1 | 0 |
| **½** | 1 | ½ | 0 |
| **0** | 0 | 0 | 0 |

| $\vee$ | **1** | **½** | **0** |
|---|---|---|---|
| **1** | 1 | 1 | 1 |
| **½** | 1 | ½ | 1 |
| **0** | 1 | 1 | 0 |

| $\to$ | **1** | **½** | **0** |
|---|---|---|---|
| **1** | 1 | 1 | 0 |
| **½** | 1 | ½ | 0 |
| **0** | 1 | 1 | 1 |

| | $\neg$ | $\circ$ |
|---|---|---|
| **1** | 0 | 1 |
| **½** | 1 | 0 |
| **0** | 1 | 1 |

where 1 and ½ are both distinguished. We shall here call this logic $\mathbf{P}^3$. All these three logics, $\mathbf{P}^1$, $\mathbf{P}^2$ and $\mathbf{P}^3$, are in fact **dC**-systems, and can ultimately dispense with the axiom (cb), proving it from the other axioms. If, on the other hand, we consider the logic **Ciore**, the result is a maximal three-valued logic again, **LFI2** (investigated in [44]) whose matrices differ from those of $\mathbf{P}^3$ only in the matrix of negation, assigning ½ instead of 1 as the value of $\neg(½)$.

All the above logics have some kind of 'non-structural' propagation of consistency, that is, a propagation that does not really depend on the particular connective in focus. Alternatively, one can propose other forms of propagation which do depend on the connectives being considered. Now, some reasonable symmetry conditions on inconsistency and its behavior with respect to the different connectives could suggest to us, for instance, the consideration of the following forms:

(cj1)    $\bullet(A \wedge B) \dashv\vdash ((\bullet A \wedge B) \vee (\bullet B \wedge A))$;

(cj2)    $\bullet(A \vee B) \dashv\vdash ((\bullet A \wedge \neg B) \vee (\bullet B \wedge \neg A))$;

(cj3)    $\bullet(A \to B) \dashv\vdash (A \wedge \bullet B)$.

The logic **Cij**, built from the addition of (cj1)–(cj3) to **Ci**, can now be enriched with (ce) in order to give us **Cije**, an axiomatization for the above many times mentioned maximal three-valued paraconsistent logic **LFI1** (see its matrices in FACT 3.29, and consult [44] again). It is interesting enough to note that neither **LFI1** nor **LFI2** are **dC**-systems, that is, they *cannot* define the consistency operator by way of the other connectives. Let us now just summarize, give some references and mention some properties of the five above mentioned maximal paraconsistent three-valued logics, before we proceed to show, in the next subsection, that these are just the top of the iceberg:

**THEOREM 3.69** The matrices of $\mathbf{P}^1$ (in THEOREM 3.15(ii)) are axiomatized by **Civw**; the matrices of $\mathbf{P}^2$ (in THEOREM 3.59) by **Cive**; the matrices of $\mathbf{P}^3$ (above) by **Ciorw**; the matrices of **LFI2** (above) by **Ciore**; the matrices of **LFI1** (in FACT 3.29) by **Cije**.

**Proof:** As a general reference for all these logics and all the other ones in the next section, consult [78]. As specific references for some of them, go to [103] and [84] for $\mathbf{P}^1$ (or [10], where it appeared under the name $PI^v$); notice that $\mathbf{P}^2$ has also appeared in [10], under the name $PI^m$, but was then redefined in [84] (where it actually was wrongly supposed to be characterizable using just one distinguished value, invalidating the soundness proof therein presented —note 13 of that paper shows that the author had even been informed about that) and later rediscovered in [76]; go to [44] for **LFI2** and **LFI1** (but the reader should bear in mind that this last logic is in fact equivalent to the logic called $\mathbf{J}_3$ in [60], and to the propositional fragment of the logic called **CLuNs** in [15], which is in fact identical to the logic $\Phi_v$ presented in [101], a logic that has been reappearing quite often in the literature).

**FACT 3.70** Of the rules displayed in THEOREM 3.18 on the interdefinability of the binary connectives, these are the instances validated by each of the five three-valued logics above:

(i)  in $\mathbf{P}^1$, $\mathbf{P}^2$, $\mathbf{P}^3$ and **LFI2**: parts (i), (iv), (vi), (vii), (ix), (xi);

(ii)  in **LFI1**: parts (i), (iii), (iv), (v), (vii), (viii), (ix), (x), (xi), (xii).

Also, formulas such as $\neg(A \wedge \neg A)$ and $\neg(\neg A \wedge A)$ may be easily seen to hold in **LFI1** and **LFI2**, and rules such as $(A \wedge \neg A) \vdash \neg\neg(A \wedge \neg A)$ hold in $\mathbf{P}^2$, **LFI1** and **LFI2**.

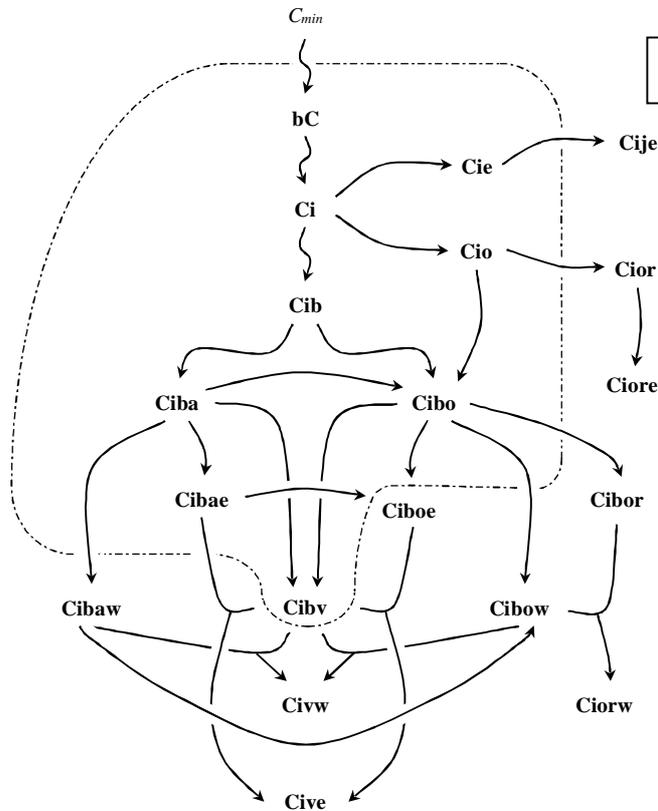

**Figure 3.1**

**Proof:** Just use the corresponding matrices to check this. Notice from FACT 3.68 that part (iv) was the only addition made by the first four logics above to the rules already validated by **Cilo**.

This last result, of course, supplies us with still some further justifications for the failure of (IpE) in all these logics: parts (iii) and (vi) of THEOREM 3.51 applies to **LFI1**, parts (ii) and (vii) of 3.51 apply to both **LFI1** and **LFI2**, parts (iv) and (v) of 3.51 apply to all of the five logics; and finally we will see in FACT 3.76 that part (i) of 3.51 also applies to all of them.

We can, at this point, try our hand at sketching a very thin slice of the great number of **C**-systems introduced so far. Doing that, something like **Figure 3.1** might eventually be obtained. In that figure, an arrow leading from a logic **L**1 into a logic **L**2 says that **L**2 deductively extends **L**1. The logic $C_{min}$, at the upper end, is the only one that does not constitute a **C**-system; the logics at the lower ends are the three-valued ones appearing in THEOREM 3.69. The logics inside the dotted lines are some of those which we can prove to be *not* many-valued, by adapting the results in [76], pp.213–216 or, better, by checking [42]. The other logics not contemplated by these results, namely **Cior**, **Ciboe**, **Cibor**, **Cibaw**, and **Cibow**, are also conjectured to be not many-valued, but we must at this moment leave the proof of this fact in the hands of our clever readers. Do remember to have a look, however, at the elegant possible-translations semantics offered to **Ciboe** in [42] (also originating from [76], section 5.3).

### 3.11 Taking it literally: the Brazilian plan completed.

The sagacious reader will have observed that all we have been doing so far, in this section on axiomatization of **C**-systems, was to basically try to explore at a very general level some of the possibilities for the formalization and understanding of the relationship between the concepts of consistency, inconsistency and contradictoriness. In particular, this research line makes it possible for us to reconsider and pursue, in an abstract perspective, a specific interpretation of da Costa's method and requisites on the construction of his first paraconsistent calculi (see **dC[i]**–**dC[iv]**, in the subsection **3.8**). Indeed, starting from the intuition that consistency should be expressible inside some classes of paraconsistent logics, and assuming furthermore that the consistency of a given formula would be enough to guarantee its explosive character (that is, assuming a Gentle Principle of Explosion, as formulated in the subsection **2.4**), we have arrived at the definition of an **LFI**, a Logic of Formal Inconsistency (see (D20), in the same subsection). To realize that (in a finitary way), we have above proposed the axiom (bc1): [∘$A$, $A$, ¬$A$ ⊢ $B$] for particular classes of **C**-systems based on classical logic. Even more than that, as we have remarked before, while **dC[ii]** simply establishes the non-explosive character of the paraconsistent negation, a general formulation of **dC[i]** is realized in a subclass of the **C**-systems, the ones in which the connectives '∘' and '•' happen to be definable from the remaining connectives, and to the members of this class we gave the name of **dC**-systems. Now, putting **dC[iii]**, the problem of providing higher-order versions of these logics, aside for a moment, we still need to provide an answer to **dC[iv]**, the requirement that 'most schemas and rules of classical logic' should hold in our logics. And that's the point we will ruminate in the present subsection.

Our proposed interpretation for **dC[iv]** will in fact be a very simple one, involving the following notion of 'maximality'. A logic **L**2 is said to be *maximal relative to* a logic **L**1 if: (i) both are written in the same language (so that they can be deductively compared); (ii) all theorems of **L**2 are provable by **L**1; (iii) given a theorem $D$ of **L**1 which is not a theorem of **L**2, if $D$ is added to **L**2 as a new schematic axiom, then all theorems of **L**1 turn to be provable. The idea, of course, is that any deductive extension of **L**2 contained in **L**1 and obtained by adding a new axiom to **L**2 would turn out to be identical to **L**1. We will call *maximal*, to simplify, any logic **L**2 which is maximal relative to some logic **L**1, previously introduced. Examples of maximal logics abound in the literature. It is widely known, for instance, that each Łukasiewicz's logic $Ł_m$ is maximal relative to **CPL**, the classical propositional logic, if and only if $(m-1)$ is a prime number. We also know that **CPL** is maximal relative to a 'trivial logic', in which all formulas are provable, but on the other hand it is also well-known that intuitionistic logic is *not* a maximal fragment of **CPL**, as the existence of an infinite number of *intermediate logics* promptly attests. As to the **C**-systems which have been introduced this far, only the five three-valued ones that were collected in the THEOREM 3.69 are maximal relative to **CPL**, or else relative to **eCPL**, the extended version of **CPL** introduced in the subsection **3.7** (so that, in particular, the calculus $C_1$, that we have presented as **Cila**, the strongest calculus introduced by da Costa on his first hierarchy of paraconsistent calculi, or the even stronger calculus $C_1^+$, that we presented as **Cilo**, proposed by da Costa and his collaborators much later, readily fail to be maximal, and to respect **dC[iv]**).

Let us explore then the idea that underlies the five three-valued maximal **C**-systems above. Suppose we are faced with this problem of finding models to contradictory, and yet non-trivial, theories. We might then intuitively start looking for non-trivial interpretations under which both some formula $A$ and its negation $\neg A$ would be simultaneously validated. A very simple such interpretation would be found in the domain of the many-valued. Suppose we try to depart from classical logic as little as possible, so that the interpretation of our connectives will still be classical if we remain inside the classical domain, and suppose we just introduce a third value (½), besides true (1) and false (0), so that this third value will also be seen as a modality of *trueness*, that is, ½ will also be a distinguished value, together with 1, while 0 will be the only non-distinguished value. There are then two possible negations which are such that there is a model for both $A$ and $\neg A$ being

| | $\neg$ |
|---|---|
| **1** | 0 |
| ½ | ½ or 1 |
| **0** | 1 |

true, for some formula $A$ (see the table to the right). One of these negations, the one that takes ½ into ½, is exactly the negation of **LFI1** and of **LFI2**, the other negation, taking ½ into 1, is exactly the negation of $\mathbf{P}^i$, $\mathbf{P}^2$ and $\mathbf{P}^3$. What about the other connectives? Let us again try to keep them as classical as possible (we want to keep on investigating **C**-systems *based on classical logic*), even at the level of the third value, that is, let us add to the requirement of coincidence of classical outputs for classical inputs the further higher-level 'classical' requirements to the effect that:

    (C∧)    $v(A \wedge B) \in \{½, 1\} \iff v(A) \in \{½, 1\}$ and $v(B) \in \{½, 1\}$;

    (C∨)    $v(A \vee B) \in \{½, 1\} \iff v(A) \in \{½, 1\}$ or $v(B) \in \{½, 1\}$;

    (C→)    $v(A \rightarrow B) \in \{½, 1\} \iff v(A) \notin \{½, 1\}$ or $v(B) \in \{½, 1\}$.

This leaves us then with the following options:

| ∧ | **1** | **½** | **0** |
|---|---|---|---|
| **1** | **1** | ½ or 1 | **0** |
| **½** | ½ or 1 | ½ or 1 | 0 |
| **0** | **0** | 0 | **0** |

| ∨ | **1** | **½** | **0** |
|---|---|---|---|
| **1** | **1** | ½ or 1 | **1** |
| **½** | ½ or 1 | ½ or 1 | ½ or 1 |
| **0** | **1** | ½ or 1 | **0** |

| → | **1** | **½** | **0** |
|---|---|---|---|
| **1** | **1** | ½ or 1 | **0** |
| **½** | ½ or 1 | ½ or 1 | 0 |
| **0** | **1** | ½ or 1 | **1** |

Thus, we have, theoretically, $2^3$ options of 'conjunctions', $2^5$ options of 'disjunctions', $2^4$ options of 'implications', and, as we saw above, $2^1$ options of 'negations', making a total of $2^{13} (= 8{,}192$, or 8K) possible 'logics' to play with. To remove the scare quotes of the previous passage we just have to show that these logics make some sense, and are worthy of being explored. To such an end, and to complete the definition of our 8K logics as **LFI**s, we will just also add to these logics the connectives for consistency and for inconsistency, implicitly assuming that the consistent models are the ones given by classical valuations, and only those (see matri-

| | ∘ | • |
|---|---|---|
| **1** | 1 | 0 |
| **½** | 0 | 1 |
| **0** | 1 | 0 |

ces to the right). Evidently, all these 8K logics will be fragments of **eCPL**, the Extended Classical Propositional Logic (recall the subsection **3.7**). It is also clear that the logic *Pac* (subsection **2.4**) is not one of these, for it cannot define the connectives ∘ and •, though its conservative extension **LFI1** can (and it is one of the 8K).

Evidently, the five three-valued logics we discussed earlier are but special cases of the above outlined 8K logics, and we already know (from THEOREM 3.69 and **Figure 3.1**) that those five are axiomatizable by way of the addition of suitable axioms to the axiomatization of **Ci**, one axiom for each connective. In fact, as shown in [78], this idea can be extended to all the 8K logics above:

**THEOREM 3.71** All the 8K three-valued sets of matrices above are axiomatizable as extensions of **Ci**.

**Proof:** In each case, one just has to add, for the negation, either the axiom $(A \rightarrow \neg\neg A)$ or the axiom ∘¬$A$, depending respectively if the negation of ½ goes to ½ or to 1. And, for each other binary connective, $\#(\in \{\wedge, \vee, \rightarrow\})$, one just has to add either ∘$(A\#B)$ or else (•$(A\#B) \leftrightarrow \sigma(A, B)$), where $\sigma(A, B)$ is a schema depending only on $A$ and on $B$ —these last axioms will evidently depend on the specific matrices of each $\#$, and act in order to describe how inconsistency (or consistency) propagates back and forth for each binary connective. Full details on how to define these axioms may be found in [78].

Moreover, one can also prove that:

**THEOREM 3.72** All the 8K three-valued logics above are distinct from each other, and they are all maximal relative to **eCPL**.

**Proof:** Again, we refer to [78] for the general proofs. The basic idea behind the proof of *distinctness* is the following: choosing any two of these 8K logics (without repetition), there will be some connective about which they differ, one of them giv-ing 1 as an output for the same input(s) that the other one gives ½. But then the negations of such matrices will not be equivalent, and all we must do then is write down a formula which describes that situation in such a way that this formula will be a theorem of one of these logics, but a non-theorem of the other (again, see [78] on how to do it). For the *maximality* proofs, the reader might mind to be informed that for at least five of those logics (the ones referred to in THEOREM 3.69) the specific proofs were already presented elsewhere. It is also interesting to remark that the connective ∘ (or •) plays a fundamental role in the general maximality proof exhibited in [78].

Now, how do these 8K logics compare with the other **C**-systems that have been studied this far? It is this simple: *every* logic investigated so far either coincides or is extended by some of the above 8K three-valued logics. So that now we have a very interesting class of (extended) solutions to the problem posed by da Costa's requirement **dC[iv]**! Furthermore, it is straightforward to check that all the above matrices do not only extend **Ci**, but also extend **Cia**, so that the original hunch by da Costa for the propagation axioms is a kind of minimal condition obeyed by every one of our 8K maximal three-valued logics. The only limitative point of the original proposal, under this approach, really rests in **dC[i]**, which is of course *not* verified by all those matrices, and in fact imposes a very restricted interpretation for the notion of consistency, limiting our sample space to only a *very* selective class of **dC**-systems, which is, however, larger than the reader might initially imagine (recall, in any case, that logics such as **LFI1** and **LFI2** are *not* **dC**-systems). Indeed:

FACT 3.73 All the 8,192 logics above are **C**-systems extending **Cia**. Of these, 7,680 are in fact **dC**-systems, being able to define ∘ and • in terms of the remaining connectives (and being maximal, thus, relative to **CPL**, and not only to **eCPL**). Of these, 4,096 are able to define ∘$A$ as ¬$(A \wedge \neg A)$, and so all of these do extend $C_1$ (that is, **Cila**). Of the 7,680 logics which are **dC**-systems, 1,680 extend **Cio**, the stronger alternative to **Cia**, and 980 of these are able to define ∘$A$ as ¬$(A \wedge \neg A)$, so that these 980 logics do extend $C_1^+$ (that is, **Cilo**).

**Proof:** This is just a combinatorial exercise on the above matrices, and we shall leave it for the reader to check.

It might well be that not all of the above 8K three-valued maximal logics will be interesting as logics. Some of them, for instance, do not have symmetric matrices for the conjunction or for the disjunction (but notice that some such logics have had their use in results such as 3.21, 3.26, 3.35 and 3.57, or in 3.58), though any conjunction / disjunction is evidently equivalent to any other conjunction / disjunction (the negations of these conjunctions / disjunctions are what may differ). The fact that all the 8K three-valued logics do extend **Cia** (FACT 3.73) informs us, as a corollary to FACT 3.64, that:

FACT 3.74 (IpE) cannot hold in any of the 8K logics above.
**Proof:** Again, just recall THEOREM 3.51(iv).

Now, if, in the next subsection, this failure of the replacement theorem will be seen to constitute a negative answer for the possibility of obtaining a Lindenbaum-Tarski-style algebraization for these logics (as already occurred for the calculi $C_n$ and all of their extensions —see COROLLARY 3.65), the following result will help us to show in the following a positive answer for the possibility of obtaining a Blok-Pigozzi-like algebraization to each one of them (as already hinted for some extensions of $C_n$, such as $C_1^+$, see FACT 3.81 and FACT 3.82):

**FACT 3.75** The following matrices of *classical negation* and *congruences* can be defined in each one of the above 8K logics:

| | ~ |
|---|---|
| **1** | 0 |
| **½** | 0 |
| **0** | 1 |

| ≡ | **1** | **½** | **0** |
|---|---|---|---|
| **1** | 1 | 0 | 0 |
| **½** | 0 | ½ or 1 | 0 |
| **0** | 0 | 0 | 1 |

**Proof:** To define the classical negation ~ one just has first to define ⊥ either as $(B \land (\neg B \land \circ B))$ or as $(\circ B \land \neg \circ B)$, for some formula $B$, and then define ~$A$ either as $(\neg A \land \circ A)$ or as $(A \rightarrow \bot)$. To define one of the above congruences one just has to set $(A \equiv B)$ as $((A \leftrightarrow B) \land (\circ A \leftrightarrow \circ B))$. If one wants to make sure that $v(A \equiv B) = 1$ when both $v(A) = $½ and $v(B) = $½, this is also possible: just set some $(A \stackrel{\cdot}{\equiv} B)$ as ~~$(A \equiv B)$.

In fact, it is not difficult to see that the above classical negation is indeed the *one and only* matrix of a strong negation that can be defined inside of these 8K three-valued logics (the paper [42] will also come back to this question). The reader will notice that this negation is indeed, in a sense, a 'highly' classical one. Indeed, it comes as a corollary, for instance, that:

**FACT 3.76** The schema $\circ$~$A$ is provable in all of the above 8K three-valued logics.

This last result is more than what one needs to confirm, by way of THEOREM 3.51(i), the fact that (IpE) cannot hold in these logics (as in FACT 3.74).

A noteworthy expressibility result that can be proved for these 8K three-valued logics is the following:

**FACT 3.77** (i) The matrices of $\mathbf{P}^I$ can be defined inside of any of the 8K three-valued logics above. (ii) All the matrices of all the 8K logics above can be defined inside of **LFI1**.

**Proof:** To check part (i), let $\land$, $\lor$, $\rightarrow$, $\neg$, $\circ$ and $\bullet$ be the connectives of any of the 8K logics above, and let ~ be the classical negation, defined inside this logic as in FACT 3.75. Then, the $\mathbf{P}^I$'s negation of a formula $A$ can be defined as ~~$\neg A$, the $\mathbf{P}^I$'s conjunction of some given formulas $A$ and $B$, in this order, can be defined either as ~~$(A \land B)$ or as $(\sim\sim A \land \sim\sim B)$, and the same we did for conjunction applies to both disjunction and implication, *mutatis mutandis*. The matrices for the connectives $\circ$ and $\bullet$ already coincide in all of these logics. Part (ii) is a particular consequence of the expressibility result that we have proven in [44], Theorem 3.6. In that result we showed, in fact, that the matrices definable in **LFI1** are all those, and exactly those, $n$-ary matrices that have classical (1 or 0) outputs for classical inputs (and that can have any output value if non-classical inputs are considered). Of course, all the above matrices, on these 8K three-valued logics, are, by definition, just 1-ary and 2-ary examples of such **LFI1**-definable *hyper-classical* matrices, as we have called them.

**COROLLARY 3.78** (i) The logic $\mathbf{P}^I$ can be conservatively translated inside any of the 8K three-valued logics above. (ii) Any of the 8K logics above can be conservatively translated inside of **LFI1**.

Are there other interpretations, besides *maximality*, of da Costa's requisite **dC[iv]** leading to yet some other classes of solutions to the problem of finding **C**-systems containing 'most rules and schemas of classical logic'? Are there non-many-valued (monotonic) solutions to that problem, or perhaps some other $n$-valued ones, for $n > 3$? And, this is an important first step and probably an easier problem to solve, are there other interesting **C**-systems based on classical logic which are *not* extended by any of the above 8K three-valued logics? We must leave these questions open at this stage. It is interesting to notice, at any rate, that this problem has already been addressed here and there, in the literature. Besides [78], from which we drew the results in the subsection **3.12**, one could also recall, for instance, the adaptive programme for the confection of paraconsistent logics aiming to represent (non-monotonically) the dynamics of scientific reasoning and of argumentation (see [15]). Roughly speaking, the basic idea behind adaptive logics is that of working in between two boundary logics, classical logic often being one of them and a paraconsistent logic being the other one, so that consistency is pressuposed by default and we try to keep on reasoning (i.e. making inferences) inside of classical logic up to the point in which an 'abnormality' (an inconsistency?) pops out, a situation in which we had better descend to the level of the complementary paraconsistent logic, and go on reasoning over there. Indeed, the ancestral motivations of this programme (see [10]) seem to have been, as it is reasonable to conceive, yet another attempt to originate logics maintaining as much of classical logic as possible, so that paraconsistency will only be needed at limit cases.

As the reader will see in [42], when we go on to provide possible-translations semantics (as in [36], [39] or [76]) to some of the above non-three-valued logics, as for instance **Ci** or **Ciboe**, the intuition behind the construction of the previous three-valued paraconsistent logics can be pushed much farther, since it can be shown that some infinite-valued logics can also be *split* in terms of suitable combinations of clusters of three-valued logics.

**3.12 Algebraic stuff.** You may think, perhaps, that logic has 'too many formulas'. There is nothing unreasonable in supposing, however, that some of these formulas can in fact be identified, and indistinctly used in all contexts. If we consider classical logic, for instance, we will promptly see that there is no reason to distinguish between any two given theorems (or top particles) with respect to their relation to the other formulas of the classical language —even though they may very well still be understood as 'expressing' quite different facts, somehow conveying different bits of information. In the classical case, also, and for the very same reason, one does not really need to distinguish between any two given bottom particles (or two different pairs of contradictory formulas); even more than that, any two formulas $A$ and $B$ which turn out to be provably equivalent (that is, such that $[A \dashv\vdash B]$, or, what in many logics amounts to be just the same, to put it in terms of a bi-implication, such that $[\vdash (A \leftrightarrow B)]$) are, in a certain sense, indistinguishable, and can be indistinctly employed in the same contexts, to attend similar purposes. The action of putting the glasses through which some 'contingent' properties of formulas are hidden and only those features related to their general behavior in relation to other formulas are exposed is the task of *algebraization*. Whenever a given logic turns out to be algebraizable, so that the logical problems can be faithfully and conservatively translated into some given well behaved algebra, then it will be possible for us to use the powerful (universal) algebraic tools to tackle those problems, so that, in the next and final step, we will be able to translate the results back into logic.

Being the above remarks all too informal, we had better strive to put them in more precise terms. If the whole activity of mathematics and logic involves 'forgetting' some things (and calling attention to others), and identifying what could otherwise, at the first look, have seemed just different, their tools by excellence, in such respect, are the key notions of *equivalence*, *congruence*, *isomorphism*, and so on. Once the very definition of a logic, as we have proposed it at the start of section **2**, can immediately be seen as some sort of algebra having the set of formulas as its domain (indeed, in the structural and propositional case, it is exactly the free algebra generated by the primitive connectives of the languages —here understood as operators— over the set of atomic propositions), the quest for dividing these formulas into disjunct packages of equivalent and indiscernible ones can easily be accomplished if one is able to define a *congruence* relation over these formulas, that is, an *equivalence* (reflexive, symmetric and transitive) relation such that any two equivalent formulas (with respect to this relation) can be just justifiably and indistinctly used in all and the same contexts. So, if some given formula $A$ appears as a component of some other formula $G$, then any formula $B$ which is congruent to $A$ should be able to do the same job, with no loss or increase in expressibility or generality. What we implicitly mean with this is that the new *quotient* algebra obtained by dividing the original algebra of formulas by this congruence relation should preserve the original 'operations', being thus homomorphic to the original algebra. So, in dividing the formulas this way into classes of congruent ones, one can go on to work and dialogue with just the (arbitrary) representants of these classes, once they are supposed to behave exactly the same as any of their congruent colleagues with respect to any operation of (any isomorph) of the quotient algebra. Any two congruent formulas are 'the same up to a congruence', and can play exactly the same roles in some specific dramas.

It comes perhaps as no surprise the confirmation that the most easy and standard way of algebraizing a given logic is obtained by way of the relation of provable equivalence induced by its underlying consequence relation. Indeed, the so-called *Lindenbaum-Tarski algebraization* sets two formulas $A$ and $B$ as congruent if $A \dashv\vdash B$ —let's denote this fact by writing $A \approx B$. Such congruence relation $\approx$ is evidently an equivalence relation, and to confirm that any two so congruent formulas 'work the same in all contexts', one has to check if they can be intersubstituted everywhere, that is, one has to prove a *replacement* theorem, to the effect that the *intersubstitutivity of provable equivalents*, (IpE), holds (recall its definition in the subsection **2.3**, and check [111]). Many logics have Lindenbaum-Tarski-like algebraizations, as it is the case for classical logic, intuitionistic logic, several normal modal logics, several many-valued logics, and so on. But not all algebraizable logics are algebraizable in the sense of Lindenbaum-Tarski, not being able for instance to prove replacement with respect to provable equivalence (or provable bi-implication). In the case of many non-normal modal logics, for example, what one needs is *strict* (that is, *necessary*) provable equivalence (that is, strict bi-implication). In the case of the paraconsistent logics studied here, frequent negative results on what concerns the validity of (IpE) —and so, on what concerns the possibility of obtaining an algebraization *à la* Lindenbaum-Tarski— have been met: In fact, *all* of the above **C**-systems have been shown at some point to lack (IpE) (recall the results 3.22, 3.35, 3.58, 3.65, and 3.74). Yet, the possibility of obtaining some positive results within some extensions of those **C**-systems was not ruled out (recall 3.53, but confront it with 3.51).

An immediate result about algebraizations that may come quite handy is the following one (cf. [25], Corollary 4.9):

FACT 3.79 Every deductive extension of an algebraizable logic is algebraizable.

A case study which was particularly well investigated is that of the logic **Cila** (the logic $C_1$ of da Costa's [49] —check the subsection **3.10**). Even though at least as early as in da Costa & Guillaume's [54] it had already been noticed that (IpE) does not hold for **Cila** (COROLLARY 3.65), so that no Lindenbaum-Tarski-like algebraization for this logic (or for any other of the weaker calculi $C_n$) can be available, several attempts have been made to find other kinds of algebraizations for this logic (check, for instance, da Costa's [48]). The intuitive idea underlying the search of other algebraizations, generalizing the idea of Lindenbaum-Tarski, has been quite often that of finding 'any' congruence on the set of formulas that could be used to produce a quotient algebra from the algebra of formulas of the logic. Furthermore, if such a congruence is no more necessarily supposed to be induced directly by way of the consequence relation associated to the logic, nor should this congruence be necessarily supposed to be expressible by way of a formula written in the very language of the logic (it may happen to be definable only metalinguistically —for instance, if you do need a metalinguistical 'and' to characterize it, but there is no adequate conjunction available to express it in the language of the logic), it is still reasonable to suppose as well that this congruence should put no distinguished and non-distinguished formulas inside the same class of equivalence (so, for instance, no class will simultaneously contain a theorem and a non-theorem), so that we will have no trouble in attributing a distinguished or a non-distinguished status to some class of equivalence (cf. [84]). The final blow to the search for congruences algebraizing the logic **Cila** was delivered by Mortensen's [82], where this author proved that:

THEOREM 3.80 No non-trivial quotient algebra is definable for **Cila**, or for any logic weaker than **Cila**.

It is never too late to remember that a *trivial quotient algebra* is an algebra defined by a congruence relation $\approx$ such that $A \approx B$ if, and only if, $A$ and $B$ are the same formula (so, all equivalence classes are singletons). Now, some authors have argued that the exclusive existence of trivial quotient relations for a given logic is a major 'defect' (cf. [84], section 3), while others do not think so (cf. [20]) —and this is the reason why we have used scare quotes in writing ''any' congruence', above. In any case, this last result can be easily remedied by extensions of **Cila**. Consider, for instance, the logic **Cilo** (the logic $C_1^+$ of da Costa, Béziau & Bueno's [57] —check again the subsection **3.10**). A non-trivial congruence can be defined within this logic by requiring, for any two given formulas, that they are not only provably equivalent, but are also both consistent. This can be put in terms of a single formula, by defining $A \approx B$ if $\vdash ((A \leftrightarrow B) \wedge (\circ A \wedge \circ B))$. One can then immediately prove that:

FACT 3.81 There is a non-trivial quotient algebra for **Cilo** (and already for **Cio**).
**Proof:** The above defined connective $\approx$ clearly sets up an equivalence relation. We have to show that it is in fact a congruence, so that given a schema $G(A)$ depending on $A$ as a component formula (and possibly on some other formulas as well), we have to show that $G(A) \approx G(B)$ whenever $A \approx B$, where $G(B)$ is obtained by replacing each occurrence of $A$ in $G(A)$ by $B$. Now, given this supposition that $A \approx B$, and recalling from THEOREM 3.67 that the consistency of any component of a complex formula, in **Cio**, is enough to guarantee the consistency of the complex formula itself, we may infer that $\vdash \circ G(A)$ and $\vdash \circ G(B)$. To check that $\vdash (G(A) \leftrightarrow G(B))$ just do a straightforward induction on the complexity of $G$. In the trivial case in which no other connectives or formulas intervene, but $A$, there is really nothing to prove. The case of conjunction, disjunction and implication is also immediate, from positive logic. For negation, just recall, as a consequence of FACT 3.17, that contraposition holds for provably consistent formulas, so that from $A \dashv\vdash B$ and both $\vdash \circ A$ and $\vdash \circ B$ one can infer $\neg A \dashv\vdash \neg B$. This concludes the proof (a similar semantical argument can already be found in [57], Theorem 3.21), and it is obvious that this congruence is non-trivial —we know for example from FACT 3.66 that $\circ\bot$ is a theorem of **Cila**, and thus of **Cilo**, so that all bottom particles will of course belong to the same equivalence class determined by $\approx$ over **Cilo**. In all other respects, except for this last particular example, the above proof is clearly valid not only in **Cilo** but also in **Cio** (that is, **Cilo** without the axiom (cl) that transforms this last **C**-system, **Cio**, into a **dC**-system).

Various other extensions of **Cila** having non-trivial quotient algebras have been proposed in the literature. In [84], for instance, Mortensen has proposed an infinite number of them, all situated of course somewhere in between **Cila** and classical logic. They were called $C_{n/(n+1)}$, for $n>0$ ($C_0$ is the name traditionally reserved for classical logic), and axiomatized by the addition to **Cila** of the following axioms, for each fixed $n>0$:

(M1n)    $\neg^{n-1}A \vdash \neg^{n+1}A$, where $\neg^n$, as usual, denotes $n$ iterations of $\neg$;

(M2n)    $\bigwedge_{i=1}^{n}(\neg^{i-1}A \leftrightarrow \neg^{i-1}B) \vdash \bigwedge_{i=1}^{n}(\neg^{i}(A\#C) \leftrightarrow \neg^{i}(B\#C)) \wedge \bigwedge_{i=1}^{n}(\neg^{i}(C\#A) \leftrightarrow \neg^{i}(C\#B))$, where # is any binary connective, and $\wedge$ abbreviates, as usual, a long conjunction.

In the section 4 of [84] the connective $\approx$ defined by letting $A \approx B$ hold whenever $\vdash \bigwedge_{i=0}^{n}(\neg^{i}A \leftrightarrow \neg^{i}B)$ is shown to constitute a non-trivial congruence, for each $n>0$. The reason for non-triviality is that, in general, each $C_{n/(n+1)}$ can be understood as providing us with $n+1$ 'negations': for any formula $A$ of this logic we have that $\neg^{m-1}A$ is congruent to $\neg^{m+1}A$ if, and only if, $m \geq n$, so that there are $n+2$ distinct equivalent classes (represented by $A$, $\neg A$, ..., $\neg^{n+1}A$) of the quotient algebra generated by $\approx$. Do any of these new **C**-systems coincide with any of the other above studied ones? Do they have any special interest in themselves (besides being equipped with a non-trivial congruence)?

How can one understand these more general algebras induced by more esoteric congruence relations, if they do not fit inside the 'classical' algebraization theory of Lindenbaum-Tarski? A neat and elegant solution to that can be found in the study of Blok & Pigozzi (cf. [25]), where a much more general theory of algebraization is developed, extending the work of other authors. Some terminology and definitions are needed to explain what is a *Blok-Pigozzi algebraization*. Fixing some logic **L**=<*For*, ⊩>, an **L**-*algebra* is any structure homomorphic to **L** (being *For* a structured set of formulas constructed over some set of connectives, the corresponding **L**-algebra will of course contain, for each connective, an operator of the same arity 'interpreting' it). An **L**-*matrix model* of an **L**-algebra **Alg** is any pair <**Alg**, **D**>, where **D** is a proper subset of the universe of **Alg**, of the so-called *distinguished elements*. Formally, let an *interpretation* of a set of formulas *For* be an assignment of terms of **Alg** to each element of *For* (an assignment which is usually defined over some primitive elements and then extended to the whole set of formulas by way of the interpretation of the building structural operators). The *semantic consequence relation* $\vDash_{\mathbf{M}}$ associated to an **L**-matrix model **M** is then defined, as usual, by setting $\Gamma \vDash_{\mathbf{M}} A$ whenever $A$ is assigned a distinguished element for every assignment of distinguished elements to all members of $\Gamma$. Matrices of finite many-valued logics are simple practical examples of *sound* and *complete* matrix models (that is, models such that $[(\Gamma \Vdash A) \Leftrightarrow (\Gamma \vDash_{\mathbf{M}} A)]$) that can be associated to some logics. In general, by a result of Wójcicki (see [111]) it is known that every structural logic can be characterized by sound and complete matrix models, in fact by $\kappa$-valued matrices, where $\kappa$ has at most the cardinality of the set of formulas of the logic.

Any pair of terms $\varphi$ and $\psi$ of the **L**-algebra will be said to constitute an *equation*, to be designated by writing $(\varphi \doteq \psi)$. Such equations are always schematic, as any usual mathematical equation, and their non-operational components are said to be its *variables*; we may accordingly write $\varphi(C) \doteq \psi(C)$ to designate an equation having $C$ as its single variable, and similarly for any number of variables. Now, what an interpretation does is exactly assigning values to these variables. One may then define an *equational consequence relation* induced by a class of **L**-algebras **KA**, to be denoted as $\vDash_{\mathbf{KA}}$, as follows: $[\Gamma \vDash_{\mathbf{KA}} (\varphi \doteq \psi)]$, where $\Gamma$ is a set of equations, whenever the equation $(\varphi \doteq \psi)$ is a semantic consequence of $\Gamma$ for every **L**-matrix model **M** of each **L**-algebra in **KA**, that is, when all those matrix models are such that $[\Gamma \vDash_{\mathbf{M}} (\varphi \doteq \psi)]$. The relation $\vDash_{\mathbf{KA}}$ is said to constitute an adequate *algebraic semantics* for a given logic **L** whenever there is a finite set of equations $\delta_i(C) \doteq \varepsilon_i(C)$, for $i<n$, such that: $(\Gamma \Vdash A) \Leftrightarrow [\{\delta_i(B) \doteq \varepsilon_i(B): \text{for all } i<n \text{ and all } B \in \Gamma\} \vDash_{\mathbf{KA}} (\delta_i(A) \doteq \varepsilon_i(A))]$. In this case, the equations $\delta_i(C) \doteq \varepsilon_i(C)$, for $i<n$, are called *defining equations* of **L**, and we shall write simply $\delta \doteq \varepsilon$ as an abbreviation of them. Finally, an algebraic semantics for a logic **L**, induced by a class of **L**-algebras **KA**, is said to be *equivalent* (or *congruential*) if there can be defined in **L** a finite set of connectives with two variables $\approx_j$, for $j<m$, such that, for every equation $\varphi \doteq \psi$, we have that $[\{(\varphi \approx_j \psi): \text{for all } j<m\} =\!\|\!\vDash_{\mathbf{KA}} \{\delta(\varphi \doteq \psi) \approx_j \varepsilon(\varphi \doteq \psi): \text{for all } j<m\}]$. This set of connectives $\approx_j$, for $j<m$, will be abbreviated simply as $\approx$ and called a *system of equivalence* (or *congruence*) connectives for **L** and **KA**. Now, a logic **L** is said to be (Blok-Pigozzi-)*algebraizable* if it has an equivalent algebraic semantics. Another way of stating this definition (in terms of the consequence relation of **L**) is by requiring, to call a logic **L** algebraizable, to have in hand a set of equations $\delta \doteq \varepsilon$ and a set of formulas $\approx$ such that: (i) $\approx$ constitutes an equivalence relation; (ii) $(A_1 \approx B_1)$, …, $(A_n \approx B_n) \Vdash \sigma(A_1, …, A_n) \approx \sigma(B_1, …, B_n)$, for each $n$-ary connective $\sigma$; and (iii) $A \dashv\!\Vdash \delta(A) \approx \varepsilon(A)$. It should by now be completely clear how this generalizes the idea of (proving (IpE) and) producing a congruence over a set of formulas.

Not all logics are algebraizable (even in this broader sense of Blok-Pigozzi). For example, most modal logics, and the system **E** of entailment are not algebraizable, though they do have non-congruential algebraizable semantics. As to the **C**-systems that we study in this paper, it has already been shown or mentioned some lines above (FACT 3.81 and below), that the logics **Cilo** and $C_{n/(n+1)}$, extensions of **Cila**, do have non-trivial congruences, being thus algebraizable in the sense of Blok-Pigozzi (though they are not algebraizable in the traditional sense of Lindenbaum-Tarski).

Also, as hinted in the last subsection, one can now prove that all the 8K three-valued maximal paraconsistent logics there presented are algebraizable (making use of and extending an argument by Lewin, Mikenberg & Schwarze, who have proved in [71] that the three-valued logic $\mathbf{P}^1$ is algebraizable):

**FACT 3.82** All the 8K three-valued logics from the last subsection are algebraizable.
**Proof:** Just consider any of the two connectives $\equiv$ or $\stackrel{\cdot}{\equiv}$ defined in the FACT 3.75, let $\delta(A)$ be defined as $((A{\to}A){\to}A)$ and $\varepsilon(A)$ be defined as $(A{\to}A)$, and check that the conditions (i)–(iii) defining an algebraizable logic two paragraphs above do hold.

It can also be shown, at this point, that some of our **C**-systems are not algebraizable. To such an intent, yet another characterization of algebraizable logics can come on handy. Let **L** be a logic and **M** be an **L**-matrix model. A *Leibniz operator* $\Lambda$ is a mapping from each arbitrary subset $S$ of **M** into the largest congruence $\approx$ of **M** compatible with $S$, where $\approx$ is *compatible* with $S$ if whenever we have that $\varphi \in S$ and $\varphi \approx \psi$ we also have that $\psi \in S$. It can be proved that a logic **L** is algebraizable if, and only if: (iv) $\Lambda$ is injective and order-preserving on the collection $CT(\mathbf{L})$ of all closed theories of **L**; (v) $\Lambda$ preserves unions of directed subsets of $CT(\mathbf{L})$, where a subset of $CT(\mathbf{L})$ is *directed* if there is a common upper limit to every finite collection of elements of $CT(\mathbf{L})$. (At this point, we had better direct the reader to [25] for details and proofs.) In any case, one might observe that a consequence of these last observations is that, for every logic **L**, the Leibniz operator produces an isomorphism between the lattice of filters of each **L**-matrix model **M** and the lattice of congruences of **M**. So, if such an operator is not an isomorphism, for some **L**-matrix model **M**, then the logic **L** is not algebraizable. This was the idea used by Lewin, Mikenberg & Schwarze in [72] (and that we extend here) to prove that:

**THEOREM 3.83** The logic **Cila** (that is, da Costa's $C_1$) is not algebraizable. The same holds even for the stronger logic **Cibaw** (see **Figure 3.1**, in the subsection **3.10**), or any weaker logics extended by **Cibaw**.
**Proof:** Consider the following set of truth-values, $\mathbf{V} = \{0, a, b, 1, u\}$, ordered as follows: $0 \le a$, $0 \le b$, $a \le 1$, $b \le 1$, $1 \le u$, and where $u$ and $1$ are the distinguished elements. Consider now the following matrices defined over them:

| $\wedge$ | $u$ | $1$ | $a$ | $b$ | $0$ |
|---|---|---|---|---|---|
| $u$ | $u$ | $1$ | $a$ | $b$ | $0$ |
| $1$ | $1$ | $1$ | $a$ | $b$ | $0$ |
| $a$ | $a$ | $a$ | $a$ | $0$ | $0$ |
| $b$ | $b$ | $b$ | $0$ | $b$ | $0$ |
| $0$ | $0$ | $0$ | $0$ | $0$ | $0$ |

| $\vee$ | $u$ | $1$ | $a$ | $b$ | $0$ |
|---|---|---|---|---|---|
| $u$ | $u$ | $u$ | $u$ | $u$ | $u$ |
| $1$ | $u$ | $1$ | $1$ | $1$ | $1$ |
| $a$ | $u$ | $1$ | $a$ | $1$ | $a$ |
| $b$ | $u$ | $1$ | $1$ | $b$ | $b$ |
| $0$ | $u$ | $1$ | $a$ | $b$ | $0$ |

| $\to$ | $u$ | $1$ | $a$ | $b$ | $0$ |
|---|---|---|---|---|---|
| $u$ | $u$ | $u$ | $a$ | $b$ | $0$ |
| $1$ | $u$ | $1$ | $a$ | $b$ | $0$ |
| $a$ | $u$ | $1$ | $1$ | $b$ | $b$ |
| $b$ | $u$ | $1$ | $a$ | $1$ | $a$ |
| $0$ | $u$ | $1$ | $1$ | $1$ | $1$ |

| | $\neg$ | $\circ$ |
|---|---|---|
| $u$ | $1$ | $0$ |
| $1$ | $0$ | $1$ |
| $a$ | $b$ | $1$ |
| $b$ | $a$ | $1$ |
| $0$ | $1$ | $1$ |

All axioms of **Cibaw** are validated by these matrices, as the reader can easily check. Now, it is also easy to check that there are no non-trivial congruences over **V**. Suppose for instance that $u \approx x$, for some $x \neq u$. In this case, as we know that $\neg\neg u \doteq 0$, and $\neg\neg x \doteq x$, then the condition (ii) above will give us $\neg\neg u \approx \neg\neg x$ from $u \approx x$, and so we conclude that $0 \approx x$, and thus $0 \approx u$. But, as we have observed before, there can be no congruence class containing both distinguished and non-distinguished values (in any case, this will violate condition (iii) above). We leave to the reader the easy exercise of showing, using the above connectives, that for any $x \approx y$, with $x \neq y$, one gets trapped at a similar predicament, namely, that of a distinguished value getting grouped with a non-distinguished one inside the same congruence class. Now, it is clear that $<A, \wedge, \vee>$ is a lattice, and that $\{0, a, 1, u\}$ and $\{0, b, 1, u\}$ are two filters over **V**. But there is just one congruence over **V** (which is of course the largest one compatible with both the filters just mentioned), and so the Leibniz operator cannot be an isomorphism. Once the logic **Cibaw** is, as a consequence, not algebraizable, FACT 3.79 informs us that none of its fragments can be algebraizable.

Evidently, THEOREM 3.80 may be proved as a corollary of this last result.

Now, even non-algebraizable logics can happen to be amenable to sensible algebraic investigation. Indeed, a class of *protoalgebraic logics* is characterized by the validity of condition (iv) above, one of the clauses on the characterization of Blok-Pigozzi algebraizability in terms of the Leibniz operator. This class includes all normal modal logics and most non-normal ones, but there still are some other logics which are not protoalgebraizable: an example is **IPL\***, the implicationless fragment of intuitionistic logic, is neither algebraizable nor protoalgebraizable (cf. [25], chapter 5). Which of our non-algebraizable **C**-systems are protoalgebraizable, and which not (if any)? We shall leave this question open at this stage. It is interesting to notice, at any rate, that some sort of algebraic counterparts to some of these non-algebraizable **C**-systems have been proposed and studied, for instance, in Carnielli & de Alcantara's [37] and Seoane & de Alcantara's [102], where a variety of 'da Costa algebras' for the logic **Cila** has been introduced and studied, and a Stone-like representation theorem was proved, to the effect that every da Costa algebra is isomorphic to a 'paraconsistent algebra of sets'. It would be interesting now not only to extend that approach to other **C**-systems, but also to check how it fits inside this more general picture given by (Blok-Pigozzi-)algebraizable and protoalgebraizable logics.

An interesting application of the above mentioned algebraic tools is the following. Consider again, for example, the FACT 2.10(ii), where strong negations were shown to be 'definable' from bottom particles. Now, it is completely clear how this definition can be stated in practice if, for instance, a suitable implication is available inside of a compact logic (this is the case in all our examples, but needs not to be). This illustrates in fact how intuitionistic negation may be defined from a bottom particle and intuitionistic implication. But is that *definition*, in the general case, an *implicit* or an *explicit* one? For example, in positive classical logic (plus bottom and top) the theory containing both the formulas $((A \wedge B) \leftrightarrow \bot)$ and $((A \vee B) \leftrightarrow \top)$ implicitly defines the formula $B$ (as the 'classical negation of $A$'). Are all implicit definitions also explicit ones? Or do we have, in some cases, to explicitly add some more structure to a logic to make explicit definitions expressable even when implicit ones are available? If, whenever a logic can implicitly define something, it can also explicitly de-fine it, then the logic is said to have the *Beth definability property*, (BDP). Now, consider any class of **L**-algebras **KA**, for some logic **L**, and pick up a set **HA** of homomorphisms between any two of these **L**-algebras. A homomorphism $f: \mathbf{Alg}_1 \to \mathbf{Alg}_2$ in **HA** is said to be an *epi* if every pair of homomorphisms $g, h: \mathbf{Alg}_2 \to \mathbf{Alg}_3$ in **HA** is such that $g \circ f = h \circ f$ only if $g = h$. Evidently, all surjective homomorphisms are epis; if the converse also holds, that is, if all epis are surjective, we say that **KA** has the property (ES). Now, by a result of I. Németi (cf. [4]), an algebraizable logic has (BDP) if, and only if, its class of algebras has (ES). This is a very interesting result, and constitutes, in fact, just one example of how algebraic approaches can help us to solve real logical problems, in this case the problem of definability. Extensions of such results to wider classes of algebraic structures associated to (wider classes of) logics are clearly desirable.

## 4    FUTUROLOGY OF **C**-SYSTEMS

> When you encounter difficulties and contradictions, do not try to break them, but bend them with gentleness and time.
> —Saint Francis de Sales.

This is *not* the end. The next and natural small step for a paper, giant leap for paraconsistency, is providing reasonable interpretations for **C**-systems. This is the theme of our [42], where semantics for **C**-systems are presented and surveyed, ranging from the already traditional *bivaluations* to the more recently proposed *possible-translations semantics*, traversing on the way a few connections to many-valued semantics (a theme that already intromitted in our subsection **3.11**), and to modal semantics. A quite diverse approach to paraconsistent logics (in general) from the semantical point of view is also soon to be found in Priest's [91] (the remarkable possible-translations semantics, according to which some complex logics are to be understood in terms of *combinations* of simpler ones, will nevertheless not be found there —see instead [76], [36], [39], and, of course, [42]—, in addition, its section on many-valued logics is unfortunately too poor to give a reasonably good idea on the topic).

In this last section of the present paper we want to point out some interesting open problems and research directions connected to what we have herein presented. For example, in the section **2** we have extensively investigated the abstract foundations of paraconsistent logic, and the possibility and interest of defining the so-called *logical principles* at a purely logical level. There is still a lot of stirring open space to work here, and we will feel happy to have stimulated the reader to try their own hand at the relations between all the alternative formulations of (PPS), that is, all the different forms of explosion (or, if they prefer, the various forms of *reductio*, as in the subsection **3.2** —recall that the *reductio* and the Pseudo-Scotus are not always equivalent, for instance, if you think about intuitionistic logics). Think about it: are there any *interesting* logics, in our sense, disrespecting the Pseudo-Scotus, while respecting *ex falso* or the Supplementing Principle of Explosion, but still disrespecting the Gentle Principle of Explosion as well? In other words, are there interesting paraconsistent logics having either bottom particles or strong negations which do not constitute **LFIs**?

Recall that our approach contributed a novel notion of *consistency*. This is the picture again: There are consistent and inconsistent logics. The inconsistent ones may be either paraconsistent or trivial, but not both. The paraconsistent ones may be either dialectical or not. The consistent logics are explosive and non-trivial. The paraconsistent logics are non-explosive, and the dialectical paraconsistent ones are contradictory as well. The trivial logics (or trivial logic, if you fix some language) are explosive and contradictory (if the underlying logic has a negation symbol). Negationless logics are trivial if and only if inconsistent. Let us say that a theory *has models* only if these are non-trivial (they do not assign distinguished values to all formulas). So, the theories of a consistent logic have models if and only if they are non-contradictory. Paraconsistent logics may have models for some of its contradictory theories, and in the dialectical case all models of all theories are contradictory. Trivial theories (of trivial logics, or of any other logics) have no models. The consistency of each formula $A$ of a logic $\mathbf{L}$ is defined exactly as what else one must say about $A$ in order to make it explosive, that is, as what one should add to an $A$-contradictory theory in order to make it trivial. If the answer is 'nothing', then $A$ is already consistent in $\mathbf{L}$ (whether the theories that derive this formula are contradictory or not). So, consistent logics are, quite naturally, those logics that have only consistent formulas. The above study sharply distinguishes the notions of non-contradictoriness and of consistency, and the model-theoretic impact of this should obviously be better appraised!

We now also have a precise definition of a large and fascinating class of paraconsistent logics, the *logics of formal inconsistency*, **LFI**s, and an important subclass of that, the **C**-*systems*. This is important to stress: according to our proposal it should be no more the case that the **C**-systems will be identified simply with the calculi $C_n$ of da Costa, or with some other logics which just happen to be axiomatized in a more or less similar way. A general idea was put forward to be explored, namely that of being able to express *consistency* inside of our paraconsistent logics, and this helped collecting inside one big single class logics as diverse as the $C_n$ and $\mathbf{P}^1$, or even $\mathbf{J}_3$ (now rephrased as **LFI1**), whose close kinship to the $C_n$ seemed to have passed unsuspected until very recently (recall the subsections **2.4**, and end of **3.11**). This is, we may suggest, a fascinating challenge that we propose to our readers: To show that many other logics in the literature on paraconsistent logics can be characterized as **C**-systems, or, in general, as **LFI**s. This exercise has been explicitly put forward in the subsection **2.6**, but even previous to that, in the end of the subsection **2.4**, we have already hinted to the fact that other logics, such as Jaśkowski's **D2**, a discussive paraconsistent logic with motivations and technical features completely different from the ones that we study here, could be recast as an **LFI** (based on the modal logic $S5$) —more precisely, it can be recast as an **LFI** if only it is presented having some necessity operator, $\square$, among its primitive or definable connectives (see note 11). Another recent example of that is the paraconsistent logic **Z** (are we perhaps running out of names?) proposed by Béziau in [24], in which a paraconsistent negation $\neg$ is defined from a primitive classical negation $\sim$ and a possibility operator $\lozenge$, by setting $\neg A \overset{\text{def}}{=} \lozenge \sim A$.[27] Again, it is easy to see that **Z** can also be seen as an **LFI** (in fact, a

----

[27] By the way, exactly the same logic was proposed by Batens in [16] under the name **A**, and appears on another of Béziau's paper, [23], section 2.8, under the appelation of 'Molière's logic'. Strangely enough, after longly attacking, in the section 2.5 of [23], those logics that he calls 'paraconsistent atomical logics', that is, those logics in which 'only atomic formulas have a paraconsistent behavior' (being

**C**-system based on *S*5), in which the consistency of a formula *A* is expressed by the formula ($\Box A \lor \sim A$). Which other paraconsistent logics constitute **C**-systems, or **LFI**s, and which *not*? Inverting the question, are there good reasons for one trying to *avoid* **LFI**s, that is, can the investigation of non-**LFI**s have good technical or philosophical justifications? And how would the **C**-systems based on intuitionistic logic (**I**-systems?) or on relevance logic (**R**-systems?) look like, and which interesting properties would they have? How would this improve our map and understanding of **C**-land? In general, how could one use the very idea of a **C**-system to build up some new interesting paraconsistent logics, what advantages would they bring and particular technical tools could they contribute to the general inquiry about **LFI**s? The point to insist here is on the remarkable *unification* of aims and techniques that **LFI**s can seemingly produce in the paraconsistency terrain!

Another related interesting route is the one of *upgrading* any given paraconsistent logic in order to turn it into an **LFI**. This is exactly what is done by the logic **LFI1** (or **CLuNs**, or **J**₃) over the logics *Pac* and *LP* (see subsection **2.4**), for which the gain in expressive power should already be obvious to the reader. Now, consider, for instance, the three-valued *closed set logic* studied by Mortensen and collaborators in [85], whose matrices of conjunction and of disjunction coincide with those of **LFI1**, and whose matrix of negation coincides with that of **P**¹, having, again, 0 as the only non-distinguished value. Now, it is easy to see that the addition of appropriate matrices of implication and of a consistency operator will turn the upgraded closed-set logic into one of our 8K three-valued maximal paraconsistent logics, discussed in the subsection **3.11** above. The motivation for such closed set logic is also to be found among some of the most striking features of the 'Brazilian approach' to paraconsistency, namely, the idea of studying paraconsistent logics which are in a sense *dual* to other *broadly intuitionistic* (also called *paracomplete*) logics. We have slightly touched on this issue at a few points above —see, for instance, subsections **3.1** and **3.2**— as this has been one of the preferred justifications used by da Costa, among other authors, for the constitution of many of his paraconsistent logics. Indeed, if classical logic is not rarely held by some authors as the 'logic of sets', particularly because of its Boolean algebraic counterpart, Heyting's Intuitionistic Logic is very naturally held, in a topological setting, as the 'logic of open sets'. The very same dualizing intuition that we have just mentioned can then lead one to study the 'logic of closed sets' as a very natural paraconsistent logic. This is done in [85], where this investigation is also lifted to the categorial space —again, if intuitionistic logic is very naturally thought of as the logic of a *topos*, the closed set logic can be thought of as the underlying logic of a categorial structure called *complemented topos*. The upgrade of

---

thus controllably explosive with respect to every complex, or 'molecular' formula), a class of logics that seems to comprise not many logics up to this moment (the logic **P**¹ —recall THEOREM 3.15(ii) and THEOREM 3.69— being among them), Béziau presents the above mentioned logic **Z** as the logic enjoying 'the best paraconsistent negation' around. But even if one concedes an enlargement of this definition of atomical logics in order to comprise all paraconsistent logics which are only 'paraconsistent up to some level of complexity', this author would still have to deal with the *bourgeois* fact that the negation of his preferred logic **Z**, exactly as what occurred with **P**¹'s negation, is such that $\{p, \neg p\}$ is not trivial, for atomic *p*, while $\{\neg p, \neg\neg p\}$ *is* trivial: why, in this largely analogous case, would that same phenomenon be 'philosophically justifiable', and not be reducible just to some more bits of 'formal nonsense', as he puts it? One interesting such a justification, we suggest, may be found in terms of the dualization obtained by logics such as Mortensen's closed set logic, also mentioned in the present section.

the closed set logic into an **LFI** may thus set up some interesting new space for the study of topological and categorial interpretations of the notion of consistency. This proposed duality has also often been pushed, in the literature, in the contrary direction, namely, into the study of paracomplete logics which are dual to some given paraconsistent logics. Some samples of this can be found, for instance, in our papers [39] (where a logic called $D_{min}$ is presented as dual to $C_{min}$, mentioned above in the subsection **3.1**), and [43] (where logics dual to slightly stronger versions of the calculi $C_n$ are studied), as well as in Marcos's [78] (where 1K three-valued maximal paracomplete logics are presented, in addition to the 8K paraconsistent ones above mentioned). A more thorough study of the **LFU**s, the *Logics of Formal Undecidability* (or *Logics of Incomplete Information*), following the standards of the investigation set up by the present paper rests yet to be done.

We of course do not, and cannot, claim to have included and studied above *all* 'interesting' **C**-systems based on classical logic. We cannot but offer here a very partial medium-altitude mapping of the region, but we strongly encourage the reader to help us expand our horizons. So, if, guided mostly by technical reasons, we have started our study from the basic logic **bC** (subsection **3.2**), constructing all the remaining **C**-systems as extensions of **bC**, this should *not* be a impediment for the reader to study still weaker and *more basic* logics, such as **mbC**, the logic axiomatized by deleting (Min11): ($\neg\neg A \to A$) from the axiomatization of **bC**, and their extensions, **mCi** (a logic studied under the name *Ci* in [17]) and so on. So, if **bC** was presented as a quite natural conservative extension of the logic $C_{min}$ ([39]), **mbC** can similarly be presented as an extension of the logic *PI* ([10]). Just keep in mind that starting your study from **mbC**, instead of **bC**, and avoiding the axiom (Min11), you will be allowing for the existence of a few more in principle uninteresting partially explosive paraconsistent extensions of your logics (THEOREM 3.13), and you may also lose a series of other results, as for instance 3.14(iii), half of 3.17 and of 3.26, and also 3.20, 3.36, 3.41, 3.51(iii), 3.56, 3.57, 3.66, as well as some derived results and comments. Notice that we do not say that the loss of some of these results cannot be positive, but some symmetry certainly seems to be lost if the logic **mbC** happens to be extendable, for instance, in such a way as to validate the schematic rule $[(A \to B) \vdash (\neg B \to \neg A)]$, though it can in no way be extended so as to validate the similar rule $[(\neg A \to B) \vdash (\neg B \to A)]$.

What effects can the **LFI**s have on the study of some general mathematical questions, such as *incompleteness results* in Arithmetic? Indeed, recall that Gödel's incompleteness theorems are based on the identification of 'consistency' and 'non-contradictoriness' —what then if we start from our present more general notion of consistency (see (D19), subsection **2.4**)? And what if we integrate these logics with such modal logics as the *logic of provability*? (See [27], where consistency is also intended as a kind of dual notion to provability —if you cannot prove the negation of a formula, it is consistent with what you can prove—, and in which the necessary environment for the study of Gödel's theorems is provided.) What effects could that have (if at all) on the investigation of (set-theoretical) paradoxes? In fact, the analogy of the logics of formal inconsistency with the logics of provability is rather striking and worthy of being further explored; connections with other powerful internalizing-metatheoretical-notions logics, such as *hybrid logics*, and *labelled deductive systems* in general, are also to be expected.

We are *not* trying to escape from da Costa's requisite **dC[iii]** (subsection **3.8**), according to which extensions of our logics to higher-order logics ought to be available. But it still seems that most interesting problems related to *paraconsistency* appear already at the propositional level! Moreover, more or less automatic processes to *first-ordify* some given propositional paraconsistent logics can be devised, by the use of combination techniques such as *fibring*, if only we choose the right abstraction level to express our logics (see [34], where the logic $C_1$ is given a first-order version —coinciding with the one it had originally received, in [49] or [50], but richer in expressibility power— by the use of the notion of *non-truth-functional rooms*). One interesting thing about first-order paraconsistent logics is that they might allow for inconsistencies at the level of its objects, opening a new panorama for ontological investigations. Another interesting thing that we can conceive about first-order versions of paraconsistent logics in general, and especially of first-order **LFI**s, and that seems to have been completely neglected in the literature up to this point, is the investigation of consistent yet ω-*inconsistent* structures, or theories (also related to Gödel's theorems). Again, a point about that is scored by the logics of provability, but studious of paraconsistency should definitely have something to say about this.

Let us further mention some more palpable specific points to which we have already drawn attention here and there, and on which more research is still to be done. For instance, is the logic **bC** (subsections **3.2** to **3.4**) *controllably explosive*? Recall from FACT 3.32 that in the case of its extension **Ci**, the formulas causing controllable explosion coincide with the consistent theorems, and that, as it happens, **bC** does not have consistent theorems (THEOREM 3.10). For **bC**, however, only one side becomes immediate: consistent theorems cause controllable explosion... Another question: Does this logic **bC** have an intuitive adequate modal interpretation? And are there also extensions of **Ci** in which (IpE) holds (see the end of the subsection **3.7**)? What about investigating some other extensions of **bC** which do not extend **Ci** as well, such as **bCe**, obtained by the direct addition to **bC** of the axiom (ce): $A \vdash \neg\neg A$, studied in subsection **3.9**? Remember that we have extended **bC** to **Ci**, in the subsection **3.5**, arguing that all paraconsistent logics in the literature *do* identify inconsistencies and contradictions, and then all the other logics that we have studied after that in fact extended **Ci**. But the logic **bCe** also seems interesting in its own right, being able to accept some simple extensions that can express dual inconsistency, differently of what had occurred in the subsection **3.4** with some other extensions of **bC**, and paralleling a result obtainable for **Ci** (subsection **3.6**). **bCe** would also presumably constitute a step further in the direction of obtaining full (IpE), but in any case the general search for interesting **C**-systems extending **bC** and not extending **Ci** can already be funny enough. As to other ways of fixing the non-duality of the consistency and inconsistency operators in **bC**, alternative to extending this logic into **Ci**, suggestions have been made, for instance, for the addition to **bC** of schematic rules such as $\circ A$, $\neg \circ A \vdash B$ or else $\vdash \circ\circ A$ (perhaps having FACT 3.32 in mind, and trying to carry it forward into **bC**). None of these will work, however, as it is easy to see if one just considers again the matrices in THEOREM 3.16, noticing that they also validate the two last rules, while still not validating schemas such as (ci1) or (ci2). But there may quite well exist some other way out of this quagmire (perhaps the reader will find it).

Now, this is a quickie: can you find any (grammatical) conservative translation from **eCPL** (the extended classical propositional logic obtained by the addition to classical logic of the then innocuous consistency operator) into **bC** (recall the subsection **3.7**), as the one we had for **Ci** (COROLLARY 3.47)? And what happens when one considers, in the construction of **dC**-systems, the addition of more general rules such as (cg) or (RG) (see the subsection **3.8**), which implement more inclusive definitions for the consistency connective? Do we obtain interesting logics from that, fixing some asymmetries observed on the calculi $C_n$ and its relatives? What effects does this move have on the semantical counterparts of these logics? Moving yet farther, we may ask about the logic $C_{Lim}$, the deductive limit to the hierarchy $C_n$ (see subsection **3.10**), whether it can be proved to be *not* finitely gently explosive. For if it is not finitely gently explosive, given that it is gently explosive *lato senso*, then it cannot be *compact*. Or are we rather obligated to abandon *strong* completeness, in this case? Would there be other interesting **LFI**s in which consistency is not *finitely* expressible?

Recall the subsection **3.11**, where da Costa's requisite **dC[iv]** (subsection **3.8**) is taken very seriously, and we search for *maximal* paraconsistent fragments of classical logic, that is, logics having 'most rules and schemas of classical logic', and a class of 8K three-valued logics is presented as a solution to this. Now, are there other (full) solutions to this 'problem of da Costa'? And are there other interpretations, besides maximality, of da Costa's requisite **dC[iv]** leading to yet some other solutions to that problem? Are there non-many-valued (monotonic) solutions to that problem, or perhaps some other $n$-valued ones, for $n>3$? And, this is an important first step and probably an easier problem to solve: are there other interesting **C**-systems based on classical logic which are *not* extended by any of the above 8K three-valued logics? We must leave these questions open at this moment. It is interesting to notice, at any rate, that this problem has already been directly addressed here and there, in the literature. Besides [78], from which we drew the results in the subsection **3.12**, one could also recall, for instance, the *adaptive* programme for the confection of paraconsistent logics aiming to represent (non-monotonically) the dynamics of scientific reasoning and of argumentation (see [15]). Roughly speaking, the basic idea behind adaptive logics is that of working in between two boundary logics, often classical logic constituting one of them and a paraconsistent logic constituting the other one, so that consistency is pressuposed by default and we try to keep on reasoning inside of classical logic up to the point in which an 'abnormality' (an inconsistency?) pops out, a situation in which we had better descend to the level of the complementary paraconsistent logic, and go on reasoning over there. Indeed, the ancestral motivations of this programme seem to have been, as it is reasonable to conceive, yet another attempt to maintain the most of classical logic as possible (see [10]), so that paraconsistency is only needed at limit cases, while (most of) classical logic is, in principle, maintained in 'normal' situations. The difficulty here, as far as da Costa's requisite **dC[iv]** is concerned, seems to be *measuring* how much some different adaptive logics will respond to that same requirement of closeness to classical logic (if there is any measurable difference at all). As it is appealing to think of adaptive logics as situations in which two logics are combined in order to produce a third one, it seems also interesting to investigate if possible-translations semantics can, after all, be applied to such an environment as well, or at least stretch the analogies there as far as we can.

Some open questions can also be drawn from the subsection **3.12**. Notice, for instance, that all of Mortensen's axioms, (M1$n$) and (M2$n$), for every $n>0$, are validated by the matrices of the three-valued paraconsistent logic $\mathbf{P}^2$ (recall THEO-REM 3.69), as that author himself have observed, and the question was left open, in [84], (and it remains still so, as far as we know) whether the logic $C_{1/2}$, the stronger of the logics $C_{n/(n+1)}$, would in fact coincide with the three-valued logic $\mathbf{P}^2$. Once $\mathbf{P}^2$ is known to be a maximal paraconsistent logic (cf. THEOREM 3.72), to show this coincidence would amount to showing that all axioms of $\mathbf{P}^2$ (and especially axioms (ca1)–(ca3), in the subsection **3.10**, specifying the consistent behavior of binary connectives) are provable from the axioms of $C_{1/2}$. Alternatively, using the fundamental FACT 3.32, one could try to show that $C_{1/2}$ is controllably explosive (or not, if what one wants is to disprove the conjectured coincidence) in contact with any formula involving binary connectives. In one way or another, it is quite interesting to note that the 8K maximal three-valued logics in subsection **3.11** show that there are several other logics different from $\mathbf{P}^1$ and from $\mathbf{P}^2$ that are next to the classical propositional logic 'in the same kind of way' (half of them extending $C_1$, as we point out in FACT 3.73), a problem that was left open in the closing paragraph of the above mentioned paper.

Now, what about extending our investigations on the *algebraizability* of **C**-systems (again, see the subsection **3.12**)? Can these algebras solve yet some other categories of logical problems? Again, notice that the problem of finding extensions of our **C**-systems which are algebraizable in the 'classical sense' was also left open (though the plausibility of the existence of such extensions was hinted) in the end of subsection **3.7**. To be precise, what was open was the existence of such extensions as fragments of some version of classical logic —the reader will have seen in the present section, however, that modal logics such as **Z** do extend **bC** and have no problem on what concerns (IpE) (given that $S5$ is algebraizable). And what to say of extending our general approach on section **2** to other 'kinds of logics' (that is, varying the *logic structure* that we defined there) so as to include other kinds of consequence relations, such as (multiple-conclusion) non-monotonic ones (see [5] and [12])? Under this new light shed by **C**-systems and **LFI**s, it is also interesting to see how one can also move on to improve our present (rather poor) proof-theoretical approach. Indeed, Hilbert-style systems, such as the ones we present here, often require too much ingenuity to be applied, leaving intuition or mechanization of proofs far behind. For some **dC**-systems we know that *sequent* systems have already been proposed (see for instance [97] and [17]), as well as *natural deduction* systems (see [58]), and *tableau* systems (see [38]). The *really* interesting cases, however, seem to be those of **C**-systems that are *not* **dC**-systems, so that the consistency connective is, in a sense, 'ineliminable'! A first step towards such a general treatment of **C**-systems in terms of tableaux has already been offered by us in [41], where the logics **bC**, **Ci** and **LFI1** were all endowed with sound and complete tableau formulations.

There is so much yet to be done!

Bucharest, Rumania, July 2000. Edited by Mircea Dumitru, under the auspices of the New Europe College of Bucharest.